\documentclass[12pt]{amsart}
\usepackage{amssymb}
\usepackage{amsmath}
\usepackage{bm}
\usepackage{amsthm}
\usepackage{amsbsy}
\usepackage{psfrag}
\usepackage{graphics}
\usepackage{graphicx}
\usepackage{pgfplots}
\usepackage{relsize}
\usepackage{color}
\usepackage{graphicx}
\usepackage{caption}
\usepackage{subcaption}

\theoremstyle{plain}
\newtheorem{theorem}{Theorem}[section]
\newtheorem{proposition}[theorem]{Proposition}

\newtheorem{lemma}[theorem]{Lemma}

\newtheorem{example}[theorem]{Example}
\newcommand{\norm}[1]{\left\Vert #1 \right\Vert}

\theoremstyle{remark}
\newtheorem{remark}[theorem]{Remark}

\begin{document}
\allowdisplaybreaks[4]
\numberwithin{figure}{section}
\numberwithin{table}{section}
 \numberwithin{equation}{section}
%
\title[Optimal control of the Stokes equations]
 {An Error analysis of Discontinuous Finite Element Methods for the Optimal Control problems governed by Stokes equation}

 \author[A. K. Dond]{Asha K Dond}
  \address{Department of Mathematics, Indian Institute of Science, Bangalore - 560012, India}
  \email{ashadond@iisc.ac.in}

\author[T. Gudi]{Thirupathi Gudi}
 \address{Department of Mathematics, Indian Institute of Science, Bangalore - 560012, India}
 \email{gudi@iisc.ac.in}

 \author[R. C. Sau]{Ramesh Ch. Sau}

\address{Department of Mathematics, Indian Institute of Science, Bangalore - 560012, India}

\email{rameshsau@iisc.ac.in}

\date{}

\begin{abstract}
In this paper, an abstract framework for the error analysis of  discontinuous finite element method is developed for the distributed and Neumann boundary control problems governed by the stationary Stokes equation with control constraints.
 {\it A~priori} error estimates of optimal order are derived for velocity and pressure in the energy norm and the $L^2$-norm, respectively. Moreover, a reliable and efficient {\it a~posteriori} error estimator is derived. The results are applicable to a variety of problems just under the minimal regularity possessed by the well-posedness of the problem. In particular, we consider the abstract results with suitable stable pairs of velocity and pressure spaces like  as the lowest-order Crouzeix-Raviart finite element and piecewise constant spaces, piecewise  linear and constant  finite element spaces. The theoretical results are illustrated by the numerical experiments.
\end{abstract}


\maketitle
\def \R{{{\Bbb R}}}
\def \P{{{\rm {\!}\cal P}}}
\def \Z{{{\rm {\!}\cal Z}}}
\def \Q{{{\rm {\!}\cal Q}}}
\allowdisplaybreaks
\def\d{\displaystyle}
\def\R{\mathbb{R}}
\def\cA{\mathcal{A}}
\def\dx{{\rm~dx}}
\def\ds{{\rm~ds}}
\def\u{\mathbf{u}}
\def\V{\mathbf{V}}
\def\w{\mathbf{w}}
\def\J{\mathbf{J}}
\def\f{\mathbf{f}}
\def\g{\mathbf{g}}
\def\div{{\rm div~}}
\def\v{\mathbf{v}}
\def\z{\mathbf{z}}
\def\b{\mathbf{b}}
\def\cB{\mathcal{B}}
\def\p{\partial}
\def\O{\Omega}
\def\s{\sigma}
\def\bbP{\mathbb{P}}
\def\vb{{\vec b}}
\def\bn{{\bf n}}
\def\cV{\mathcal{V}}
\def\V{\mathbf{V}}
\def\cR{\mathcal{R}}
\def\cM{\mathcal{M}}
\def\cT{\mathcal{T}}
\def\cE{\mathcal{E}}
\def\bbE{\mathbb{E}}
\def\ssT{{\scriptscriptstyle T}}
\def\HT{{H^2(\O,\cT_h)}}
\def\mean#1{\left\{\hskip -5pt\left\{#1\right\}\hskip -5pt\right\}}
\def\jump#1{\left[\hskip -3.5pt\left[#1\right]\hskip -3.5pt\right]}
\def\smean#1{\{\hskip -3pt\{#1\}\hskip -3pt\}}
\def\sjump#1{[\hskip -1.5pt[#1]\hskip -1.5pt]}
\def\jumptwo{\jump{\frac{\p^2 u_h}{\p n^2}}}

\section{Introduction}\label{sec:Intro}
We consider the following distributed control and Neumann boundary control problems governed by Stokes equations	
$$\text{min} ~J(\mathbf{v},\mathbf{x})= \frac{1}{2}\norm{\mathbf{v}-\mathbf{u}_d}_{0,\Omega}^2 + \frac{\lambda}{2}\norm{\mathbf{x}}_{A}^2,$$
 where $\norm{.}_{A}=\norm{.}_{0,\Omega}$  for distributed control problem  and
 $\norm{.}_{A}=\norm{.}_{0,\partial\Omega}$ for Neumann boundary control problem,
subject to,

\vspace{0.2cm}
\noindent\begin{minipage}{.5\linewidth}
	For distributed control problem	
	\begin{align*}
	-\Delta \mathbf{v}+\nabla{q}&=\mathbf{x}+\mathbf{f}\;\; \text{in}\;\Omega, \\
	\nabla\cdot{\mathbf{v}}&=0\;\qquad \text{in}\; \Omega,\\
	\bf{ v}&={\bf 0}\;\qquad\text{on}\;\partial\Omega,\\
	\mathbf{y}_a\leq \mathbf{x}(x)&\leq \mathbf{y}_b\;\; \text{a.e.}\;x\in \Omega.
		\end{align*}
\end{minipage} %
\begin{minipage}{.5\linewidth}
	For Neumann boundary control problem
	\begin{align*}
	-\Delta \mathbf{v}+\nabla{q}&=\mathbf{f}\;\; \text{in}\;\Omega,\\
	\nabla\cdot{\mathbf{v}}&=0\;\; \text{in}\;\Omega,\\
	{\partial \mathbf{v}}/{\partial \mathbf{n}}-q \mathbf{n}&=\mathbf{x}\;\; \text{on} \;\;\partial\Omega,\\
	\mathbf{y}_a\leq \mathbf{x}(x)&\leq \mathbf{y}_b \;\;\text{a.e.}\;x\in \partial\Omega.
	\end{align*}
\end{minipage}

\noindent
This paper investigates the discretization of the above systems  based on a finite element approximation of the state and the control variable and also develop an abstract framework for the error analysis of the above problem under minimal regularity. The discussion of discretizations of optimal control problems governed by partial differential equations started with papers of Falk \cite{Falk:1973}, Gevici \cite{geveci:1979}. Subsequently, there are many significant contributions to this field. It is difficult to list all the results in this introduction; we
refer to some of the articles and references therein for the development of numerical
methods and their error analysis. Refer to the monograph \cite{Troltzzsch:2005:book} for the theory of
optimal control problems and  the development of numerical methods. The primal-dual active set
algorithm has been developed in \cite{Kunish 2003}, and also it has been discussed in the context of the optimal
control problems. Apart from this, we refer
to \cite{Meyer & Rosch 2004} for a super-convergence result using a post-processed control
for constrained control problems. A variational discretization method has been
introduced in \cite{Hinze 2005} to derive optimal error estimates by exploiting the
relationship between the control and the adjoint state. For the numerical
approximation of Neumann boundary control problem with graded mesh
refinement refer to \cite{Apel & Rosche 2012} and for the numerical treatment of the Dirichlet
boundary control problems refer to \cite{Casas 2006,Hinze 2009,Swobodny 1991,May & Vexler 2013,Steinbach et. al. 2015} and references
therein.
On the other hand, while
the adaptive finite element methods based on \textit{a posteriori} error estimators
have grown in popularity, the study of \textit{a posteriori} error analysis for optimal
control problems has also gained much interest in the recent years.
In particular, the control in control constrained problem can exhibit
kinks and hence lacks smoothness. In this context, adaptive finite element
methods would be useful to enhance accuracy. An \textit{a posteriori} error analysis
of a conforming finite element method for control constrained problems has been
derived in \cite{Kiewag et al 2008}. Recently, a general framework for \textit{a priori} and \textit{a posteriori} energy norm error analysis for
Neumann and distributed control problems by discontinuous Galerkin discretization
can be found in \cite{S.C et al.:2015:IP} for scalar problems. The results therein are obtained by
the help of appropriate auxiliary problems. Local error analysis of discontinuous Galerkin methods for the
distributed control problem for the advection-diffusion equation has been studied in \cite{Leykekhman 2012}.

R\"osch and Vexler have applied the post-processing technique to a linear-quadratic optimal control problem governed by the Stokes equations \cite{Rosch & vexler 2006}. They have proved second order convergence under the assumption that the velocity field admits full regularity, which means it is contained in $H^2(\Omega)\cap W^{1,\infty}(\Omega).$ Nicaise and Sirch \cite{Nicaise & Sirch 2012} have extended the results of R\"osch and Vexler \cite{Rosch & vexler 2006} and Apel et al. \cite{Apel  & Sirch Winkler 2008, Apel & Winkler 2009}  to the conforming and nonconforming finite element methods for the optimal control of the Stokes equations under weaker regularity assumptions. This
means, they did not assume that the velocity field is contained in $H^2(\Omega)\cap W^{1,\infty}(\Omega),$ but only in some weighted space $H^2_{\omega}(\Omega)$. The analysis in 
\cite{Nicaise & Sirch 2012} is focused on the super-convergence result with regularity of the solution in some weighted Sobolev spaces. Our aim in this article is to derive the best approximation result in the energy norms under the weak regularity of the solution obtained through the weak formulation. It is natural to expect that the conforming methods exhibit this best approximation properties but it is not immediate for the nonconforming and discontinuous finite element methods. The nonconforming methods and discontinuous Galerkin methods are particularly attractive for the Stokes problem as they provide discrete inf-sup condition easily as compared to the conforming methods.

In this article, we consider a general optimality system of both the distributed and Neumann boundary control problem governed by Stokes equation.
We develop an abstract framework
for both \textit{a priori} and \textit{a~posteriori} error analysis of the general optimality system under some abstract assumptions.
We introduce two auxiliary problems: one dealing with an elliptic projection in \textit{a~priori} analysis and
the other is based on a reconstruction in \textit{a~posteriori} error analysis. Subsequently,
Theorem \ref{thm2.2} and Theorem \ref{thm2.6}  are proved, which play an  essential role in
the analysis. In particular, we consider the abstract results with  the lowest-order Crouzeix-Raviart finite element and piecewise constant spaces stable pair for velocity and pressure approximations respectively, and also  discontinuous Galerkin formulation with piecewise linear and constant  finite element spaces.
The outcome of the
result is the best approximation result for the numerical method.
Furthermore, we derive the optimal order of convergence for control, state, and adjoint state variables.
This framework for the error analysis of finite element methods for control
problems have been presented under limited regularity assumptions. It is worth noting that the standard error analysis of DG methods
require additional regularity which does not exist in several cases, for
example in mixed boundary value problems or simply supported plates,
example, see the discussions in \cite{Gudi 2010}. Therefore, the error analysis of
DG methods has to be treated carefully.
Here, the best approximation error estimates are derived under the minimum regularity on the state and the adjoint state variables for DG methods.
 Moreover, \textit{a posteriori} error estimators are derived for model problems, which are
useful in adaptive mesh refinement algorithms. It is important to note
that the best approximation results are key estimates in establishing the
optimality of adaptive finite element methods.
To the authors' best knowledge this is the first attempt of  discussing the error analysis of DG methods under minimal regularity for the optimal control of the Stokes equations with pointwise control
constraints.

This paper is organized as follows. Section \ref{sec:Abstract Setting} set up the abstract framework for the error analysis of discontinuous finite element methods and derives therein some abstract error estimates that form the basis for \textit{a priori}. Subsection 2.2 deals with \textit{a posteriori} error analysis. Section \ref{sec3} introduces two model examples that are under discussion. In section \ref{sec4}, we develop the discrete setting and discuss the applications to the model problems introduced in section \ref{sec3}. Section \ref{sec5} presents some numerical examples to illustrate the theoretical results. 

\section{Abstract Setting}\label{sec:Abstract Setting}

In this section, we develop an abstract framework for the error analysis of discontinuous and nonconforming finite element methods for approximating the solutions of optimal control problems with either boundary control or distributed control. We will assume all the vector spaces are real.

Let $X$ and $M$ are Hilbert spaces with the norm $\lVert\cdot\rVert_{X}$ and $\lVert\cdot\rVert_{M}$, respectively. We denote $V:=X\times M$ the admissible pair of spaces for state variables (velocity and pressure, respectively) and adjoint states. Let $V':=X'\times M'$, where $X'$ and $M'$ are dual of $X$ and $M$, respectively. Let $W$ be a Hilbert space such that $X \subset W \subset X'$ and the inclusions are continuous. The inner product and the norm on $W$ are denoted by $\langle \cdot,\cdot \rangle_{W}$ and $\lVert \cdot \rVert_{W}$, respectively. Let $Q$ be a Hilbert space that will be used for seeking the control variable. The norm and  inner product on $Q$ will be respectively denoted by $\lVert\cdot \rVert_{Q}$ and $\langle \cdot , \cdot \rangle_{Q}$. Let $E:X\rightarrow Q$ be a continuous linear operator. Let $Q_{ad}\subset Q$ be a nonempty closed convex subset.

Assume that $\big((\mathbf{u},p),(\bm{\phi},r),\mathbf{y}\big)\in V \times V\times Q_{ad}$ solves the following optimality system
\begin{subequations}\label{eq1}
\begin{align}
a(\mathbf{u},\mathbf{z})+b(\mathbf{z},p) &={\langle \mathbf{y},E\mathbf{z} \rangle}_{Q} +{\langle \mathbf{f},\mathbf{z} \rangle}_{W} \;\;\;
{\rm for~all}\;\mathbf{z} \in X,\label{eq:BF} \\ 
b(\mathbf{u},w)&=0\;\;\;\hspace{3.1 cm}{\rm for~all}\;w \in M,\label{2.1b}\\
a(\mathbf{z},\bm{\phi})-b(\mathbf{z},r)&={\langle \mathbf{u}-\mathbf{u}_d,\mathbf{z} \rangle}_{W} \;\;\;\hspace{1.0 cm}{\rm for~all}\;\mathbf{z} \in X, \label{eq:ABF}\\
b(\bm{\phi},w)&=0\;\;\;\hspace{3.1 cm}{\rm for~all}\;w \in M,\label{2.1d}\\
{\langle E\bm{\phi}+\lambda \mathbf{y}, \mathbf{x}-\mathbf{y} \rangle }_{Q} &\geq 0 \;\;\;\hspace{3.1 cm}{\rm for~all} ~\mathbf{x}\in Q_{ad},\label{eq:VI}
\end{align}
\end{subequations}
where  $\mathbf{f}\in W,~{\bf u}_d\in W $, $\lambda>0$ are given and $a:X\times X \rightarrow \mathbb{R}$, $b:X\times M \rightarrow \mathbb{R}$ are continuous bilinear forms in the sense that there exist $C_{1}, C_{2}>0$ such that
\begin{align*}
|a(\mathbf{v},\mathbf{z})|&\leq C_{1}\norm{\mathbf{v}}_{X} \lVert \mathbf{z} \rVert_{X},\;\;
|b(\mathbf{z},p)|\leq C_{2}\lVert \mathbf{z} \rVert_{X} \lVert p \rVert_{M},
\end{align*}
for all $\mathbf{v},$ $\mathbf{z} \in X$ and $p \in M.$ The set $Z=\left\{\mathbf{v}\in X :  \forall w \in M,~~ b(\mathbf{v},w)=0\right\}$ and $a$ is $Z$-elliptic, that is, there exists a constant $\alpha>0$ such that $\alpha \norm{\mathbf{v}}_X^2 \leq a(\mathbf{v},\mathbf{v})\; {\rm for~all}~ \mathbf{v}\in Z,$ and $b$ satisfies the inf-sup condition
 which is given by
\begin{align*}
\inf_{p\in M}\sup_{\mathbf{z}\in X} \frac{b(\mathbf{z},p)}{\lVert \mathbf{z} \rVert_{X} \lVert p \rVert}_M& \geq C,~~{\rm for ~some} \;C>0.
\end{align*}
\begin{remark} We have assumed that the system (\ref{eq1}) is well-posed and in Section \ref{sec3}, we have  verified the well-posedness of  (\ref{eq1}) for two model problems.
\end{remark}
Now we introduce corresponding discrete setting. Let $X_{h}\subset W$ be a finite dimensional subspace and $\norm{\cdot}_{h}$ be a norm on $X+X_{h}$ such that $\norm{\mathbf{v}}_{h}=\norm{\mathbf{v}}_{X}$ for all $\mathbf{v}\in X.$ Also let $M_{h}$ be a finite dimensional subspace of $M$ and the norm on $M_{h}$ is $\norm{\cdot}_{M}.$ $a_{h}:X_{h}\times X_{h} \rightarrow \mathbb{R}$, $b_h:X_{h}\times M_{h} \rightarrow \mathbb{R}$ are continuous bilinear forms in the sense that there exist $\bar{C_{1}},\bar{C_{2}}>0$ such that
\begin{align*}
|a_{h}(\mathbf{v}_{h},\mathbf{z}_{h})|&\leq \bar{C_{1}}\lVert  \mathbf{v}_h \rVert_{h} \lVert \mathbf{z}_h \rVert_{h},~~~
|b_{h}(\mathbf{z}_{h},p_{h})|\leq \bar{C_{2}}\lVert \mathbf{z}_h \rVert_{h} \lVert  p_h \rVert_{M},
\end{align*}
for all $\mathbf{v}_h,$ $\mathbf{z}_h \in X_h$ and $p_h \in M_h.$
The set
\begin{equation}\label{zhh}
Z_{h}=\left\{\mathbf{v}_{h}\in X_{h} : \forall w_{h} \in M_{h},~~ b_{h}(\mathbf{v}_{h},w_{h})=0 \right\}
\end{equation}
 and $a_{h}$ is $Z_{h}$-elliptic, that is, there exists $\bar{\alpha}>0$ independent of mesh-size such that $\bar{\alpha} \lVert  \mathbf{v}_h \rVert_h^2 \leq a_{h}(\mathbf{v}_{h},\mathbf{v}_{h})\; {\rm for~all~} \mathbf{v}_{h} \in Z_{h}$  and $b_{h}$ satisfies the inf-sup condition $$ \inf_{p_{h}\in M_{h}}\sup_{\mathbf{z}_{h}\in X_{h}}\frac{ b_{h}(\mathbf{z}_{h},p_{h})}{\lVert \mathbf{z}_h \rVert_{h} \lVert  p_h \rVert_{M}}\geq \bar{C},~~{\rm ~for~ some}~\bar{C}>0,$$
  which is independent of mesh-size.
Similarly, assume that $Q_{h}\subseteq Q$ is a finite dimensional subspace and $Q_{ad}^h \subset Q_{ad}$ is nonempty closed convex subset of $Q_{h}.$ We  denote $V_{h}:=X_{h}\times M_{h}.$
$\newline$
Assume that
 $((\mathbf{u}_{h},p_{h}),(\bm{\phi}_{h},r_{h}),\mathbf{y}_h)\in V_{h} \times V_{h}\times Q_{ad}^h$ solves the following optimality system:
\begin{subequations}\label{eq2}
\begin{align}
a_{h}(\mathbf{u}_{h},\mathbf{z}_{h})+b_{h}(\mathbf{z}_{h},p_{h}) &={\langle \mathbf{y}_h,E_{h}\mathbf{z}_{h} \rangle}_{Q} +{\langle \mathbf{f},\mathbf{z}_{h} \rangle}_{W} \;\;{\rm~for~all}\;\mathbf{z}_{h} \in X_{h},\label{eq:DBF} \\
b_{h}({\bf u}_{h},w_{h})&=0\;\;\;\hspace{3.7 cm}{\rm~for~all}\;w_{h} \in M_{h},\label{2.2b} \\
a_{h}(\mathbf{z}_{h},{\bm \phi}_{h})-b_h(\mathbf{z}_{h},r_{h})& ={\langle \mathbf{u}_{h}-\mathbf{u}_d,\mathbf{z}_{h} \rangle}_{W} \;\;\;\hspace{1.3 cm}{\rm~for~all}\;\mathbf{z}_{h} \in X_{h},\label{eq:DABF}\\
b_{h}(\bm{\phi}_{h},w_{h})&=0\;\;\;\hspace{3.8 cm}{\rm~for~all}\;w_{h} \in M_{h},\label{2.2d}\\
{\langle E_{h}\bm{\phi}_{h}+\lambda \mathbf{y}_h, \mathbf{x}_h-\mathbf{y}_h \rangle }_{Q}& \geq 0 \;\;\;\hspace{3.8 cm}{\rm~for~all} ~\mathbf{x}_h\in Q_{ad}^h,\label{eq:DVI}
\end{align}
\end{subequations}
where $E_{h}:X_{h}+X\rightarrow Q$ is a discrete counterpart of $E$ such that $E_{h}\mathbf{v}=E\mathbf{v}$ for all $\mathbf{v}\in X.$
\newline
Throughout this section, we assume that the following hold true:\\
\textbf{Assumption I:} For all $\mathbf{v}\in X+X_{h}$,
\begin{align}
\norm{\mathbf{v}}_{W}&\leq C \norm{\mathbf{v}}_{h},\;\label{Assmption1}
\end{align}
where $C$ is independent of mesh-size.\\
\textbf{Assumption II:} For all $ \mathbf{v}\in X,$ and $\mathbf{v}_{h} \in X_{h}$,
\begin{align}
\norm{E_{h}(\mathbf{v}-\mathbf{v}_{h})}_{Q}&\leq \norm{\mathbf{v}-\mathbf{v}_{h}}_{h}.\label{Assmption2}
\end{align}
\textbf{Assumption III:}
 The $Q$-projection defined as: For given $\mathbf{x}\in Q$, let $\Pi_{h}\mathbf{x} \in Q_{h}$ be the solution of
\begin{equation}
\langle\Pi_{h}\mathbf{x}- \mathbf{x}, \mathbf{x}_{h}\rangle_{Q}=0\;\;\; {\rm for~all}~ \mathbf{x}_{h}\in Q_{h}.\label{Assmption 3}
\end{equation}
Assume that $\Pi_{h}\mathbf{x} \in Q_{ad}^h,$ whenever $\mathbf{x}\in Q_{ad}.$

\subsection{A priori Error Analysis}
To derive some abstract \textit{a priori} error analysis, we introduce some projections as follows: Let $P_{h}\mathbf{u} \in X_{h}$ , $\bar{P}_{h}\bm{\phi} \in X_{h}$, $R_{h}p \in M_{h}$ and $\bar{R}_{h}r \in M_{h}$ solve
\begin{subequations}
\begin{align}
a_{h}(P_{h}\mathbf{u},\mathbf{z}_{h})+b_{h}(\mathbf{z}_{h},R_{h}p) &={\langle \mathbf{y},E_{h}\mathbf{z}_{h} \rangle}_{Q} +{\langle \mathbf{f},\mathbf{z}_{h} \rangle}_{W} \;{\rm for~all}\;\mathbf{z}_{h} \in X_{h},\label{APr} \\
b_{h}(P_{h}\mathbf{u},w_{h})&=0\;\hspace{3.7 cm}{\rm for~all}\;w_{h} \in M_{h},\label{2.6b}\\
a_{h}(\mathbf{z}_{h},\bar{P_{h}}\bm{\phi})-b_{h}(\mathbf{z}_{h},\bar{R}_{h}r) &={\langle \mathbf{u}-\mathbf{u}_d,\mathbf{z}_{h} \rangle}_{W} \;\hspace{1.4 cm}{\rm for~all}\;\mathbf{z}_{h} \in X_{h},\label{AAPr}\\
b_{h}(\bar{P}_{h}\bm{\phi},w_{h})&=0\;\hspace{3.7 cm}{\rm for~all}\;w_{h} \in M_{h}\label{2.6d}.
\end{align}
\end{subequations}
Here,  we assumed that the bilinear forms $a_h$ and $b_h$ are continuous, $a_h$ is $Z_h$-elliptic, and $b_h$ is inf-sup stable. Also, the right-hand side of (\ref{APr}) is a bounded linear functional on $X_h$. Hence the system  (\ref{APr})-(\ref{2.6b}) has a unique solution \cite[pp. 112]{Girault:1979}. Similarly, the system (\ref{AAPr})-(\ref{2.6d}) is well-posed.\\

\noindent
To derive \textit{a priori} error estimate of control, we need the following lemma.
\begin{lemma}\label{lem2.1} For all $\mathbf{x}_h\in Q_{ad}^h$, it holds
\begin{align}\label{lm21}
{\langle E_{h}(\bm{\phi}_{h}-\bar{P_{h}}\bm{\phi}),\mathbf{y}-\mathbf{y}_h \rangle}_{Q} \geq {\langle E_{h}(\bm{\phi}-\bar P_{h}\bm{\phi}),\mathbf{y}-\mathbf{y}_h \rangle}_{Q} + {\langle E_{h}\bm{\phi}_{h}+\lambda \mathbf{y}_h,\mathbf{y}-\mathbf{x}_h \rangle}_{Q}+ \lambda \lVert   \mathbf{y}-\mathbf{y}_h \rVert_{Q}^2.\end{align}
\end{lemma}
\begin{proof} For all $\mathbf{x}_h\in Q_{ad}^h$, we have $${\langle E_{h}\bm{\phi}_{h}+\lambda \mathbf{y}_h,\mathbf{y}-\mathbf{y}_h \rangle}_{Q} = {\langle E_{h}\bm{\phi}_{h}+\lambda \mathbf{y}_h,\mathbf{y}-\mathbf{x}_h \rangle}_{Q} + {\langle E_{h}\bm{\phi}_{h}+\lambda \mathbf{y}_h,\mathbf{x}_h-\mathbf{y}_h \rangle}_{Q}.$$
From (\ref{eq:DVI}), 
	\begin{equation}\label{1}
	{\langle E_{h}\bm{\phi}_{h}+\lambda \mathbf{y}_h,\mathbf{y}-\mathbf{y}_h \rangle}_{Q} \geq {\langle E_{h}\bm{\phi}_{h}+\lambda \mathbf{y}_h,\mathbf{y}-\mathbf{x}_h \rangle}_{Q}.
	\end{equation}
	The substitution $\mathbf{x}= \mathbf{y}_h \in Q_{ad}^h (\subseteq Q_{ad})$ and $E_{h}\bm{\phi}=E\bm{\phi}\;\; {\rm for}\;\; \bm{\phi} \in X$ in (\ref{eq:VI}), imply
	\begin{equation}\label{2}
	{-\langle E_{h}\bm{\phi}+\lambda \mathbf{y}, \mathbf{y}-\mathbf{y}_h \rangle }_{Q} \geq 0.
	\end{equation}
  An addition of (\ref{1}) and (\ref{2}), yields
	$${\langle E_{h}(\bm{\phi}_{h}-\bm{\phi})+\lambda(\mathbf{y}_h-\mathbf{y}),\mathbf{y}-\mathbf{y}_h \rangle}_{Q} \geq {\langle E_{h}\bm{\phi}_{h}+\lambda \mathbf{y}_h,\mathbf{y}-\mathbf{x}_h \rangle}_{Q}.$$ 
	%
	Hence,
	$${\langle E_{h}(\bm{\phi}_{h}-\bar{P_{h}}\bm{\phi}),\mathbf{y}-\mathbf{y}_h \rangle}_{Q} \geq \lambda \lVert   \mathbf{y}-\mathbf{y}_h \rVert_{Q}^2 + {\langle E_{h}(\bm{\phi}-\bar P_{h}\bm{\phi}),\mathbf{y}-\mathbf{y}_h \rangle}_{Q} + {\langle E_{h}\bm{\phi}_{h}+\lambda \mathbf{y}_h,\mathbf{y}-\mathbf{x}_h \rangle}_{Q}.$$
	This concludes the proof
\end{proof}
Following theorem gives an \textit{a priori} error estimate for the control variable.
\begin{theorem}\label{thm2.2} Let $ \big((\mathbf{u},p),(\bm{\phi},r),\mathbf{y}\big)$ and $((\mathbf{u}_{h},p_{h}),(\bm{\phi}_{h},r_{h}),\mathbf{y}_h) $ solve systems (\ref{eq1}),  and (\ref{eq2}), respectively.
 Then it holds
	\begin{align}\label{thm22}
\lVert   \mathbf{y}-\mathbf{y}_h \rVert_{Q}^2+\lVert \mathbf{u}-\mathbf{u}_{h}\rVert_{W}^2 \leq C(\lVert E\bm{\phi}-\Pi_{h}(E\bm{\phi})\rVert_{Q}^2 +\lVert \mathbf{y}-\Pi_{h}\mathbf{y}\rVert_{Q}^2+\lVert \bm{\phi}-\bar{P_{h}}\bm{\phi}\rVert_{h}^2+\lVert \mathbf{u}-P_{h}\mathbf{u} \rVert_{W}^2).
\end{align}
\end{theorem}

\begin{proof}
	From (\ref{APr})-(\ref{eq:DBF}) and (\ref{2.6b})-(\ref{2.2b}) for all $(\mathbf{z}_{h},w_h )\in X_{h}\times M_h$ we have
	\begin{align}
	a_{h}(P_{h}\mathbf{u}-\mathbf{u}_{h},\mathbf{z}_{h})+b_{h}(\mathbf{z}_{h},R_{h}p-p_{h}) &={\langle \mathbf{y}-\mathbf{y}_h,E_{h}\mathbf{z}_{h} \rangle}_{Q},\label{3}\\
	b_{h}(P_{h}\mathbf{u}-\mathbf{u}_{h},w_{h})&=0. \label{3a}
	\end{align}
	Similarly, from (\ref{AAPr})-(\ref{eq:DABF}) and (\ref{2.6d})-(\ref{2.2d}) for all $(\mathbf{z}_{h},w_h )\in X_{h}\times M_h$, we have
	\begin{align}
a_{h}(\mathbf{z}_{h},\bar{P_{h}}\bm{\phi}-\bm{\phi}_{h})+b_{h}(\mathbf{z}_{h},-\bar{R_{h}}r+r_{h}) &={\langle \mathbf{u}-\mathbf{u}_{h},\mathbf{z}_{h} \rangle}_{W},\label{4}\\
	b_{h}(\bar{P_{h}}\bm{\phi}-\bm{\phi}_{h},w_{h})&=0.\label{4b}
	\end{align}
The substitution of
	 $\mathbf{z}_{h}= \bar{P_{h}}\bm{\phi}-\bm{\phi}_{h}$ in (\ref{3}) and $\mathbf{z}_{h}= P_{h}\mathbf{u}-\mathbf{u}_{h}$ in (\ref{4}), and use of the fact that 	 $b_{h}(P_{h}\mathbf{u}-\mathbf{u}_{h},-\bar{R_{h}}r+r_{h})=0$ in  (\ref{3}) and $b_{h}(\bar{P_{h}}\bm{\phi}-\bm{\phi}_{h},R_{h}p-{p_{h}})=0 $ in (\ref{4}), and finally  subtraction of the resulting equations lead to
	$${\langle \mathbf{y}-\mathbf{y}_h,E_h(\bar{P_{h}}\bm{\phi}-\bm{\phi}_{h}) \rangle}_{Q} - {\langle \mathbf{u}-\mathbf{u}_{h},P_{h}\mathbf{u}-\mathbf{u}_{h} \rangle}_{W}=0.$$
Further, we have
\begin{align}\label{lm22}
{\langle \mathbf{y}-\mathbf{y}_h,E_h(\bm{\phi}_{h}-\bar{P_{h}}\bm{\phi}) \rangle}_{Q} +\lVert P_{h}\mathbf{u}-\mathbf{u}_{h}\rVert_{W}^2= {\langle \mathbf{u}-P_{h}\mathbf{u},\mathbf{u}_{h}-P_{h}\mathbf{u} \rangle}_{W}.
\end{align}
An addition of  (\ref{lm22}) and (\ref{lm21}) from  Lemma \ref{lem2.1} with $\mathbf{x}_h= \Pi_{h}\mathbf{y}$ yields
	\begin{align}\label{k1}
	\lambda \lVert   \mathbf{y}-\mathbf{y}_h \rVert_{Q}^2 + \lVert P_{h}\mathbf{u}-\mathbf{u}_{h}\rVert_{W}^2 \leq &- {\langle E_h\bm{\phi}_{h}+\lambda \mathbf{y}_h,\mathbf{y}-\mathbf{x}_h \rangle}_{Q}\nonumber\\
		& - {\langle E_h(\bm{\phi}-\bar P_{h}\bm{\phi}),\mathbf{y}-\mathbf{y}_h \rangle}_{Q}
	 + {\langle \mathbf{u}-P_{h}\mathbf{u},\mathbf{u}_{h}-P_{h}\mathbf{u} \rangle}_{W}.
	\end{align}
	An addition and subtraction of some terms in the first term on the right-hand side of (\ref{k1}) shows
	\begin{align}\label{k2}
 {\langle E_h\bm{\phi}_{h}+\lambda \mathbf{y}_h,\mathbf{y}-\mathbf{x}_h \rangle}_{Q}&= -{\langle E_h(\bm{\phi}_{h}-\bm{\phi}+\bm{\phi})+\lambda( \mathbf{y}_h-\mathbf{y}+\mathbf{y}),\mathbf{y}-\mathbf{x}_h \rangle}_{Q}\nonumber \\
	=& -\langle E_h\bm{\phi}-\Pi_{h}(E\bm{\phi})+\Pi_{h}(E\bm{\phi})+\lambda(\mathbf{y}-\mathbf{x}_h+\mathbf{x}_h), \mathbf{y}-\mathbf{x}_h\rangle_{Q}\nonumber\\
	& -{\langle E_h(\bm{\phi}_{h}-\bm{\phi})+\lambda( \mathbf{y}_h-\mathbf{y}),\mathbf{y}-\mathbf{x}_h \rangle}_{Q}\nonumber\\
	=&  -\langle E\bm{\phi}-\Pi_{h}(E\bm{\phi})+\lambda(\mathbf{y}-\mathbf{x}_h), \mathbf{y}-\mathbf{x}_h\rangle_{Q} -\langle \Pi_{h}(E\bm{\phi})+\lambda \mathbf{x}_h, \mathbf{y}-\mathbf{x}_h\rangle_{Q}\nonumber\\
	& -{\langle E_h(\bm{\phi}_{h}-\bm{\phi})+\lambda( \mathbf{y}_h-\mathbf{x}_h+\mathbf{x}_h-\mathbf{y}),\mathbf{y}-\mathbf{x}_h \rangle}_{Q}\nonumber\\
	=& -\langle E\bm{\phi}-\Pi_{h}(E\bm{\phi})+\lambda(\mathbf{y}-\mathbf{x}_h), \mathbf{y}-\mathbf{x}_h\rangle_{Q}\nonumber\\
	& -{\langle E_h(\bm{\phi}_{h}-\bm{\phi})+\lambda(\mathbf{x}_h-\mathbf{y}),\mathbf{y}-\mathbf{x}_h \rangle}_{Q}.
	\end{align}
	Here, the selection of $\mathbf{x}_h= \Pi_{h}\mathbf{y}$ and the $Q$-orthogonal projection property imply
	$$\langle \Pi_{h}(E\bm{\phi})+\lambda \mathbf{x}_h, \mathbf{y}-\mathbf{x}_h\rangle_{Q}-{\langle \lambda( \mathbf{y}_h-\mathbf{x}_h),\mathbf{y}-\mathbf{x}_h \rangle}_{Q}=0.$$
	The substitution of (\ref{k2}) in (\ref{k1}) and a use of the Young inequality imply
	\begin{align}
	\lambda \lVert   \mathbf{y}-\mathbf{y}_h \rVert_{Q}^2 + \lVert P_{h}\mathbf{u}-\mathbf{u}_{h}\rVert_{W}^2	\leq &~ \big(\frac{1}{2}\lVert E\bm{\phi}-\Pi_{h}(E\bm{\phi})\rVert_{Q}^2+\frac{1}{2}\lVert \mathbf{y}-\mathbf{x}_h\rVert_{Q}^2+\frac{1}{2\lambda}\norm{E_h(\bm{\phi}-\bar{P_{h}}\bm{\phi})}_{Q}^2\big)\nonumber\\
	&+ \big(\norm{E_h(\bm{\phi}_{h}-\bm{\phi})}_{Q} \lVert \mathbf{y}-\mathbf{x}_h\rVert_{Q}+\frac{1}{2}\lVert \mathbf{u}-P_{h}\mathbf{u} \rVert_{W}^2\big)\nonumber\\
	&+\frac{1}{2}\lVert \mathbf{u}_{h}-P_{h}\mathbf{u}\rVert_{W}^2+\frac{\lambda}{2} \lVert   \mathbf{y}-\mathbf{y}_h \rVert_{Q}^2. \label{5}
	\end{align}
In order to estimate  $\norm{E_h(\bm{\phi}_{h}-\bm{\phi})}_{Q}$, the choice of $\mathbf{z}_{h}=\bar{P_{h}}\bm{\phi}-\bm{\phi}_{h}$ in (\ref{4}) shows
	$$a_{h}(\bar{P_{h}}\bm{\phi}-\bm{\phi}_{h},\bar{P_{h}}\bm{\phi}-\bm{\phi}_{h})+b_{h}(\bar{P_{h}}\bm{\phi}-\bm{\phi}_{h},-\bar{R_{h}}r+r_{h}) ={\langle \mathbf{u}-\mathbf{u}_{h},\bar{P_{h}}\bm{\phi}-\bm{\phi}_{h} \rangle}_{W}.$$
	Since $\bar{P_{h}}\bm{\phi}-\bm{\phi}_{h} \in Z_h$, the $Z_h$-ellipticity of $a_h$ yields
	$$C\lVert \bar{P_{h}}\bm{\phi}-\bm{\phi}_{h}\rVert_{h}^2 \leq \lVert \mathbf{u}-\mathbf{u}_{h}\rVert_{W} \lVert \bar{P_{h}}\bm{\phi}-\bm{\phi}_{h}\rVert_{W}.$$
	From  Assumption I, 
	we have
	$\lVert \bar{P_{h}}\bm{\phi}-\bm{\phi}_{h}\rVert_{W}\leq \lVert \bar{P_{h}}\bm{\phi}-\bm{\phi}_{h}\rVert_{h}.$ Hence
	\begin{equation}
	 \lVert \bar{P_{h}}\bm{\phi}-\bm{\phi}_{h}\rVert_{h} \leq C \lVert \mathbf{u}-\mathbf{u}_{h}\rVert_{W}\label{6}.
	\end{equation}
	Also, from Assumption II, we have 
	\begin{align}
	\lVert E_h(\bm{\phi}-\bm{\phi}_{h})\rVert_{Q}& \leq C\lVert \bm{\phi}-\bm{\phi}_{h}\rVert_h \leq C \big(\lVert \bm{\phi}-\bar{P_{h}}\bm{\phi}\rVert_{h}+\lVert \bar{P_{h}}\bm{\phi}-\bm{\phi}_{h}\rVert_{h})\nonumber\\
	&\leq C \big(\lVert \bm{\phi}-\bar{P_{h}}\bm{\phi}\rVert_{h}+\lVert \mathbf{u}-\mathbf{u}_{h}\rVert_{W}\big)\nonumber\\
	& \leq C \big(\lVert \bm{\phi}-\bar{P_{h}}\bm{\phi}\rVert_{h}+\lVert \mathbf{u}-P_{h}\mathbf{u} \rVert_{W}+{\lVert P_{h}\mathbf{u}-\mathbf{u}_{h}\rVert_{W}}\big)\label{7}.
	\end{align}
	A substitution of (\ref{7}) in (\ref{5}), then  a use of the Young inequality
	with  kick-back the term $\lVert P_{h}\mathbf{u}-\mathbf{u}_{h}\rVert^2_{W}$ (to LHS) and  the orthogonality  property of  $Q$-projection from Assumption III in the estimates of $ \norm{\mathbf{y}-\mathbf{x}_h}^2_Q$ result in 
	\begin{align}\label{k4}
	\lVert   \mathbf{y}-\mathbf{y}_h \rVert_{Q}^2+\lVert P_h \mathbf{u}-\mathbf{u}_{h}\rVert_{W}^2 \leq& ~C(\lVert E\bm{\phi}-\Pi_{h}(E\bm{\phi})\rVert_{Q}^2 +\lVert \mathbf{y}-\Pi_{h}\mathbf{y}\rVert_{Q}^2\nonumber \\&+\lVert \bm{\phi}-\bar{P_{h}}\bm{\phi}\rVert_{h}^2+\lVert \mathbf{u}-P_{h}\mathbf{u} \rVert_{W}^2).
	\end{align}
Finally, the triangle inequality in the term  $\norm{\mathbf{u}-\mathbf{u}_{h}}_{W}$ and (\ref{k4}) lead to (\ref{thm22}). This concludes the proof. \end{proof}
\begin{theorem}\label{thm 2.3}
	It  holds,
	$$ \lVert \bm{\phi}-\bm{\phi}_{h}\rVert_{h}\leq C \big(\lVert \bm{\phi}-\bar{P_{h}}\bm{\phi}\rVert_{h}+\lVert E\bm{\phi}-\Pi_{h}(E\bm{\phi})\rVert_{Q} +\lVert \mathbf{y}-\Pi_{h}\mathbf{y}\rVert_{Q}+\lVert \mathbf{u}-P_{h}\mathbf{u} \rVert_{W}\big)$$ and
	\begin{equation}\label{thm232}
	\lVert \mathbf{u}-\mathbf{u}_{h}\rVert_{h} \leq C \big(\lVert \mathbf{u}-P_{h}\mathbf{u} \rVert_{h}+\lVert \bm{\phi}-\bar{P_{h}}\bm{\phi}\rVert_{h}+\lVert E\bm{\phi}-\Pi_{h}(E\bm{\phi})\rVert_{Q} +\lVert \mathbf{y}-\Pi_{h}\mathbf{y}\rVert_{Q}\big).
	\end{equation}
\end{theorem}
\begin{proof}
	From (\ref{6}) and (\ref{thm22}), we have
	\begin{align*}
	\lVert \bm{\phi}-\bm{\phi}_{h}\rVert_{h}&\leq \lVert \bm{\phi}-\bar{P_{h}}\bm{\phi}\rVert_{h}+\lVert \bar{P_{h}}\bm{\phi}-\bm{\phi}_{h}\rVert_{h}
	\leq \lVert \bm{\phi}-\bar{P_{h}}\bm{\phi}\rVert_{h}+ C \lVert \mathbf{u}-\mathbf{u}_{h}\rVert_W\\
	&\leq  C \big(\lVert \bm{\phi}-\bar{P_{h}}\bm{\phi}\rVert_{h}+\lVert E\bm{\phi}-\Pi_{h}(E\bm{\phi})\rVert_{Q} +\lVert \mathbf{y}-\Pi_{h}\mathbf{y}\rVert_{Q}+\lVert \mathbf{u}-P_{h}\mathbf{u} \rVert_{W}\big).\;\;
	\end{align*}
	For the estimate of $\lVert \mathbf{u}-\mathbf{u}_{h}\rVert_{h},$ the substitution $\mathbf{z}_{h}=P_{h}\mathbf{u}-\mathbf{u}_{h}$ in (\ref{3})  and use of (\ref{3a}) imply
	$$a_{h}(P_{h}\mathbf{u}-\mathbf{u}_{h},P_{h}\mathbf{u}-\mathbf{u}_{h})={\langle \mathbf{y}-\mathbf{y}_h,E_h(P_{h}\mathbf{u}-\mathbf{u}_{h}) \rangle}_{Q}.$$
	The $Z_h$-ellipticity and Assumption II give
	$\lVert P_{h}\mathbf{u}-\mathbf{u}_{h}\rVert_{h}\leq C \lVert   \mathbf{y}-\mathbf{y}_h \rVert_{Q}.$
	Hence
	\begin{align*}
	\lVert \mathbf{u}-\mathbf{u}_{h}\rVert_{h}&\leq \lVert \mathbf{u}-P_{h}\mathbf{u} \rVert_{h}+\lVert P_{h}\mathbf{u}-\mathbf{u}_{h}\rVert_{h},
	\leq \lVert \mathbf{u}-P_{h}\mathbf{u} \rVert_{h}+C \lVert   \mathbf{y}-\mathbf{y}_h \rVert_{Q},\\
	&\leq C  \big(\lVert \mathbf{u}-P_{h}\mathbf{u} \rVert_{h}+\lVert \bm{\phi}-\bar{P_{h}}\bm{\phi}\rVert_{h}+\lVert E\bm{\phi}-\Pi_{h}(E\bm{\phi})\rVert_{Q} +\lVert \mathbf{y}-\Pi_{h}\mathbf{y}\rVert_{Q}\big).
	\end{align*}
 Here, the last inequality follows by using Theorem \ref{thm2.2}. This concludes the proof.	
\end{proof}
The following theorem  gives an error estimate for the pressure.
\begin{theorem}\label{thm 2.4} It holds
\begin{align}\label{thm2.41}
\lVert p-p_{h}\rVert_{M}\leq &C\big( \lVert p-R_{h}p\rVert_{M}+ \lVert \mathbf{u}-P_{h}\mathbf{u} \rVert_{h}+\lVert \bm{\phi}-\bar{P_{h}}\bm{\phi}\rVert_{h}\nonumber\\
&+\lVert E\bm{\phi}-\Pi_{h}(E\bm{\phi})\rVert_{Q} +\lVert \mathbf{y}-\Pi_{h}\mathbf{y}\rVert_{Q}\big).
\end{align}
\end{theorem}

\begin{proof} The error term  $\lVert p-p_{h}\rVert_{M}$ can be bounded by the sum of $\lVert p-R_{h}p\rVert_{M}$ and $\lVert R_{h}p-p_{h}\rVert_{M}$ with   $R_{h}$ projection defined in (\ref{APr})-(\ref{2.6b}). Here it is sufficient to estimate $\lVert R_{h}p-p_{h}\rVert_{M}$ to get (\ref{thm2.41}).
%
	From the inf-sup condition we have
	$$\bar{C} \lVert R_{h}p-p_{h}\rVert_{M}\leq \sup_{\mathbf{z}_{h}\in X_{h}} \frac{b_{h}(\mathbf{z}_{h},R_{h}p-p_{h})}{\lVert \mathbf{z}_h \rVert_{h}}.$$
	Use of (\ref{3})  and Assumption II imply
	\begin{align*}
	\bar{C} \lVert R_{h}p-p_{h}\rVert_{M}&\leq \sup_{\mathbf{z}_{h}\in X_{h}} \frac{a_{h}(\mathbf{u}_{h}-P_{h}\mathbf{u},\mathbf{z}_{h})+{\langle \mathbf{y}-\mathbf{y}_h,E_h\mathbf{z}_{h} \rangle}_{Q}}{\lVert \mathbf{z}_h \rVert_{h}},\\
	&\leq C (\lVert \mathbf{u}_{h}-P_{h}\mathbf{u}\rVert_{h}+ \lVert   \mathbf{y}-\mathbf{y}_h \rVert_{Q}).
	\end{align*}
The split of term $\lVert \mathbf{u}_{h}-P_{h}\mathbf{u}\rVert_{h}$  as  $\lVert \mathbf{u}_{h}-\mathbf{u}\rVert_{h} +\lVert \mathbf{u}-P_{h}\mathbf{u}\rVert_{h}$ and the use of estimates from (\ref{thm22}) and (\ref{thm232}) lead to (\ref{thm2.41}).
\end{proof}

Similarly we can derive the following error estimates of adjoint pressure.
\begin{align*}
\norm{r-r_{h}}_{M}\leq C\big( \norm{r-\bar{R}_{h}r}_{M}+ \norm{\mathbf{u}-P_{h}\mathbf{u}}_{h}&+\norm{\bm{\phi}-\bar{P}_{h}\bm{\phi}}_{h}\\&+\norm{E\bm{\phi}-\Pi_{h}(E\bm{\phi})}_{Q} +\norm{\mathbf{y}-\Pi_{h}\mathbf{y}}_{Q}\big).
\end{align*}	
\subsection{{A~posteriori} Error Analysis}
This subsection is devoted to \textit{a~posteriori} error analysis.
 Define reconstructions $R\mathbf{u}\in X$, $\bar{R}\bm{\phi}\in X $ and $R_{0}p\in M$, $\bar{R_{0}}r\in M$ by
\begin{subequations}\label{apo1}
\begin{align}
a(R\mathbf{u},\mathbf{z})+b(\mathbf{z},R_{0}p) &={\langle \mathbf{y}_h,E\mathbf{z} \rangle}_{Q} +{\langle \mathbf{f},\mathbf{z} \rangle}_{W} \;~~\hspace{0.6 cm}{\rm for~all}\;\mathbf{z} \in X,\label{8}\\
b(R\mathbf{u},w)&=0\;~~\hspace{3.9 cm}{\rm for~all}\;w \in M,\label{3.1b}\\
a(\mathbf{z},\bar{R}\bm{\phi})-b(\mathbf{z},\bar{R_{0}}r) &={\langle \mathbf{u}_h-\mathbf{u}_d,\mathbf{z} \rangle}_{W} \;~~\hspace{1.6 cm}{\rm for~all}\;\mathbf{z} \in X,\label{9} \\
b(\bar{R}\bm{\phi},w)&=0\;~~\hspace{3.9 cm}{\rm for~all}\;w \in M.\label{3.1d}
\end{align}
\end{subequations}
The well-posedness of the above system \eqref{8}-\eqref{3.1b} follows from the facts that
 the right-hand side of \eqref{8} is a bounded linear functional on $X$,  the bilinear forms $a$ and $b$ are continuous, $a$ is $Z$-elliptic and $b$ is inf-sup stable, and hence the system  \eqref{8}-\eqref{3.1b} has a unique solution  \cite[pp. 81]{Girault:1979}. Similarly, the system \eqref{9}-\eqref{3.1d} is well-posed.\\

\noindent
From  the systems of equations (\ref{eq1}) and (\ref{apo1}), we have
\begin{subequations}
\begin{align}
a(\mathbf{u}-R\mathbf{u},\mathbf{z})+b(\mathbf{z},p-R_{0}p) &={\langle \mathbf{y}-\mathbf{y}_h,E\mathbf{z} \rangle}_{Q}\; \;~~{\rm for~all}\;\mathbf{z} \in X,\label{10} \\
b(\mathbf{u}-R\mathbf{u},w)&=0\;~~\hspace{2.4 cm}{\rm for~all}\;w \in M,\label{3.2b}\\
a(\mathbf{z},\bm{\phi}-\bar{R}\bm{\phi})+b(\mathbf{z},-r+\bar{R_{0}}r)& ={\langle \mathbf{u}-\mathbf{u}_{h},\mathbf{z} \rangle}_{W} \;~~\hspace{0.3 cm}{\rm for~all}\;\mathbf{z} \in X,\label{11} \\
b(\bm{\phi}-\bar{R}\bm{\phi},w)& =0\;~~\hspace{2.4 cm}{\rm for~all}\;w \in M.\label{3.2d}
\end{align}	
\end{subequations}

\begin{lemma}\label{lem 2.5} For all $\mathbf{x}_h\in Q_{ad}^h$, it holds
	\begin{align}\label{eqn4}
{\langle E_h(\bar{R}\bm{\phi}-\bm{\phi}),\mathbf{y}-\mathbf{y}_h \rangle}_{Q} \geq -{\langle E_h(\bm{\phi}_{h}-\bar{R}\bm{\phi}),\mathbf{y}-\mathbf{y}_h \rangle}_{Q} + {\langle E_h\bm{\phi}_{h}+\lambda \mathbf{y}_h,\mathbf{y}-\mathbf{x}_h \rangle}_{Q}+ \lambda \lVert   \mathbf{y}-\mathbf{y}_h \rVert_{Q}^2.
	\end{align}
	\end{lemma}
\begin{proof} For any $\mathbf{x}_h\in Q_{ad}^h$, we have
	 $${\langle E_h\bm{\phi}_{h}+\lambda \mathbf{y}_h,\mathbf{y}-\mathbf{y}_h \rangle}_{Q} = {\langle E_h\bm{\phi}_{h}+\lambda \mathbf{y}_h,\mathbf{y}-\mathbf{x}_h \rangle}_{Q} + {\langle E_h\bm{\phi}_{h}+\lambda \mathbf{y}_h,\mathbf{x}_h-\mathbf{y}_h \rangle}_{Q}.\; $$
Since  ${\rm for~all} ~\mathbf{x}_h\in Q_{ad}^h$ ~~${\langle E_h\bm{\phi}_{h}+\lambda \mathbf{y}_h, \mathbf{x}_h-\mathbf{y}_h \rangle }_{Q} \geq 0 $, we have
	\begin{equation}\label{12}
	{\langle E_h\bm{\phi}_{h}+\lambda \mathbf{y}_h,\mathbf{y}-\mathbf{y}_h \rangle}_{Q} \geq {\langle E_h\bm{\phi}_{h}+\lambda \mathbf{y}_h,\mathbf{y}-\mathbf{x}_h \rangle}_{Q}. 
	\end{equation}
	Substituting $\mathbf{x}_h= \mathbf{y}_h \in Q_{ad}^h (\subseteq Q_{ad})$ in (\ref{eq:VI}) and $E_h\bm{\phi}=E\bm{\phi}\;\; \text{for}~ \bm{\phi} \in X$ , we get
	\begin{equation}\label{13}
	{-\langle E_h\bm{\phi}+\lambda \mathbf{y}, \mathbf{y}-\mathbf{y}_h \rangle }_{Q} \geq 0.
	\end{equation}
	Adding (\ref{12}) and (\ref{13}) yields
	\begin{equation}\label{eqn1}
{\langle E_h(\bm{\phi}_{h}-\bm{\phi})+\lambda(\mathbf{y}_h-\mathbf{y}),\mathbf{y}-\mathbf{y}_h \rangle}_{Q} \geq {\langle E_h\bm{\phi}_{h}+\lambda \mathbf{y}_h,\mathbf{y}-\mathbf{x}_h \rangle}_{Q}.
	\end{equation}
	Consider the left-hand side of (\ref{eqn4}): we have
	\begin{align}\label{eqn2}
	{\langle E_h(\bar{R}\bm{\phi}-\bm{\phi}),\mathbf{y}-\mathbf{y}_h \rangle}_{Q} \geq -  {\langle E_h(\bm{\phi}_{h}-\bar{R}\bm{\phi}),\mathbf{y}-\mathbf{y}_h \rangle}_{Q}
	+  {\langle E_h(\bm{\phi}_{h}-\bm{\phi}),\mathbf{y}-\mathbf{y}_h \rangle}_{Q}.
	\end{align}
	A use of (\ref{eqn1}) in (\ref{eqn2})  concludes the proof. 
\end{proof}

\begin{theorem}\label{thm2.6}
	It holds
	\begin{align}\label{thm32}
	\lVert   \mathbf{y}-\mathbf{y}_h \rVert_{Q}+\lVert \mathbf{u}-R\mathbf{u}\rVert_{W} \leq C(\lVert E_h\bm{\phi}_{h}-\Pi_{h}(E_h\bm{\phi}_{h})\rVert_{Q} +\lVert \bm{\phi}_{h}-\bar{R}\bm{\phi}\rVert_{h}+\lVert R\mathbf{u}-\mathbf{u}_{h}\rVert_{W}).
	\end{align}
\end{theorem}
\begin{proof}
	With the substitutions $\mathbf{z}=\bm{\phi}-\bar{R}\bm{\phi}$ in (\ref{10}) and $\mathbf{z}= \mathbf{u}-R\mathbf{u}$ in (\ref{11}), we have
	\begin{align}
	a(\mathbf{u}-R\mathbf{u},\bm{\phi}-\bar{R}\bm{\phi})+b(\bm{\phi}-\bar{R}\bm{\phi},p-R_{0}p) &={\langle \mathbf{y}-\mathbf{y}_h,E(\bm{\phi}-\bar{R}\bm{\phi}) \rangle_{Q}}, \label{14} \\
	a(\mathbf{u}-R\mathbf{u},\bm{\phi}-\bar{R}\bm{\phi})+b(\mathbf{u}-R\mathbf{u},-r+\bar{R_{0}}r) &={\langle \mathbf{u}-\mathbf{u}_{h},\mathbf{u}-R\mathbf{u} \rangle}_{W}.\label{15}
	\end{align}
 Since $b(\bm{\phi}-\bar{R}\bm{\phi},p-R_{0}p)=0$ and $b(\mathbf{u}-R\mathbf{u},-r+\bar{R_{0}}r)=0$  from (\ref{3.2d}) and (\ref{3.2b}), the subtraction of
  (\ref{15}) from  (\ref{14}) yields
	$${\langle \mathbf{y}-\mathbf{y}_h,E(\bm{\phi}-\bar{R}\bm{\phi}) \rangle}_{Q}-{\langle \mathbf{u}-\mathbf{u}_{h},\mathbf{u}-R\mathbf{u} \rangle}_{W}=0,$$
and
\begin{align}\label{thm321}
{\langle \mathbf{y}-\mathbf{y}_h,E(\bar{R}\bm{\phi}-\bm{\phi}) \rangle}_{Q}+\lVert \mathbf{u}-R\mathbf{u}\rVert_{W}^2=-{\langle R\mathbf{u}-\mathbf{u}_{h},\mathbf{u}-R\mathbf{u} \rangle}_{W}.
\end{align}
From Lemma \ref{lem 2.5} and (\ref{thm321}), we have
	\begin{align}\label{18}
	\lambda \lVert \mathbf{y}-\mathbf{y}_h \rVert_{Q}^2+&\lVert \mathbf{u}-R\mathbf{u}\rVert^2_{W}\leq  -{\langle E_h\bm{\phi}_{h}+\lambda \mathbf{y}_h,\mathbf{y}-\mathbf{x}_h \rangle}_{Q} + {\langle E_h(\bar{R}\bm{\phi}-\bm{\phi}),\mathbf{y}-\mathbf{y}_h \rangle}_{Q}\nonumber\\
	& + {\langle E_h(\bm{\phi}_{h}-\bar{R}\bm{\phi}),\mathbf{y}-\mathbf{y}_h \rangle}_{Q}-{\langle \mathbf{y}-\mathbf{y}_h,E(\bar{R}\bm{\phi}-\bm{\phi}) \rangle}_{Q}-
	{\langle R\mathbf{u}-\mathbf{u}_{h},\mathbf{u}-R\mathbf{u} \rangle}_{W}.
	\end{align}
	Taking the first term on the right-hand side of (\ref{18}) with $\mathbf{x}_h=\Pi_{h}\mathbf{y}\in Q_{ad}^h$ and using  Assumption III show
	\begin{align}\label{19}
	{\langle E_h\bm{\phi}_{h}+\lambda \mathbf{y}_h,\mathbf{y}-\Pi_{h}\mathbf{y} \rangle}_{Q} &= {\langle E_h\bm{\phi}_{h},\mathbf{y}-\Pi_{h}\mathbf{y}\rangle}_{Q},={\langle E_h\bm{\phi}_{h}-\Pi_{h}(E_h\bm{\phi}_{h}),\mathbf{y}-\Pi_{h}\mathbf{y}\rangle}_{Q},\nonumber\\
	&={\langle E_h\bm{\phi}_{h}-\Pi_{h}(E_h\bm{\phi}_{h}),\mathbf{y}-\mathbf{y}_h\rangle}_{Q}.
	\end{align}
	The substitution of (\ref{19}) in (\ref{18}), a use of Assumption II and the Young inequality result in
	\begin{align}\label{thm322}
	\lambda \lVert   \mathbf{y}-\mathbf{y}_h \rVert_{Q}^2+\lVert \mathbf{u}-R\mathbf{u}\rVert_{W}^2\leq& {\langle E_h\bm{\phi}_{h}-\Pi_{h}(E_h\bm{\phi}_{h}),\mathbf{y}-\mathbf{y}_h\rangle}_{Q} + {\langle E_h(\bm{\phi}_{h}-\bar{R}\bm{\phi}),\mathbf{y}-\mathbf{y}_h \rangle}_{Q}\nonumber\\
	&-{\langle R\mathbf{u}-\mathbf{u}_{h},\mathbf{u}-R\mathbf{u} \rangle}_{W},\nonumber\\
	\leq&~ \frac{1}{2 \lambda}\lVert E_h\bm{\phi}_{h}-\Pi_{h}(E_h\bm{\phi}_{h})\rVert_{Q}^2+\frac{\lambda}{4} \lVert   \mathbf{y}-\mathbf{y}_h \rVert_{Q}^2+\frac{C}{2\lambda}\lVert \bm{\phi}_{h}-\bar{R}\bm{\phi}\rVert_{h}^2,\nonumber\\
	&+\frac{\lambda}{4}\lVert   \mathbf{y}-\mathbf{y}_h \rVert_{Q}^2+\frac{1}{2}\lVert R\mathbf{u}-\mathbf{u}_{h}\rVert_{W}^2+\frac{1}{2}\lVert \mathbf{u}-R\mathbf{u}\rVert_{W}^2.
	\end{align}
	Rearrangement of the terms in (\ref{thm322}) leads to the estimates (\ref{thm32}).
\end{proof}
\begin{theorem}\label{thm2.7}
It holds
\begin{align}
\lVert \mathbf{u}-\mathbf{u}_{h}\rVert_{h} \leq C \big(\lVert R\mathbf{u}-\mathbf{u}_{h}\rVert_{h}+\lVert E_h\bm{\phi}_{h}-\Pi_{h}(E_h\bm{\phi}_{h})\rVert_{Q}+\lVert\bm{\phi}_{h}-\bar{R}\bm{\phi}\rVert_{h}\big),\label{331}
\end{align}
		and
		\begin{align}
		\lVert \bm{\phi}-\bm{\phi}_{h}\rVert_{h}\leq C \big(\lVert\bm{\phi}_{h}-\bar{R}\bm{\phi}\rVert_{h}+\lVert R\mathbf{u}-\mathbf{u}_{h}\rVert_{h}+\lVert E_h\bm{\phi}_{h}-\Pi_{h}(E_h\bm{\phi}_{h})\rVert_{Q}\big).\label{332}
		\end{align}
\end{theorem}
\begin{proof} To prove (\ref{331}), it is sufficient to  estimate  $ \lVert \mathbf{u}-R\mathbf{u}\rVert_{h}$ as
	 the error term  $\lVert \mathbf{u}-\mathbf{u}_{h}\rVert_{h}$ can be bounded by the sum of $\lVert \mathbf{u}-R\mathbf{u}\rVert_{h}$ and $\lVert R\mathbf{u}-\mathbf{u}_{h}\rVert_{h}$. 
	Putting $\mathbf{z}= \mathbf{u}-R\mathbf{u}$ in (\ref{10}), gives
	\begin{align*}
	a(\mathbf{u}-R\mathbf{u},\mathbf{u}-R\mathbf{u})&={\langle \mathbf{y}-\mathbf{y}_h,E(\mathbf{u}-R\mathbf{u}) \rangle}_{Q}.
	\end{align*}
	The Z-ellipticity of $a(\cdot,\cdot)$ and Assumption II imply
\begin{align}\label{326}
\lVert \mathbf{u}-R\mathbf{u}\rVert_{X}\leq C \lVert   \mathbf{y}-\mathbf{y}_h \rVert_{Q}.
\end{align}
Since  $\mathbf{u}-R\mathbf{u}\in X$ and $\lVert \mathbf{u}-R\mathbf{u}\rVert_{X}=\lVert \mathbf{u}-R\mathbf{u}\rVert_{h},$  we have
	$\lVert \mathbf{u}-R\mathbf{u}\rVert_{h}\leq C \lVert   \mathbf{y}-\mathbf{y}_h \rVert_{Q}.$
	Hence
	$$\lVert \mathbf{u}-\mathbf{u}_{h}\rVert_{h} \leq C\lVert   \mathbf{y}-\mathbf{y}_h \rVert_{Q}+ \lVert R\mathbf{u}-\mathbf{u}_{h}\rVert_{h}. $$	
 The estimates for $\lVert   \mathbf{y}-\mathbf{y}_h \rVert_{Q}$ from Theorem \ref{thm2.6} and Assumption I lead to
\begin{align}\label{thm334}
\lVert \mathbf{u}-\mathbf{u}_{h}\rVert_{h} \leq C \big(\lVert R\mathbf{u}-\mathbf{u}_{h}\rVert_{h}+\lVert E_h\bm{\phi}_{h}-\Pi_{h}(E_h\bm{\phi}_{h})\rVert_{Q}+\lVert\bm{\phi}_{h}-\bar{R}\bm{\phi}\rVert_{h}\big).\end{align}
	Now $$\lVert \bm{\phi}-\bm{\phi}_{h}\rVert_{h} \leq	\lVert \bm{\phi}-\bar{R}\bm{\phi}\rVert_{h}+\lVert \bar{R}\bm{\phi}-\bm{\phi}_{h}\rVert_{h}.$$
	Putting $\mathbf{z}=\bm{\phi}-\bar{R}\bm{\phi}$ in (\ref{11}) and using assumption (\ref{Assmption2}), we conclude that
	$$\lVert \bm{\phi}-\bar{R}\bm{\phi}\rVert_{h} \leq C \lVert \mathbf{u}-\mathbf{u}_{h}\rVert_{h},$$
	which implies that
	\begin{align}\label{thm3.3}
	\lVert \bm{\phi}-\bm{\phi}_{h}\rVert_{h} \leq C \lVert \mathbf{u}-\mathbf{u}_{h}\rVert_{h}	+\lVert \bar{R}\bm{\phi}-\bm{\phi}_{h}\rVert_{h}.
	\end{align}
	Finally, substitution of the bounds of  $\lVert \mathbf{u}-\mathbf{u}_{h}\rVert_{h}$ from (\ref{thm334}) in (\ref{thm3.3}) leads to (\ref{332}).
\end{proof}
\begin{theorem}\label{thm2.8} It holds,
	\begin{align*}
	\lVert p-p_{h}\rVert_{M}&\leq C \big(\lVert R_{0}p-p_{h}\rVert_{M}+\lVert E_h\bm{\phi}_{h}-\Pi_{h}(E_h\bm{\phi}_{h})\rVert_{Q}+\lVert R\mathbf{u}-\mathbf{u}_{h}\rVert_{W}+\lVert\bm{\phi}_{h}-\bar{R}\bm{\phi}\rVert_{h} \big).
	\end{align*}
\end{theorem}
\begin{proof}
	From the inf-sup condition, we have
\begin{align}\label{342}
\beta \lVert p-R_{0}p\rVert_{M}\leq \sup_{\mathbf{z}\in X}\frac{b(\mathbf{z}, p-R_{0}p)}{\lVert \mathbf{z} \rVert_{X}}.
\end{align}
	Use of (\ref{10}) on the right-hand side of (\ref{342}) shows
	\begin{align}\label{thm345}
	\beta \lVert p-R_{0}p\rVert_{M}&\leq \sup_{\mathbf{z}\in X}\frac{b(\mathbf{z}, p-R_{0}p)}{\lVert \mathbf{z} \rVert_{X}}
\leq \sup_{\mathbf{z}\in X} \frac{{\langle \mathbf{y}-\mathbf{y}_h,E\mathbf{z} \rangle}_{Q}+a(R\mathbf{u}-\mathbf{u},\mathbf{z})}{\lVert \mathbf{z}\rVert_{X}}\nonumber\\
	&\leq C \big(\lVert \mathbf{y}-\mathbf{y}_h\rVert_{Q}+\lVert R\mathbf{u}-\mathbf{u}\rVert_{X}\big).
	\end{align}
Since from (\ref{326}) we have $\lVert R\mathbf{u}-\mathbf{u} \rVert_{X}=\lVert R\mathbf{u}-\mathbf{u} \rVert_{h} \leq C \lVert \mathbf{y}-\mathbf{y}_h\rVert_{Q}.$  The triangle inequality and (\ref{thm345}) show
	$$\lVert p-p_{h}\rVert_{M}\leq \lVert p-R_{0}p\rVert_{M}+\lVert R_{0}p-p_{h}\rVert_{M} \leq C \big(\lVert \mathbf{y}-\mathbf{y}_h\rVert_{Q}+\lVert R_{0}p-p_{h}\rVert_{M} \big).$$
	Substitution of the expression  $\lVert \mathbf{y}-\mathbf{y}_h\rVert_{Q}$ from (\ref{thm32}) in the above equation results in
	\begin{align*}
	\lVert p-p_{h}\rVert_{M}&\leq C \big(\lVert R_{0}p-p_{h}\rVert_{M}+\lVert E_h\bm{\phi}_{h}-\Pi_{h}(E_h\bm{\phi}_{h})\rVert_{Q}+\lVert R\mathbf{u}-\mathbf{u}_{h}\rVert_{h}+\lVert\bm{\phi}_{h}-\bar{R}\bm{\phi}\rVert_{h} \big).
	\end{align*}
This concludes the proof.
\end{proof}
\begin{theorem}\label{thms1} There holds,
	\begin{align*}
	\lVert r-r_{h}\rVert_{M}&\leq C \big(\lVert \bar{R}_{0}r-r_{h}\rVert_{M}+\lVert E_h\bm{\phi}_{h}-\Pi_{h}(E_h\bm{\phi}_{h})\rVert_{Q}+\lVert R\mathbf{u}-\mathbf{u}_{h}\rVert_{h}+\lVert\bm{\phi}_{h}-\bar{R}\bm{\phi}\rVert_{h} \big).
	\end{align*}
	\end{theorem}
\begin{proof}  The proof follows by the similar steps as in Theorem \ref{thm2.8}.
\end{proof}
	\begin{remark} Note that the upper bounds of the error terms in Theorem  \ref{thm2.7},  \ref{thm2.8} and \ref{thms1} contain some
	non-computable terms, for example, $\norm{R\mathbf{u}-\mathbf{u}_h}_h, \lVert\bm{\phi}_{h}-\bar{R}\bm{\phi}\rVert_{h}, \lVert R_{0}p-p_{h}\rVert_{M}$  and  $\lVert \bar{R}_{0}r-r_{h}\rVert_{M}$. These terms can be estimated in terms of the computable quantities so-called estimator, as it has been shown in  Section \ref{sec4}.
\end{remark}
\section{Model problems}\label{sec3}
This section deals with two model problems: A distributed control problem and a Neumann boundary control problem. We will see the application of the abstract framework from Section \ref{sec:Abstract Setting} to these model problems.
We start with some notation used throughout this article. Let $\Omega \subset \mathbb R^2$ be a  bounded polyhedral domain with boundary  $\Gamma$. The spaces $L^2(\Omega )$ and $H^1(\Omega)$ are standard Sobolev spaces.
The vector
valued version of $L^2(\Omega)$ and $H^1(\Omega)$ are denoted by  $[L^2(\Omega)]^2$ and $[H^1(\Omega)]^2$, respectively, and the $L^2$-norm on $\Omega$ is denoted by $\norm{\cdot}_{0,\Omega}$.
The subspace of $L^2(\Omega )$ with zero mean functions is defined by
$L_0^2(\Omega) := \{q \in L^2(\Omega) : \int_{\Omega}  q ~dx = 0\} $
and $H_0^1(\Omega)$ is the subspace of $H^1(\Omega)$ with zero trace functions.

\subsection{Distributed control problem}\label{sec4.2}
Set $$X=[H^1_{0}(\Omega)]^2,~W=[L^2(\Omega)]^2,~M=L^2_{0}(\Omega)~~{\rm and~~}  Q=[L^2(\Omega)]^2.$$
The map $E:X\rightarrow Q$ is the inclusion map.
 Given $\mathbf{u}_d$ in $[{L^{2}(\Omega)}]^2,$ define the quadratic functional $ J :[H_{0}^1(\Omega)]^2 \times [{L^{2}(\Omega)}]^2\rightarrow \mathbb{R}$ by
\begin{equation}\label{21}
J(\mathbf{v},\mathbf{x})= \frac{1}{2}\norm{\mathbf{v}-\mathbf{u}_d}_{0,\Omega}^2 + \frac{\lambda}{2}\norm{\mathbf{x}}_{0,\Omega}^2\;\; \mathbf{v}\in [H_{0}^1(\Omega)]^2,\; \mathbf{x}\in[L^{2}(\Omega)]^2.
\end{equation}
For given $\mathbf{y}_a,\mathbf{y}_b\in [\mathbb{R}\cup \{\pm \infty\}]^2$ with $\mathbf{y}_a<\mathbf{y}_b,$ define the admissible set of controls by
\begin{equation}\label{qb1}
 Q_{d} = \{\mathbf{x}\in [L^{2}(\Omega)]^2 |\; \mathbf{y}_a\leq \mathbf{x}\leq \mathbf{y}_b\},\;\; Q_{ad}:=Q_{d}.
\end{equation}
Consider the optimal control problem of finding $(\mathbf{u},\mathbf{y})\in X\times Q_d$ such that,
\begin{equation}\label{22}
J(\mathbf{u},\mathbf{y})=\inf_{\mathbf{v}\in [H_{0}^1(\Omega)]^2,~\mathbf{x}\in Q_{d} } J(\mathbf{v},\mathbf{x}),
\end{equation}
subject to the condition that $\mathbf{v}$ and $\mathbf{x}$ are such that $(\mathbf{v},p,\bf{x})$ satisfies:  for all $(\mathbf{z},w)\in X \times M$
\begin{align}\label{23}
a(\mathbf{v},\mathbf{z})+b(\mathbf{z},p) &={\langle \mathbf{x},\mathbf{z} \rangle}_{[L^{2}(\Omega)]^2} +{\langle \mathbf{f},\mathbf{z} \rangle}_{[L^{2}(\Omega)]^2} \;\;\;\text{and}\;\;\;
b(\mathbf{v},w)=0,
\end{align}
where $a(\mathbf{v},\mathbf{z})=  \int_{\Omega} \nabla{\mathbf{v}}:\nabla{\mathbf{z}} \dx$ , $b(\mathbf{z},p)= - \int_{\Omega} p\nabla\cdot{\mathbf{z}} \dx,$ and the matrix product $A:B := \sum_{i,j=1}^{n}a_{ij}b_{ij}$ when $A=(a_{ij})_{1\leq i,j \leq n}$ and $B=(b_{ij})_{1\leq i,j \leq n}.$

\smallskip
\noindent
The optimal solution $(\mathbf{u},p,\mathbf{y}) \in X\times M\times Q_{d},$ satisfy the following: for all $(\mathbf{z},w)\in X \times M$
\begin{align}\label{24}
a(\mathbf{u},\mathbf{z})+b(\mathbf{z},p) ={\langle \mathbf{y},\mathbf{z} \rangle}_{[L^{2}(\Omega)]^2} +{\langle \mathbf{f},\mathbf{z} \rangle}_{[L^{2}(\Omega)]^2} \;\; \text{and}\;\;
b(\mathbf{u},w)=0.
\end{align}
Note that the model problem (\ref{23}) has a unique solution $(\mathbf{v},p)\in [H_{0}^1(\Omega)]^2 \times L_{0}^2(\Omega) $ for given $\mathbf{x} \in [L^2(\Omega)]^2$ \cite[pp. 81]{Girault:1979}. We can set this correspondence as $S^{\mathbf{v}}\mathbf{x}=\mathbf{v}$ and $S^{p}\mathbf{x}=p$ and using the stability estimates of $\mathbf{v}$ and $p$ one can show that $S^{\mathbf{v}}:[L^2(\Omega)]^2\rightarrow [L^2(\Omega)]^2$ and $ S^{p}:[L^2(\Omega)]^2\rightarrow L_{0}^2(\Omega)$ are continuous affine operators. Then the minimization problem (\ref{22}) becomes \\
\begin{equation}\label{25}
j(\mathbf{y})=\inf_{\mathbf{x}\in Q_{d} } j(\mathbf{x}),
\end{equation}
where
\begin{equation}\label{26}
j(\mathbf{x})= \frac{1}{2}\norm{S^{\mathbf{v}}\mathbf{x}-\mathbf{u}_d}_{0,\Omega}^2 + \frac{\lambda}{2}\norm{\mathbf{x}}_{0,\Omega}^2.
\end{equation}
Using the theory of elliptic optimal control problems, the following proposition on the existence and uniqueness of the solution can be proved and the optimality condition can be derived \cite[pp. 50]{Troltzzsch:2005:book}.

\begin{proposition}\label{proposition 3.1}
	The control problem (\ref{25}) has a unique solution $\mathbf{y}$ and correspondingly there exist a unique state $(\mathbf{u},p)=(S^{\mathbf{u}}\mathbf{y},S^{p}\mathbf{y})$ of (\ref{24}). Furthermore, by introducing the adjoint state $(\bm{\phi},r)\in X \times M$ such that, ${\rm for~all}\;(\mathbf{z},w )\in X \times M$
	\begin{align}\label{27}
	a(\mathbf{z},\bm{\phi})-b(\mathbf{z},r)&={\langle \mathbf{u}-\mathbf{u}_d,\mathbf{z} \rangle}_{[{L^{2}(\Omega)}]^2}, \qquad
	b(\bm{\phi},w)=0,
	\end{align}
	the optimality condition that $j'(\mathbf{y})(\mathbf{x}-\mathbf{y})\geq 0,$  $\forall \mathbf{x}\in Q_{d},$ can be expressed as
	\begin{equation}\label{28}
	{\langle \bm{\phi}+\lambda \mathbf{y}, \mathbf{x}-\mathbf{y} \rangle }_{[L^2(\Omega)]^2}  \geq 0 \;\;\; {\rm for ~all}~~\mathbf{x}\in Q_{d}.
	\end{equation}
\end{proposition}

\vspace{.5cm}
The strong formulation of the optimality conditions satisfied by $((\mathbf{u},p),(\bm{\phi},r),\mathbf{y})$ is given by the following system of equations:
\begin{subequations}
\begin{align}
-\Delta \mathbf{u}+\nabla{p}&=\mathbf{y}+\mathbf{f},&\nabla\cdot{\mathbf{u}}&=0\;~~ \text{in}\; \Omega, &\bf{ u}&={\bf 0}\;\; \text{on} \;\;\partial \Omega,\label{mp1}\\
-\Delta \bm{\phi}-\nabla{r}&=\mathbf{u}-\mathbf{u}_d,
&\nabla\cdot{\bm{\phi}}&=0\;\; \text{in}\; \Omega,& \bm{\phi}&={\bf 0}\;\; \text{on} \;\;\partial\Omega,\\
\mathbf{y}&= {\Pi_{[\mathbf{y}_a,\mathbf{y}_b]}}\big(-\frac{\bm{\phi}}{\lambda}\big)\;\;\; \text{on}\; \Omega,\label{mp3}
\end{align}
\end{subequations}
where $\Pi_{[a,b]}g(x)=\text{min}\{b,\text{max}\{a,g(x)\}\}.$
\subsection{Neumann boundary control problem}\label{Model problem 2}
Set
$$X=[H^1(\Omega)\cap L_0^2(\Omega)]^2,~ W=[L^2(\Omega)]^2, ~M=L^2(\Omega)~ {\rm and} ~Q=[L^2(\Gamma)]^2.$$
The map $E:X\rightarrow Q$ is the trace map.
Define the quadratic functional $ J : X\times [{L_{2}(\Gamma)}]^2\rightarrow \mathbb{R}$ by\\
\begin{equation}\label{29}
J(\mathbf{v},\mathbf{x})= \frac{1}{2}\norm{\mathbf{v}-\mathbf{u}_d}_{0,\Omega}^2 + \frac{\lambda}{2}\norm{\mathbf{x}}_{0,\Gamma}^2\;\; \mathbf{v}\in X,\; \mathbf{x}\in[{L^{2}(\Gamma)}]^2.
\end{equation}
For given $\mathbf{y}_a,~\mathbf{y}_b\in [\mathbb{R}\cup \{\pm \infty\}]^2$ with $\mathbf{y}_a<\mathbf{y}_b,$ and $\mathbf{y}_a\leq -\frac{1}{|\Gamma|}\int_\Omega {\bf f}\dx\leq \mathbf{y}_b$ define the admissible set of controls by
\begin{equation}\label{qb}
 Q_{b} = \{\mathbf{x}\in [L_{2}(\Gamma)]^2 |\;\int_{\Gamma}{\mathbf{x}\ds}+\int_{\Omega}{\mathbf{f}\dx}=\mathbf{0},\;\; \mathbf{y}_a\leq \mathbf{x}\leq \mathbf{y}_b\},\;\; Q_{ad}:=Q_{d}.
\end{equation}
Consider the optimal control problem of finding $(\mathbf{u},\mathbf{y})\in X\times Q_b$ such that,
\begin{equation}\label{30}
J(\mathbf{u},\mathbf{y})=\inf_{\mathbf{v}\in X,~\mathbf{x}\in Q_{b} } J(\mathbf{v},\mathbf{x})
\end{equation}
subject to the condition that $\mathbf{v}$ and $\mathbf{x}$ are such that $(\mathbf{v},p, \mathbf{x})$ satisfies: for all $({\bf z},w)\in X \times M$
\begin{align}\label{31}
a(\mathbf{v},\mathbf{z})+b(\mathbf{z},p) &={\langle \mathbf{x},\mathbf{z} \rangle}_{[L^{2}(\Gamma)]^2} +{\langle \mathbf{f},\mathbf{z} \rangle}_{[L^{2}(\Omega)]^2},~~
b(\mathbf{v},w)=0.
\end{align}
where both the bilinear forms $a(\mathbf{v},\mathbf{z})$ and $b(\mathbf{z},p)$ are same as that of the distributed case.

\vspace{.5cm}
\noindent
The optimal solution $(\mathbf{u},p,\mathbf{y}) \in X\times M \times Q_{b}$ satisfies
\begin{align}\label{32}
a(\mathbf{u},\mathbf{z})+b(\mathbf{z},p) &={\langle \mathbf{y},\mathbf{z} \rangle}_{[L_{2}(\Gamma)]^2} +{\langle \mathbf{f},\mathbf{z} \rangle}_{[{L^{2}(\Omega)}]^2} \;{\rm for~all}\;\mathbf{z} \in X, \nonumber\\
b(\mathbf{u},w)&=0\;{\rm for~all}\;w \in M.
\end{align}
Note that the model problem (\ref{31}) has a unique solution $(\mathbf{v},p)\in X \times M $ for given $\mathbf{x} \in [L^2(\Gamma)]^2$ \cite{Girault:1979}. We can set this correspondence as $S^{\mathbf{v}}\mathbf{x}=\mathbf{v}$ and $S^{p}\mathbf{x}=p$ and using the stability estimates of $\mathbf{v}$ and $p,$ one can show that $S^{\mathbf{v}}:[L^2(\Gamma)]^2\rightarrow [L^2(\Omega)]^2$ and $ S^{p}:[L^2(\Gamma)]^2\rightarrow L^2(\Omega)$ are continuous affine operators. Then the minimization problem (\ref{30}) becomes
\begin{equation}\label{33}
j(\mathbf{y})=\inf_{\mathbf{x}\in Q_{b} } j(\mathbf{x}),
\end{equation}
where
\begin{equation}\label{34}
j(\mathbf{x})= \frac{1}{2}\norm{S^{\mathbf{v}}\mathbf{x}-\mathbf{u}_d}^2_{0,\Omega} + \frac{\lambda}{2}\norm{\mathbf{x}}^2_{0,\Gamma}.
\end{equation}
 The following proposition on the existence and uniqueness of the solution  and the optimality condition can be derived using the theory of elliptic optimal control problems \cite[pp. 50]{Troltzzsch:2005:book}.

\begin{proposition}\label{proposition 3.2}
	The control problem (\ref{33}) has a unique solution $\mathbf{y}$ and correspondingly there exist a unique state $(\mathbf{u},p)=(S^{\mathbf{v}}\mathbf{y},S^{p}\mathbf{y})$ of (\ref{32}). Furthermore, by introducing the adjoint state $(\bm{\phi},r)\in X\times M$ such that, for all $(\mathbf{z},w) \in X \times M $
	\begin{align}\label{35}
	a(\mathbf{z},\mathbf{u})-b(\mathbf{z},r) &={\langle \mathbf{u}-\mathbf{u}_d,\mathbf{z} \rangle}_{[L^{2}(\Omega)]^2}, \qquad
	b(\bm{\phi},w)=0,
	\end{align}
	the optimality condition that $j'(\mathbf{y})(\mathbf{x}-\mathbf{y})\geq 0,$  ${\rm for~all} ~\mathbf{x}\in Q_{b},$ can be expressed as
	\begin{equation}\label{36}
	{\langle \bm{\phi}+\lambda \mathbf{y}, \mathbf{x}-\mathbf{y} \rangle }_{[L_{2}(\Gamma)]^2} \geq 0 \; \;{\rm for~all} ~\mathbf{x}\in Q_{b}.
	\end{equation}
\end{proposition}
The strong formulation of the optimality conditions satisfied by $((\mathbf{u},p),(\bm{\phi},r),\mathbf{y})$ is given by the following system of equations:
\begin{subequations}
	\begin{align}
	-\Delta \mathbf{u}+\nabla{p}&=\mathbf{f},\;\;\qquad\;\; \nabla\cdot{\mathbf{u}}=0\;\; \text{in}\;\Omega,\qquad\frac{\partial \mathbf{u}}{\partial \mathbf{n}}-p \mathbf{n}=\mathbf{y}\;\; \text{on} \;\;\partial\Omega,\\
	-\Delta \bm{\phi}-\nabla{r}&=\mathbf{u}-\mathbf{u}_d,\;\; \nabla\cdot{\bm{\phi}}=0\;\; \text{in}\; \Omega,\qquad \frac{\partial \bm{\phi}}{\partial \mathbf{n}}+r \mathbf{n}=0\;\; \text{on} \;\;\partial\Omega,\\
	\mathbf{y}&={\Pi_{[\mathbf{y}_a,\mathbf{y}_b]}}\big(-\frac{\bm{\phi}}{\lambda}\big)\;\; \text{on}\; \partial\Omega,
\end{align}
	\end{subequations}
where $\Pi_{[a,b]}g(x)=\text{min}\{b,\text{max}\{a,g(x)\}\}.$

\section{Discrete Problems}\label{sec4}
In this section, we will discuss finite element formulations for the model problems studied in the last section.  First, we will consider
Crouzeix-Raviart finite element space for velocity and piecewise constant polynomial space for pressure.  Secondly, we will discuss
the discontinuous Galerkin (DG) method with piecewise linear space and piecewise constant space  for velocity and pressure approximation, respectively. We start this section with the notation which has used throughout  the article.
\subsection{Notation}\label{notations 4.1}
\begin{figure}[!hh]
\begin{center}
\setlength{\unitlength}{0.7cm}
\begin{picture}(8,6)
\put(2,1){\line(1,0){8}} \put(2,1){\line(-1,1){4}}
\put(6,5){\line(-1,0){8}} \put(10,1){\line(-1,1){4}}
\thicklines\put(2,1){\vector(1,1){4}}\thinlines
\put(1.39,0.39){$A$} \put(5.95,5.25){$B$} \put(10,0.29){$P_{+}$}
\put(-2.75,4.85){$P_{-}$} \put(1.25,3.25){$T_{-}$}
\put(6,2.5){$T_{+}$} \put(4,3){\vector(1,-1){1}}
\put(5.10,1.90){$\mathbf{n}_{e}$}
\thicklines\put(4,3){\vector(1,1){1}}\thinlines
\put(5.10,3.50){$\tau_{e}$} \put(3.25,2.75){$e$}
\qbezier(4.2,2.8)(4.4,3.0)(4.2,3.2) \put(4.1,2.95){.}
\end{picture}
\caption{\footnotesize{Here $T_-$ and $T_{+}$ are the two neighboring triangles that share the edge $e=\partial T_-\cap\partial T_+$ with initial node $A$
and end node $B$ and unit normal $\mathbf{n}_e$. The orientation of $\mathbf{n}_e = \mathbf{n}_{-} = -\mathbf{n}_{+}$ equals the outer normal of $T_{-}$, and hence, points into $T_{+}$.}} \label{Fig1}
\end{center}
\end{figure}
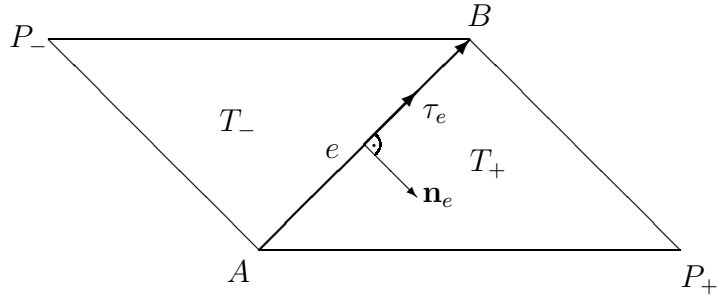

 Let
 $\mathcal{T}_{h}$ be a regular triangulation of $\Omega$ into  triangles such that $\cup_{T\in \mathcal T_h} T=\bar\Omega$. Denote the set of all interior
edges of $\mathcal{T}_{h}$ by $\mathcal{E}^i_h$, the set of boundary edges by $\mathcal{E}^b_{h}$, and define
$\mathcal{E}_h = \mathcal{E}_{h}^i \cup \mathcal{E}_h ^b$. Let
$h_{T}$:=diam($T$) and $h=\text{max} \{h_{T} : T\in \mathcal{T}_h\}$. The length of
any edge $e \in \mathcal{E}_h$ will be denoted by $h_{e}$.
Let us define a broken Sobolev space
$$H^1(\Omega,\mathcal{T}_{h})=\{v\in L^2(\Omega): v|_{T}\in H^1(T) \;\;{\rm for~all~} T\in \mathcal{T}_h \}.$$
Denote the norm and the semi-norm on $H^k(D)$ $(k\geq0)$ for any
domain $D\subseteq \mathbb{R}^2$  by $\norm{v}_{k,D}$ and $|v|_{k,D}.$
In the problem setting, we require jump and mean definitions of discontinuous functions, vector functions and tensors. For any $e \in \mathcal{E}^{i}_{h}$, there are two triangles $T_{+}$ and $T_{-}$ such that $e=\partial T_{+} \cap \partial T_{-}$. Let
$\mathbf{n}_{+}$ be the unit normal of $e$ pointing from  $T_{+}$ to $T_{-}$ and let $\mathbf{n}_{-} =-\mathbf{n}_{+}$ (cf. Fig.\ref{Fig1}). For any $v \in H^1(\Omega,\mathcal{T}_{h})$, we define
the jump and mean of $v$ on an edge $e$ by
$$\sjump{v}=v_+\mathbf{n}_+ +v_-\mathbf{n}_-\;\; \text{and} \;\;\smean{v}= \frac{1}{2}(v_++v_-) \;\; \text{respectively},$$
where $v_{\pm}=v|_{T_{\pm}}.$

\noindent
For $\mathbf{v} \in [H^1(\Omega,\mathcal{T}_{h})]^2$ we define the jump and mean of $\mathbf{v}$ on $e \in \mathcal{E}^{i}_{h}.$ by
$$\sjump{\mathbf{v}}=\mathbf{v}_+ \cdot \mathbf{n}_+ +\mathbf{v}_-\cdot \mathbf{n}_-\;\; \text{and} \;\;\smean{\mathbf{v}}= \frac{1}{2}(\mathbf{v}_++\mathbf{v}_-) \;\; \text{respectively}.$$
\noindent
We also require the full jump of vector valued functions. For $\mathbf{v} \in [H^1(\Omega,\mathcal{T}_{h})]^2$, we define the full
jump by
$$\underline{\sjump{\mathbf{v}}}=\mathbf{v}_+ \otimes \mathbf{n}_+ +\mathbf{v}_- \otimes \mathbf{n}_- ,$$
where for two vectors in Cartesian coordinates $\mathbf{a}=(a_i)$ and $\mathbf{b}=(b_j)$, we define the matrix
$\mathbf{a}\otimes \mathbf{b}=[a_ib_j]_{1\leq i,j\leq 2}$ . Similarly, for tensors $\tau  \in [H^1(\Omega,\mathcal{T}_{h})]^{2\times 2}$, the jump and mean on $ e\in \mathcal{E}^i_{h}$ are defined by
$$\sjump{\tau}=\tau_+\mathbf{n}_+ +\tau_-\mathbf{n}_-\;\; \text{and} \;\;\smean{\tau}= \frac{1}{2}(\tau_++\tau_-).$$

\noindent
For notational convenience, we also define the jump and mean on the boundary faces $e\in \mathcal{E}_{h}^b$ by modifying them appropriately. We use the definition of jump by understanding that $v_{-} = 0$ (similarly, $\mathbf{v}_{-} = 0$
and $\tau_{-} = 0$) and the definition of mean by understanding that $v_{-} = v _{+}$ (similarly, $\mathbf{v}_{-} = \mathbf{v}_{+}$  and $\tau_{-} = \tau_{+}$).

 Define $\mathbb{P}_m(T)$ to be the space of polynomials of degree at most  $m (\geq 0)$ defined on the triangle $T$. The discontinuous finite element space is
 $$ \mathbb{P}_m(\mathcal T_h):= \{ v_h\in L^2(\Omega)| ~\forall T\in \mathcal T_h, v_h|_T\in  \mathbb{P}_m(T)\}.$$
The  lowest-order Crouzeix-Raviart (CR) spaces are  defined as
\begin{align*}
{\rm CR}^1(\mathcal T_h):=\{\mathbf{v}_{h}\in [L^2(\Omega)]^2: &~ \forall T\in \mathcal{T}_{h}, \mathbf{v}_{h}|_{T}\in [\mathbb{P}_{1}(T)]^2,\\
& \;~v_{h}
\text{\;continuous at midpoint of edge}\; e \; ~{\rm for~all~} e \in \mathcal{E}^i_{h}\}.
\end{align*}
and
 ${\rm CR}^1_{0}(\mathcal T_h):=\{\mathbf{v}_{h}\in {\rm CR}^1(\mathcal T_h): \mathbf{v}_h({\rm mid}~ e)={\bf 0}~ {\rm for~all~} e \in \mathcal{E}^b_{h}\}.$\\
\noindent
 Define the oscillation of given functions $\mathbf{f},~ \mathbf{u}_d\in [L^2(T)]^2$ by
 $$osc(\mathbf{f},T)=h_{T}\min_{\mathbf{f}_h\in [\mathbb{P}_{0}(T)]^2 }\lVert \mathbf{f}-\mathbf{f}_h\rVert_{0,T}\;\; \text{and}\;\;osc(\mathbf{u}_d,T)=h_{T}\min_{\mathbf{g}_h\in [\mathbb{P}_{0}(T)]^2}\lVert \mathbf{u}_d-\mathbf{g}_h\rVert_{0,T}.$$

We will also need the following inverse estimates \cite{BScott:2008:FEM}:
\begin{lemma}
There exist a constant $C_m$ such that for all $v\in \mathbb{P}_m(T)$ one has
\begin{equation*}
\norm{\nabla v}_{0,T}\leq C_m h_T^{-1}\norm{v}_{0,T}
\end{equation*}	
and
\begin{equation}
\norm{v}_{0,\partial T}\leq C_m h_T^{-{\frac{1}{2}}}\norm{v}_{0,T}.
\end{equation}
\end{lemma}

\smallskip

\subsection{Nonconforming FEM (${\rm CR}_{0}^1(\mathcal T_h)$/$\mathbb P_0(\mathcal T_h)$)} This subsection is devoted to the \textit{a priori} and \textit{a posteriori} error analysis for the pair Crouzeix-Raviart finite element space $({\rm CR}_{0}^1(\cT_h))$ for velocity and piecewise constant polynomial space $\mathbb P_0(\mathcal T_h)$ for pressure.
\subsubsection{\bf Discrete distributed control problem:}
Consider the model problem from Subsection \ref{sec4.2}. Set $$X=[H^1_{0}(\Omega)]^2,~W=[L^2(\Omega)]^2,~M=L^2_{0}(\Omega)~ {\rm and} ~Q=[L^2(\Omega)]^2.$$
The set $Q_{ad}=Q_{d}$, where $Q_{d}$ is defined in (\ref{qb1}). The discrete spaces are defined by
$$X_{h}={\rm CR}_{0}^1(\mathcal T_h),~~
M_{h}:=\mathbb P_0(\mathcal T_h)\cap  L_{0}^2(\Omega)~{\rm and}~~ Q_{h}=\mathbb P_0(\mathcal T_h).$$
The admissible control set is $Q_{ad}^h=\{\mathbf{x}_h\in Q_{h}: \mathbf{y}_a\leq \mathbf{x}_h\leq \mathbf{y}_b\}.$ It is clear that $Q_{ad}^h \subseteq Q_{ad}$ and $\Pi_{h}\mathbf{u}\in Q_{ad}^h$ for $\mathbf{u}\in Q_{ad}.$ The operators $E: X\rightarrow Q$ and $E_h:X_{h}\rightarrow Q_{h}$ are inclusion maps.
The bilinear forms  for the diffusion term  and  pressure term are given by for all  $\mathbf{v}_{h}$,~$\mathbf{z}_{h}\in X_{h}$
\begin{equation}\label{52}
a_{h}(\mathbf{v}_{h},\mathbf{z}_{h})= \sum_{T\in \mathcal{T}_{h}} \int_{T} \nabla{\mathbf{v}_{h}}:\nabla{\mathbf{z}_{h}} \dx,\qquad b_{h}(\mathbf{z}_{h},p_{h})= -\sum_{T\in \mathcal{T}_{h}} \int_{T} p_{h}\nabla\cdot{\mathbf{z}_{h}} \dx.
\end{equation}
\noindent
Energy norm on $X_{h}$ is defined by $$\lVert \mathbf{v}_{h}\rVert_{h}^2=\sum_{T\in \mathcal{T}_{h}} \int_{T} \rvert\nabla{\mathbf{v}_{h}}\lvert^2 \dx,$$
and $b_{h}(\cdot,\cdot)$ satisfies the inf-sup condition can be found in \cite{Ern:2012:DGBook}. Assumptions (\ref{Assmption1}) and (\ref{Assmption2}) are the Poincar\'e-Friedrichs  type inequality in \cite[pp. 301]{BScott:2008:FEM}.

\begin{theorem}\label{thm 4.9}
{\bf (Best approximation for velocity and control)} It holds
		\begin{align*}
	\lVert   \mathbf{y}-\mathbf{y}_h \rVert_{0,\Omega}+&\lVert \bm{\phi}-\bm{\phi}_{h}\rVert_{h}+ \lVert \mathbf{u}-\mathbf{u}_{h}\rVert_{h} \leq C \big(\inf_{\mathbf{w}_h\in Z_{h}}\lVert \mathbf{u}-\mathbf{w}_h\rVert_{h}+\inf_{\mathbf{w}_{h}\in Z_{h}}\lVert\bm{\phi}-\mathbf{w}_{h}\rVert_{h}\\&+\inf_{s_{h}\in M_{h}}\lVert p-s_{h}\rVert_{0,\Omega}+\inf_{s_{h}\in M_{h}}\lVert r-s_{h}\rVert_{0,\Omega}+h\lVert \mathbf{f} \rVert_{0,\Omega}+ h\inf_{\mathbf{x}_h\in Q_{h}}\lVert \mathbf{y}-\mathbf{x}_h\rVert_{0,\Omega}\\
	& + h\lVert \mathbf{u}-\mathbf{u}_d\rVert_{0,\Omega} +\norm{\bm{\phi}-\Pi_{h}\bm{\phi}}_{0,\Omega} +\lVert \mathbf{y}-\Pi_{h}\mathbf{y}\rVert_{0,\Omega} \big).
	\end{align*}
	\end{theorem}
\begin{proof}
	From Theorem \ref{thm 2.3} we know that
	$$ \lVert \bm{\phi}-\bm{\phi}_{h}\rVert_{h}\leq C \big(\lVert \bm{\phi}-\bar{P_{h}}\bm{\phi}\rVert_{h}+\norm{\bm{\phi}-\Pi_{h}\bm{\phi}}_{0,\Omega} +\lVert \mathbf{y}-\Pi_{h}\mathbf{y}\rVert_{0,\Omega}+\lVert \mathbf{u}-P_{h}\mathbf{u} \rVert_{h}\big),$$
	$$\lVert \mathbf{u}-\mathbf{u}_{h}\rVert_{h} \leq C \big(\lVert \mathbf{u}-P_{h}\mathbf{u} \rVert_{h}+\norm{\bm{\phi}-\Pi_{h}\bm{\phi}}_{0,\Omega} +\lVert \mathbf{y}-\Pi_{h}\mathbf{y}\rVert_{0,\Omega}+\lVert \bm{\phi}-\bar{P_{h}}\bm{\phi}\rVert_{h}\big).$$
	Also from the best approximation results \cite[Theorem 3.1]{BCGG:2014:DGStokes}, we have
	\begin{equation}\label{54}
	\lVert\bm{\phi}-\bar{P_{h}}\bm{\phi}\rVert_{h} \leq C\big( \inf_{\mathbf{w}_{h}\in Z_{h}}\lVert\bm{\phi}-\mathbf{w}_{h}\rVert_{h}+\inf_{s_{h}\in M_{h}}\lVert r-s_{h}\rVert_{0,\Omega}+ h\lVert \mathbf{u}-\mathbf{u}_d\rVert_{0,\Omega}\big)
	\end{equation}
	and
	\begin{equation}\label{55}
	\lVert \mathbf{u}-P_{h}\mathbf{u}\rVert_{h} \leq C\big( \inf_{\mathbf{v}_{h}\in Z_{h}}\lVert \mathbf{u}-\mathbf{v}_h\rVert_{h}+\inf_{s_{h}\in M_{h}}\lVert p-s_{h}\rVert_{0,\Omega}+h\lVert \mathbf{f} \rVert_{0,\Omega}+ h\inf_{\mathbf{x}_h\in Q_{h}}\lVert \mathbf{y}-\mathbf{x}_h\rVert_{0,\Omega}\big).
	\end{equation}
 Theorem \ref{thm2.2} shows that
	$$\lVert   \mathbf{y}-\mathbf{y}_h \rVert_{0,\Omega} \leq C\big(\norm{\bm{\phi}-\Pi_{h}\bm{\phi}}_{0,\Omega} +\lVert \mathbf{y}-\Pi_{h}\mathbf{y}\rVert_{0,\Omega}+\lVert \bm{\phi}-\bar{P_{h}}\bm{\phi}\rVert_{h}+\lVert \mathbf{u}-P_{h}\mathbf{u} \rVert_{h}\big).$$
	Using the above estimates, we find
	\begin{align*}
	\lVert   \mathbf{y}-\mathbf{y}_h \rVert_{0,\Omega}+&\lVert \bm{\phi}-\bm{\phi}_{h}\rVert_{h}+ \lVert \mathbf{u}-\mathbf{u}_{h}\rVert_{h} \leq C \big(\inf_{\mathbf{w}_h\in Z_{h}}\lVert \mathbf{u}-\mathbf{w}_h\rVert_{h}+\inf_{\mathbf{w}_{h}\in Z_{h}}\lVert\bm{\phi}-\mathbf{w}_{h}\rVert_{h}\\&+\inf_{s_{h}\in M_{h}}\lVert p-s_{h}\rVert_{0,\Omega}+\inf_{s_{h}\in M_{h}}\lVert r-s_{h}\rVert_{0,\Omega}+h\lVert \mathbf{f} \rVert_{0,\Omega}+ h\inf_{\mathbf{x}_h\in Q_{h}}\lVert \mathbf{y}-\mathbf{x}_h\rVert_{0,\Omega}\\
	& + h\lVert \mathbf{u}-\mathbf{u}_d\rVert_{0,\Omega} +\norm{\bm{\phi}-\Pi_{h}\bm{\phi}}_{0,\Omega} +\lVert \mathbf{y}-\Pi_{h}\mathbf{y}\rVert_{0,\Omega} \big).
	\end{align*}
This concludes the proof. \end{proof}
	
\begin{theorem}\label{thm 4.99}
{\bf (Best approximation for pressure)} It holds
	\begin{align}
\lVert p-p_{h}\rVert_{0,\Omega}\leq C& \big(\inf_{\mathbf{v}_{h}\in Z_{h}}\lVert \mathbf{u}-\mathbf{v}_{h}\rVert_{h}+\inf_{\mathbf{w}_{h}\in Z_{h}}\lVert\bm{\phi}-\mathbf{w}_{h}\rVert_{h}+\inf_{s_{h}\in M_{h}}\lVert p-s_{h}\rVert_{0,\Omega}\nonumber\\
& +\inf_{s_{h}\in M_{h}}\lVert r-s_{h}\rVert_{0,\Omega}+h\lVert \mathbf{f} \rVert_{0,\Omega}+h\inf_{\mathbf{x}_h\in Q_{h}}\lVert \mathbf{y}-\mathbf{x}_h\rVert_{0,\Omega} + h\lVert \mathbf{u}-\mathbf{u}_d\rVert_{0,\Omega}\nonumber\\
&+\lVert p-P^{0}p\rVert_{0,\Omega}+\norm{\bm{\phi}-\Pi_{h}\bm{\phi}}_{0,\Omega} +\lVert \mathbf{y}-\Pi_{h}\mathbf{y}\rVert_{0,\Omega} \big).\label{bestp}
\end{align}
\end{theorem}
\begin{proof}
	From Theorem \ref{thm 2.4}, the pressure estimate is given by
	\begin{align*}
	\lVert p-p_{h}\rVert_{0,\Omega}\leq C\big(\lVert p-R_{h}p\rVert_{0,\Omega}+ \lVert \mathbf{u}-P_{h}\mathbf{u} \rVert_{h}&+\lVert \bm{\phi}-\bar{P_{h}}\bm{\phi}\rVert_{h}\\&+\norm{\bm{\phi}-\Pi_{h}\bm{\phi}}_{0,\Omega} +\lVert \mathbf{y}-\Pi_{h}\mathbf{y}\rVert_{0,\Omega}\big).
	\end{align*}
	From the best approximation result \cite[Theorem 4.1]{BCGG:2014:DGStokes} for pressure we have
	\begin{align*}
	\lVert p-R_{h}p\rVert_{0,\Omega} \leq C \big(\inf_{\mathbf{v}_h\in Z_{h}}\lVert \mathbf{u}-\mathbf{v}_{h}\rVert_{h}+\inf_{q_{h}\in M_{h}}\lVert p-q_{h}\rVert_{0,\Omega}&+\lVert p-P^{0}p\rVert_{0,\Omega}\\
	&+ h\inf_{\mathbf{x}_h\in Q_{h}}\lVert \mathbf{y}-\mathbf{x}_h\rVert_{0,\Omega}+h\lVert \mathbf{f} \rVert_{0,\Omega}\big),
	\end{align*}
	where $P^{0}$ is the $L^2$-projection onto the space of piecewise constant polynomials. A	
	use the above estimates with (\ref{54}) and (\ref{55}) leads to (\ref{bestp}). This concludes the proof.
\end{proof}	

	Now we use elliptic regularity to derive concrete error estimates. Note that by well-posedness of the problem, $\mathbf{u},\bm{\phi}\in [H_{0}^1(\Omega)]^2$ and $p\in L^2_{0}(\Omega)$. The elliptic regularity on polygonal domains implies that $\bm{\phi} \in [H^{1+s}(\Omega) \cap H_{0}^1(\Omega)]^2$ and $r \in H^s(\Omega) \cap L^2_{0}(\Omega)$  for some $ s\in (0,1]$, which depends on the interior angles of the domain $\Omega.$ We know that,
$\mathbf{y}={\Pi_{[\mathbf{y}_a,\mathbf{y}_b]}}\big(-\frac{\bm{\phi}}{\lambda}\big)\; \text{on}\; \Omega.$
Hence, the control $\mathbf{y} \in H^{1}(\Omega)$ and $\mathbf{u}\in [H^{1+s}(\Omega)\cap H_{0}^1(\Omega)]^2$ and $p \in H^s(\Omega) \cap L^2_{0}(\Omega)$.
Using the N\'ed\'elec interpolation \cite[Eq. 3.4]{BCGG:2014:DGStokes},~\cite{nedelec} we have the following estimates
\begin{equation}\label{65}
\inf_{\mathbf{v}_{h}\in Z_{h}}\lVert \mathbf{u}-\mathbf{v}_{h}\rVert_{h}\leq Ch^{s}\lVert \mathbf{u} \rVert_{1+s}\;\;\text{and}\;\;
\inf_{\mathbf{w}_{h}\in Z_{h}}\lVert\bm{\phi}-\mathbf{w}_{h}\rVert_{h}\leq C h^{s}\lVert \bm{\phi} \rVert_{1+s}
\end{equation}
and 
\begin{align}\label{66}
\inf_{s_{h}\in M_{h}}\lVert p-s_{h}\rVert_{0,\Omega}\leq C h^{s}\lVert p \rVert_{s}\;\;\text{and}\;\;\;
\inf_{s_{h}\in M_{h}}\lVert r-s_{h}\rVert_{0,\Omega}\leq C h^{s}\lVert r \rVert_{s}.
\end{align}
Also we have the following estimates,
\begin{equation}\label{69}
\norm{\bm{\phi}-\Pi_{h}\bm{\phi}}_{0,\Omega}\leq C h \lVert \bm{\phi}\rVert_{1+s,\Omega}, ~~~~~~~~ \text{and}~~~~~~~~~
\norm{\mathbf{y}-\Pi_{h}\mathbf{y}}_{0,\Omega}\leq C h \lVert \mathbf{y}\rVert_{1,\Omega}.
\end{equation}

\noindent
Substitution of the above approximation estimates in Theorem \ref{thm 4.9} and \ref{thm 4.99} results in


\begin{theorem}
	Let $s\in (0,1]$ be the elliptic regularity index. Then, it holds
	\begin{align*}
	\lVert   \mathbf{y}-\mathbf{y}_h \rVert_{0,\Omega}+\lVert \bm{\phi}-\bm{\phi}_{h}\rVert_{h}+& \lVert \mathbf{u}-\mathbf{u}_{h}\rVert_{h} \leq C \big( h^s \lVert \mathbf{u}\rVert_{1+s,\Omega}+h^s \lVert \bm{\phi}\rVert_{1+s,\Omega}+h^s \lVert p\rVert_{s,\Omega}+h^s \lVert r\rVert_{s,\Omega}\\&+h^2\lVert \mathbf{y}\rVert_{1,\Omega} +h\lVert \mathbf{f} \rVert_{0,\Omega}+h\lVert \mathbf{u}-\mathbf{u}_d \rVert_{0,\Omega}+ h\lVert \bm{\phi}\rVert_{1+s,\Omega}+ h \lVert \mathbf{y}\rVert_{1,\Omega} \big).
	\end{align*}
\end{theorem}

\begin{theorem}
	Let $s\in (0,1]$ be the elliptic regularity index. Then there holds
	\begin{align*}
	\lVert p-p_{h}\rVert_{0,\Omega} \leq& C \big( h^s \lVert \mathbf{u}\rVert_{1+s,\Omega}+h^s \lVert \bm{\phi}\rVert_{1+s,\Omega}+h^s \lVert p\rVert_{s,\Omega}+h^s \lVert r\rVert_{s,\Omega}+h^2 \lVert \mathbf{y}\rVert_{1,\Omega}\\
	& +h\lVert \mathbf{f} \rVert_{0,\Omega}+h\lVert \mathbf{u}-\mathbf{u}_d \rVert_{0,\Omega}+ h\lVert \bm{\phi}\rVert_{1+s,\Omega}+ h \lVert \mathbf{y}\rVert_{1,\Omega} \big).
	\end{align*}
\end{theorem}

The following theorems deduce the reliable and efficient a posteriori error estimator.

\begin{theorem}\label{thm 4.11}
{\bf (A~posteriori error estimator) } It holds,
	\begin{align}\label{terror}
	\lVert \mathbf{y}-\mathbf{y}_h\rVert_{0,\Omega}+\lVert \mathbf{u}-\mathbf{u}_{h}\rVert_{h} +\lVert p-p_{h}\rVert_{0,\Omega}&+\lVert \bm{\phi}-\bm{\phi}_{h}\rVert_{h}+\lVert r-r_{h}\rVert_{0,\Omega}\nonumber\\
	&\leq C \big( \eta_{(\mathbf{u},p)}+\eta_{(\bm{\phi},r)}+\lVert\bm{\phi}_{h}-\Pi_{h}\bm{\phi}_{h}\rVert_{0,\Omega}\big)
	\end{align}
	where the estimators are defined as
	\begin{align*}
	\eta_{(\mathbf{u},p)}^2&=\sum_{T\in \mathcal{T}_{h}}\big( h_{T}^2 \lVert \mathbf{f}+\mathbf{y}_h\rVert_{0,T}^2\big)+\sum_{e\in \mathcal{E}_{h}^{i}} \lVert h^{\frac{1}{2}}\sjump{p_{h}I-\nabla{\mathbf{u}_{h}}}\rVert_{0,e}^2+ \sum_{e\in \mathcal{E}_{h}} \lVert h^{-\frac{1}{2}}\underline{\sjump{\mathbf{u}_{h}}}\rVert_{0,e}^2,
	\end{align*}
	and
	\begin{align*}
	\eta_{(\bm{\phi},r)}^2&=\sum_{T\in \mathcal{T}_{h}}\big( h_{T}^2 \lVert \mathbf{u}_{h}-\mathbf{u}_d\rVert_{0,T}^2\big)+\sum_{e\in \mathcal{E}_{h}^{i}} \lVert h^{\frac{1}{2}}\sjump{r_{h}I+\nabla{\bm{\phi}_{h}}}\rVert_{0,e}^2+ \sum_{e\in \mathcal{E}_{h}} \lVert h^{-\frac{1}{2}}\underline{\sjump{\bm{\phi}_{h}}}\rVert_{0,e}^2.
	\end{align*}
\end{theorem}
\begin{proof}From Theorem \ref{thm2.7},~ \ref{thm2.8} and \ref{thms1}, we have
	\begin{align}
	\lVert \mathbf{y}-\mathbf{y}_h\rVert_{0,\Omega}+&\lVert \mathbf{u}-\mathbf{u}_{h}\rVert_{h} +\lVert p-p_{h}\rVert_{0,\Omega}+\lVert \bm{\phi}-\bm{\phi}_{h}\rVert_{h}+\lVert r-r_{h}\rVert_{0,\Omega}\leq C \big(\lVert R\mathbf{u}-\mathbf{u}_{h}\rVert_{h}\nonumber\\
	&+\lVert R_{0}p-p_{h}\rVert_{0,\Omega}
	+\norm{\bm{\phi}_{h}-\bar{R}\bm{\phi}}_{h}+\lVert \bar{R}_{0}r-r_{h}\rVert_{0,\Omega}+\lVert\bm{\phi}_{h}-\Pi_{h}\bm{\phi}_{h}\rVert_{0,\Omega}\big).\label{111a}
	\end{align}
	The {\it a~posteriori} error analysis in  \cite[Theorem 3.1]{HSW:2005:Err.Stokes}, \cite[Section 7]{CarstensenEigelLoebhardetal.2010} gives the following error estimates:	
	\begin{equation}\label{56}
	\lVert R\mathbf{u}-\mathbf{u}_{h} \rVert_{h}+\lVert R_{0}p-p_{h}\rVert_{0,\Omega}\leq C \eta_{(\mathbf{u},p)},
	\end{equation}
	\begin{equation}\label{57}
	\lVert\bar{R} \bm{\phi}-\bm{\phi}_{h} \rVert_{h}+\lVert \bar{R_{0}}r-r_{h}\rVert_{0,\Omega}\leq C \eta_{(\bm{\phi},r)}.
	\end{equation}
By substituting (\ref{56}) and (\ref{57}) in (\ref{111a}), we conclude the proof.
\end{proof}

\begin{theorem}{\bf (Efficiency)}\label{4.12}
	Let $\mathcal{T}_{e}$ be the set of two triangles sharing the edge $e\in \mathcal{E}_{h}^{i}.$ Then, It hold
	\begin{align*}
	h_{T}\lVert \mathbf{f}+\mathbf{y}_h\rVert_{0,T} &\leq C\big(\lVert \nabla(\mathbf{u}-\mathbf{u}_{h})\rVert_{0,T}+\lVert p-p_{h}\rVert_{0,T}+ osc(\mathbf{f},T)+\lVert \mathbf{y}-\mathbf{y}_h\rVert_{0,T}\big),\\
	h_{T} \lVert \mathbf{u}_{h}-\mathbf{u}_d\rVert_{0,T}&\leq C\big(\lVert \mathbf{u}-\mathbf{u}_{h}\rVert_{0,T}+\lVert \nabla( \bm{\phi}-\bm{\phi}_{h})\rVert_{0,T}+\lVert r-r_{h}\rVert_{0,T}+ osc(\mathbf{u}_d,T)\big),\\
	\lVert h^{\frac{1}{2}}\sjump{p_{h}I-\nabla{\mathbf{u}_{h}}}\rVert_{0,e}&\leq C \sum_{T\in \mathcal{T}_{e}}\big(\lVert \nabla(\mathbf{u}-\mathbf{u}_{h})\rVert_{0,T}+\lVert p-p_{h}\rVert_{0,T}+ osc(\mathbf{f},T)+\lVert \mathbf{y}-\mathbf{y}_h\rVert_{0,T}\big),\\
	\lVert h^{\frac{1}{2}}\sjump{r_{h}I+\nabla{\bm{\phi}_{h}}}\rVert_{0,e}&\leq C \sum_{T\in\mathcal{T}_{e}}\big(\lVert \mathbf{u}-\mathbf{u}_{h}\rVert_{0,T}+\lVert \nabla( \bm{\phi}-\bm{\phi}_{h})\rVert_{0,T}+\lVert r-r_{h}\rVert_{0,T}+ osc(\mathbf{u}_d,T)\big),\\
	\lVert\bm{\phi}_{h}-\Pi_{h}\bm{\phi}_{h}\rVert_{0,T}&\leq C \big(\lVert\bm{\phi}_{h}-\bm{\phi}\rVert_{0,T}+\lVert\bm{\phi}-\Pi_{h}\bm{\phi}_{h}\rVert_{0,T}\big).
	\end{align*}
\end{theorem}
\begin{proof}
The local efficiency can be deduced by the standard bubble function techniques \cite{Verfurth 1995}.
\end{proof}

\subsubsection{\bf Discrete boundary control problem:}
The model problem in this section is the Neumann boundary control problem  introduced in Subsection \ref{Model problem 2}. Set
$$X=[H^1(\Omega)\cap L^2_{0}(\Omega)]^2,~W=[L^2(\Omega)]^2,~M=L^2(\Omega)~~{\rm and}~~Q=[L^2(\Gamma)]^2.$$
 The continuous admissible control set is $Q_{b}$, where $Q_{b}$ is defined in (\ref{qb}). The discrete spaces are
 $$X_{h}={\rm CR}^1(\cT_h)\cap [L^2_{0}(\Omega)]^2~~ {\rm and}~ M_{h}= \mathbb{P}_{0}(\mathcal T_h). $$
Define discrete control space as $$Q_{h}=\{\mathbf{x}_h\in [L^2(\Gamma)]^2 : \mathbf{x}_h|_{e}\in [\mathbb{P}_{0}(e)]^2 \;\;\text{for all}\;\; e\in \mathcal{E}^b_{h}\},$$
and the admissible discrete control set is $$Q_{ad}^h=\{\mathbf{x}_h\in Q_{h}:\int_{\Gamma}{\mathbf{x_h}}\ds+\int_{\Omega}{\mathbf{f}}\dx=\mathbf{0},\;\; \mathbf{y}_a\leq \mathbf{x}_h \leq \mathbf{y}_b\}.$$ It is clear that $Q_{ad}^h \subseteq Q_{ad}$ and $\Pi_{h}\mathbf{u}\in Q_{ad}^h$ for $\mathbf{u}\in Q_{ad}.$ The operator $E: X\rightarrow Q$ is the trace map  and $E_h:X_{h}\rightarrow Q_{h}$ is defined
by the piecewise (edge-wise) trace, i.e., $E_h
 \mathbf{w}|_e=\mathbf{w}_T|_e$ where $\mathbf{w}_T=\mathbf{w}|_T$ and $T$ be the triangle having the edge $e$ on boundary.
The bilinear forms are given by
\begin{equation}\label{58}
a_{h}(\mathbf{v}_{h},\mathbf{z}_{h})= \sum_{T\in \mathcal{T}_{h}} \int_{T} \nabla{\mathbf{v}_{h}}:\nabla{\mathbf{z}_{h}}\dx,~~~~~~~~b_{h}(\mathbf{z}_{h},p_{h})= -\sum_{T\in \mathcal{T}_{h}} \int_{T} p_{h}\nabla\cdot{\mathbf{z}_{h}}\dx,
\end{equation}
for all  $\mathbf{v}_{h}$,~ $\mathbf{z}_{h}\in X_{h}$ and $p_{h}\in M_{h}.$
The energy norm on $X_{h}$ is defined by $$\lVert \mathbf{v}_{h}\rVert_{h}^2=\sum_{T\in \mathcal{T}_{h}} \int_{T} \rvert\nabla{\mathbf{v}_{h}}\lvert^2 \dx.$$
 The inequality (\ref{Assmption1}) follows
 from the results Poincar\'e-Friedrichs  type inequalities in \cite{Brenner2003}. The estimate in (\ref{Assmption2})
 follows from the well-known trace inequality and \cite{Brenner2003}.

\begin{theorem}\label{thm4.13}
	Let $s\in (0,1]$ be the elliptic regularity index. Then, it holds
	\begin{align*}
	\lVert   \mathbf{y}-\mathbf{y}_h& \rVert_{0,\Gamma}+\lVert \bm{\phi}-\bm{\phi}_{h}\rVert_{h}+\lVert \mathbf{u}-\mathbf{u}_{h}\rVert_{h} \leq C \big( h^s \lVert \mathbf{u}\rVert_{1+s,\Omega}+h^s \lVert \bm{\phi}\rVert_{1+s,\Omega}+h^s \lVert p\rVert_{s,\Omega}+h^s \lVert r\rVert_{s,\Omega}\\&+\sum_{i=1}^{m}h^{1+s} \lVert \mathbf{y}\rVert_{\frac{1}{2}+s,\Gamma_{i}}+h\lVert \mathbf{f} \rVert_{0,\Omega}+h\lVert \mathbf{u}-\mathbf{u}_d \rVert_{0,\Omega}+ h\lVert \bm{\phi}\rVert_{1+s,\Omega}+\sum_{i=1}^{m} h^{\frac{1}{2}+s} \lVert \mathbf{y}\rVert_{\frac{1}{2}+s,\Gamma_{i}} \big).
	\end{align*}
	
\end{theorem}
\begin{proof}
	From Theorem \ref{thm 2.3} we have
	$$ \lVert \bm{\phi}-\bm{\phi}_{h}\rVert_{h}\leq C \big(\lVert \bm{\phi}-\bar{P_{h}}\bm{\phi}\rVert_{h}+\norm{\bm{\phi}-\Pi_{h}\bm{\phi}}_{0,\Gamma} +\lVert \mathbf{y}-\Pi_{h}\mathbf{y}\rVert_{0,\Gamma}+\lVert \mathbf{u}-P_{h}\mathbf{u} \rVert_{h}\big),$$
	and
	$$\lVert \mathbf{u}-\mathbf{u}_{h}\rVert_{h} \leq C \big(\lVert \mathbf{u}-P_{h}\mathbf{u} \rVert_{h}+\norm{\bm{\phi}-\Pi_{h}\bm{\phi}}_{0,\Gamma}+\lVert \mathbf{y}-\Pi_{h}\mathbf{y}\rVert_{0,\Gamma}+\lVert \bm{\phi}-\bar{P_{h}}\bm{\phi}\rVert_{h}\big).$$
	As a consequence of the error analysis in \cite[Theorem 3.1]{BCGG:2014:DGStokes} for Neumann boundary problem we have,
	\begin{equation}\label{61}
	\lVert\bm{\phi}-\bar{P_{h}}\bm{\phi}\rVert_{h} \leq C\big( \inf_{{\bf w}_{h}\in Z_{h}}\lVert\bm{\phi}-\mathbf{w}_{h}\rVert_{h}+\inf_{s_{h}\in M_{h}}\lVert r-s_{h}\rVert_{0,\Omega}+ h\lVert \mathbf{u}-\mathbf{u}_d\rVert_{0,\Omega}\big),
	\end{equation}
	\begin{equation}\label{62}
	\lVert \mathbf{u}-\bar{P_{h}}\mathbf{u}\rVert_{h} \leq C\big( \inf_{{\bf v}_{h}\in Z_{h}}\lVert \mathbf{u}-\mathbf{v}_{h}\rVert_{h}+\inf_{r_{h}\in M_{h}}\lVert p-r_{h}\rVert_{0,\Omega}+h\lVert \mathbf{f} \rVert_{0,\Omega}+ h^{\frac{1}{2}} \inf_{\mathbf{x}_h\in Q_{h}}\lVert \mathbf{y}-\mathbf{x}_h\rVert_{0,\Gamma}\big).
	\end{equation}
	From the Theorem \ref{thm2.2} we have
	$$\lVert   \mathbf{y}-\mathbf{y}_h \rVert_{0,\Gamma} \leq C\big(\norm{\bm{\phi}-\Pi_{h}\bm{\phi}}_{0,\Gamma} +\lVert \mathbf{y}-\Pi_{h}\mathbf{y}\rVert_{0,\Gamma}+\lVert \bm{\phi}-\bar{P_{h}}\bm{\phi}\rVert_{h}+\lVert \mathbf{u}-P_{h}\mathbf{u} \rVert_{h}\big).$$
	Using the above estimates, we find
	\begin{align}\label{see1}
	\lVert   \mathbf{y}-\mathbf{y}_h \rVert_{0,\Gamma}+\lVert \bm{\phi}-\bm{\phi}_{h}\rVert_{h}+& \lVert \mathbf{u}-\mathbf{u}_{h}\rVert_{h} \leq C \big(\inf_{\mathbf{v}_{h}\in Z_{h}}\lVert \mathbf{u}-\mathbf{v}_{h}\rVert_{h}+\inf_{\mathbf{w}_{h}\in Z_{h}}\lVert\bm{\phi}-\mathbf{w}_{h}\rVert_{h}\nonumber\\&+\inf_{s_{h}\in M_{h}}\lVert p-s_{h}\rVert_{0,\Omega} +\inf_{s_{h}\in M_{h}}\lVert r-s_{h}\rVert_{0,\Omega}+ h^{\frac{1}{2}} \inf_{\mathbf{x}_h\in Q_{h}}\lVert \mathbf{y}-\mathbf{x}_h\rVert_{0,\Gamma}\nonumber\\
	& +h\lVert \mathbf{f} \rVert_{0,\Omega}+h\lVert \mathbf{u}-\mathbf{u}_d\rVert_{0,\Omega} +\norm{\bm{\phi}-\Pi_{h}(\bm{\phi})}_{0,\Gamma} +\lVert \mathbf{y}-\Pi_{h}\mathbf{y}\rVert_{0,\Gamma} \big).
	\end{align}
	 The elliptic regularity on polygonal domains implies that $\bm{\phi} \in [H^{1+s}(\Omega)]^2\cap X$ and $r \in H^s(\Omega)$  for some $ s\in (0,1]$, which depends on the interior angles of the domain $\Omega.$ We know that,  $\mathbf{y}={\Pi_{[\mathbf{y}_a,\mathbf{y}_b]}}\big(-\frac{\bm{\phi}}{\lambda}\big)\;\; \text{on}\; \partial\Omega.$ Hence the control $\mathbf{y}|_{\Gamma_{i}} \in [H^{\frac{1}{2}+s}(\Gamma_{i})]^2$ for $1\leq i \leq m$ and $\mathbf{u}\in [H^{1+s}(\Omega)]^2\cap X$ and $p \in H^s(\Omega),$ where $m$ is the number of boundary edges.
	Also we have the following estimates,
	\begin{equation}\label{691}
	\norm{\bm{\phi}-\Pi_{h}(\bm{\phi})}_{0,\Gamma}\leq C h^{\frac{1}{2}+s} \lVert \bm{\phi}\rVert_{1+s,\Omega},~~~~\text{and}~~~~
	\norm{\mathbf{y}-\Pi_{h}(\mathbf{y})}_{0,\Gamma}\leq C\sum_{i=1}^{m} h^{\frac{1}{2}+s} \lVert \mathbf{y}\rVert_{\frac{1}{2}+s,\Gamma_{i}}.
	\end{equation}
	Using (\ref{65})-(\ref{69}) and (\ref{691}) in (\ref{see1}), we have proved the theorem.
\end{proof}

\begin{theorem}
	Let $s\in (0,1]$ be the elliptic regularity index. Then there holds
	\begin{align*}
	\lVert p-p_{h}\rVert_{0,\Omega} \leq& C \big( h^s \lVert \mathbf{u}\rVert_{1+s,\Omega}+h^s \lVert \bm{\phi}\rVert_{1+s,\Omega}+h^s \lVert p\rVert_{s,\Omega}+h^s \lVert r\rVert_{s,\Omega}+\sum_{i=1}^{m}h^{1+s} \lVert \mathbf{y}\rVert_{\frac{1}{2}+s,\Gamma_{i}
	}\\
	& +h\lVert \mathbf{f} \rVert_{0,\Omega}+h\lVert \mathbf{u}-\mathbf{u}_d \rVert_{0,\Omega}+ h\lVert \bm{\phi}\rVert_{1+s,\Omega}+\sum_{i=1}^{m} h^{\frac{1}{2}+s} \lVert \mathbf{y}\rVert_{\frac{1}{2}+s, \Gamma_{i}} \big).
	\end{align*}
	
\end{theorem}
\begin{proof}
	From Theorem \ref{thm 2.4}, we have the pressure estimate as below
	\begin{align*}
	\lVert p-p_{h}\rVert_{0,\Omega}\leq C\big(\lVert p-R_{h}p\rVert_{0,\Gamma}+& \lVert \mathbf{u}-P_{h}\mathbf{u} \rVert_{h}+\lVert \bm{\phi}-\bar{P_{h}}\bm{\phi}\rVert_{h}\\&+\norm{\bm{\phi}-\Pi_{h}\bm{\phi}}_{0,\Gamma} +\lVert \mathbf{y}-\Pi_{h}\mathbf{y}\rVert_{0,\Gamma}\big).
	\end{align*}
		As a consequence of the error analysis in \cite[Theorem 4.1]{BCGG:2014:DGStokes} for Neumann boundary value  problem, we have
	\begin{align*}
	\lVert p-R_{h}p\rVert_{0,\Omega} \leq & C \big(\inf_{\mathbf{v}_{h}\in Z_{h}}\lVert \mathbf{u}-\mathbf{v}_{h}\rVert_{h}+\inf_{q_{h}\in M_{h}}\lVert p-q_{h}\rVert_{0,\Omega}+\lVert p-P^{0}p\rVert_{0,\Omega}\\&+ h^{\frac{1}{2}} \inf_{\mathbf{x}_h\in Q_{h}}\lVert \mathbf{y}-\mathbf{x}_h\rVert_{0,\Gamma}+h\lVert \mathbf{f} \rVert_{0,\Omega}\big).
	\end{align*}
	Using the above estimates and (\ref{61})-(\ref{62}) we get
	\begin{align*}
	\lVert p-p_{h}\rVert_{0,\Omega}\leq& ~C \big(\inf_{\mathbf{v}_{h}\in Z_{h}}\lVert \mathbf{u}-\mathbf{v}_{h}\rVert_{h}+\inf_{\mathbf{w}_{h}\in Z_{h}}\lVert\bm{\phi}-\mathbf{w}_{h}\rVert_{h}+\inf_{s_{h}\in M_{h}}\lVert p-s_{h}\rVert_{0,\Omega}\\
	& +\inf_{s_{h}\in M_{h}}\lVert r-s_{h}\rVert_{0,\Omega}+h\lVert \mathbf{f} \rVert_{0,\Omega}+ h^{\frac{1}{2}} \inf_{\mathbf{x}_h\in Q_{h}}\lVert \mathbf{y}-\mathbf{x}_h\rVert_{0,\Gamma}
	+h\lVert \mathbf{u}-\mathbf{u}_d\rVert_{0,\Omega}\\
	&+\lVert p-P^{0}p\rVert_{0,\Omega}+\norm{\bm{\phi}-\Pi_{h}\bm{\phi}}_{0,\Gamma} +\lVert \mathbf{y}-\Pi_{h}\mathbf{y}\rVert_{0,\Gamma} \big)
	\end{align*}
Finally, use (\ref{65}), (\ref{66}) and (\ref{69}) in above estimate concludes the proof.
\end{proof}
\begin{theorem}\label{thm 4.15}
{\bf (A~posteriori error estimator) }
	It holds,
	\begin{align}
	\lVert \mathbf{y}-\mathbf{y}_h\rVert_{0,\Gamma}+\lVert \mathbf{u}-\mathbf{u}_{h}\rVert_{h} &+\lVert p-p_{h}\rVert_{0,\Omega}+\lVert \bm{\phi}-\bm{\phi}_{h}\rVert_{h}+\lVert r-r_{h}\rVert_{0,\Omega} \nonumber\\
	&	\leq C \big( \eta_{(\mathbf{u},p)}+\eta_{(\bm{\phi},r)}
	+\lVert\bm{\phi}_{h}-\Pi_{h}\bm{\phi}_{h}\rVert_{0,\Gamma}\big),\label{esti_2}
	\end{align}
	where
	the estimators are  defined as
	\begin{align*}
	\eta_{(\mathbf{u},p)}^2=\sum_{T\in \mathcal{T}_{h}} h_{T}^2 \lVert \mathbf{f}\rVert_{0,T}^2&+\sum_{e\in \mathcal{E}_{h}^{int}} \lVert h^{\frac{1}{2}}\sjump{p_{h}I-\nabla{\mathbf{u}_{h}}}\rVert_{0,e}^2\\
	&+\sum_{e\in \mathcal{E}_{h}^{b}} \lVert h^{\frac{1}{2}}(p_{h}I-\nabla{\mathbf{u}_{h}}+\mathbf{y}_hI)\rVert_{0,e}^2+\sum_{e\in \mathcal{E}^i_{h}} \lVert h^{-\frac{1}{2}}\underline{\sjump{\mathbf{u}_{h}}}\rVert_{0,e}^2,
	\end{align*}
	and
	\begin{align*}
	\eta_{(\bm{\phi},r)}^2=\sum_{T\in \mathcal{T}_{h}} h_{T}^2 \lVert \mathbf{u}_{h}-\mathbf{u}_d\rVert_{0,T}^2&+\sum_{e\in \mathcal{E}_{h}^{int}} \lVert h^{\frac{1}{2}}\sjump{r_{h}I+\nabla{\bm{\phi}_{h}}}\rVert_{0,e}^2\\
	&+\sum_{e\in \mathcal{E}_{h}^{b}} \lVert h^{\frac{1}{2}}\sjump{r_{h}I+\nabla{\bm{\phi}_{h}}}\rVert_{0,e}^2 +\sum_{e\in \mathcal{E}^i_{h}} \lVert h^{-\frac{1}{2}}\underline{\sjump{\bm{\phi}_{h}}}\rVert_{0,e}^2.
	\end{align*}
\end{theorem}
\begin{proof}
From Theorem \ref{thm2.7},~ \ref{thm2.8} and \ref{thms1}, we have
	\begin{align}
	\lVert \mathbf{y}-\mathbf{y}_h\rVert_{0,\Gamma}+&\lVert \mathbf{u}-\mathbf{u}_{h}\rVert_{h} +\lVert p-p_{h}\rVert_{0,\Omega}+\lVert \bm{\phi}-\bm{\phi}_{h}\rVert_{h}+\lVert r-r_{h}\rVert_{0,\Omega}\leq C \big(\lVert R\mathbf{u}-\mathbf{u}_{h}\rVert_{h}\nonumber\\
	&+\lVert R_{0}p-p_{h}\rVert_{0,\Omega}
	+\lVert\bm{\phi}_{h}-\bar{R}\bm{\phi}\rVert_{h}+\lVert \bar{R}_{0}r-r_{h}\rVert_{0,\Omega}+\lVert\bm{\phi}_{h}-\Pi_{h}\bm{\phi}_{h}\rVert_{0,\Gamma}\big).\label{4146}
	\end{align}
	As a consequence of the error analysis in \cite[\text{Theorem 3.1}]{HSW:2005:Err.Stokes} for Neumann boundary value problem, concludes the following error estimates:
	\begin{equation}\label{63}
	\lVert R\mathbf{u}-\mathbf{u}_{h} \rVert_{h}+\lVert R_{0}p-p_{h}\rVert_{0,\Omega}\leq C \eta_{(\mathbf{u},p)},
	\end{equation}
	\begin{equation}\label{64}
	\lVert\bar{R} \bm{\phi}-\bm{\phi}_{h} \rVert_{h}+\lVert \bar{R_{0}}r-r_{h}\rVert_{0,\Omega}\leq C \eta_{(\bm{\phi},r)}.
	\end{equation}
	The substitution of (\ref{63}) and (\ref{64}) in  (\ref{4146}) completes the proof of (\ref{esti_2}).
\end{proof}

\begin{theorem}\label{thm 4.16}{\bf (Efficiency)}
	Let $\mathcal{T}_{e}$ be the set of two triangles sharing the edge $e\in \mathcal{E}_{h}^{i}.$ Then, it hold
	\begin{align*}
	h_{T} \lVert \mathbf{f}\rVert_{0,T}&\leq C \big(\lVert \nabla(\mathbf{u}-\mathbf{u}_{h})\rVert_{0,T}+\lVert p-p_{h}\rVert_{0,T}+ osc(\mathbf{f},T)\big),\\
	h_{T} \lVert \mathbf{u}_{h}-\mathbf{u}_d\rVert_{0,T}&\leq C\big(\lVert \mathbf{u}-\mathbf{u}_{h}\rVert_{0,T}+\lVert \nabla( \bm{\phi}-\bm{\phi}_{h})\rVert_{0,T}+\lVert r-r_{h}\rVert_{0,T}+ osc(\mathbf{u}_d,T)\big),\\
	\lVert h^{\frac{1}{2}}\sjump{p_{h}I-\nabla{\mathbf{u}_{h}}}\rVert_{0,e}&\leq C \sum_{T\in \mathcal{T}_{e}}\big(\lVert \nabla(\mathbf{u}-\mathbf{u}_{h})\rVert_{0,T}+\lVert p-p_{h}\rVert_{0,T}+ osc(\mathbf{f},T)+\lVert \mathbf{y}-\mathbf{y}_h\rVert_{0,T}\big),\\
	\lVert h^{\frac{1}{2}}\sjump{r_{h}I+\nabla{\bm{\phi}_{h}}}\rVert_{0,e}&\leq C \sum_{T\in\mathcal{T}_{e}}\big(\lVert \mathbf{u}-\mathbf{u}_{h}\rVert_{0,T}+\lVert \nabla( \bm{\phi}-\bm{\phi}_{h})\rVert_{0,T}+\lVert r-r_{h}\rVert_{0,T}+ osc(\mathbf{u}_d,T)\big),\\
	\lVert\bm{\phi}_{h}-\Pi_{h}\bm{\phi}_{h}\rVert_{0,T}&\leq C \big(\lVert\bm{\phi}_{h}-\bm{\phi}\rVert_{0,T}+\lVert\bm{\phi}-\Pi_{h}\bm{\phi}_{h}\rVert_{0,T}\big).
	\end{align*}
	Further, for any boundary edge $e\in \mathcal{E}_{h}^{b}$, it hold
	\begin{align*}
	\lVert h^{\frac{1}{2}}(p_{h}I-\nabla{\mathbf{u}_{h}}+\mathbf{y}_hI)\rVert_{0,e}&\leq C\big(\lVert \nabla(\mathbf{u}-\mathbf{u}_{h})\rVert_{0,T}+\lVert p-p_{h}\rVert_{0,T}+ osc(\mathbf{f},T)+\lVert \mathbf{y}-\mathbf{y}_h\rVert_{0,T}\big),\\
	\lVert h^{\frac{1}{2}}(r_{h}I+\nabla{\bm{\phi}_{h}})\rVert_{0,e}&\leq C \big(\lVert \mathbf{u}-\mathbf{u}_{h}\rVert_{0,T}+\lVert \nabla( \bm{\phi}-\bm{\phi}_{h})\rVert_{0,T}+\lVert r-r_{h}\rVert_{0,T}+ osc(\mathbf{u}_d,T)\big).
	\end{align*}
\end{theorem}
\begin{proof}
The above theorem on local efficiency can be deduced by the standard bubble functions technique \cite{Verfurth 1995}.
\end{proof}

\subsection{Discontinuous Galerkin Method $\mathbb{P}_{1}(\cT_h)/ \mathbb{P}_{0}(\cT_h)$}
In this subsection, we will discuss about discrete  problem for  DG pair $\mathbb{P}_{1}/ \mathbb{P}_{0}$ of velocity and pressure, respectively.
\subsubsection{\bf Discrete distributed control problem:}
Set $$X=[H^1_{0}(\Omega)]^2,~W=[L^2(\Omega)]^2,~M=L^2_{0}(\Omega)~~{\rm and~~}  Q=[L^2(\Omega)]^2.$$ The set $Q_{ad}=Q_{d}$, where $Q_{d}$ is defined in (\ref{qb1}). The discrete spaces are defined by
$$X_{h}:=[\mathbb{P}_1(\mathcal T_h)]^2,~{\rm and~~}
M_{h}:=  \mathbb{P}_0(\mathcal T_h)\cap L^2_{0}(\Omega).$$
The admissible control set $Q_{ad}^h=\{\mathbf{x}_h\in [\mathbb{P}_0(\mathcal T_h)]^2: \mathbf{y}_a\leq \mathbf{x}_h\leq \mathbf{y}_b\}.$
 It is clear that $Q_{ad}^h \subseteq Q_{ad}$ and $\Pi_{h}\mathbf{u}\in Q_{ad}^h$ for $\mathbf{u}\in Q_{ad}.$ The operators $E: X\rightarrow Q$ and $E_h:X_{h}\rightarrow Q_{h}$ are inclusion maps.

The interior penalty DG bilinear form for the diffusion term is given by
\begin{align}\label{37}
a_{h}(\mathbf{v}_{h},\mathbf{z}_{h})=& \sum_{T\in \mathcal{T}_{h}} \int_{T} \nabla{\mathbf{v}_{h}}:\nabla{\mathbf{z}_{h}} \dx-\sum_{e\in \mathcal{E}_{h}}\int_{e}\smean{\nabla{\mathbf{v}_{h}}}: \underline{\sjump{\mathbf{z}_{h}}}\ds -\sum_{e\in \mathcal{E}_{h}}\int_{e}\smean{\nabla{\mathbf{z}_{h}}} : \underline{ \sjump{\mathbf{v}_{h}}}\ds\nonumber\\&+\sum_{e\in \mathcal{E}_{h}} \frac{\sigma}{h_{e}}\int_{e} \underline{\sjump{\mathbf{v}_{h}}}: \underline{\sjump{\mathbf{z}_{h}}}\ds,
\end{align}
for all  $\mathbf{v}_{h}$,~$\mathbf{z}_{h}\in X_{h}$ and $\sigma> 0$ be a real number. The DG bilinear form for the pressure term is given by
\begin{equation}\label{38}
b_{h}(\mathbf{z}_{h},p_{h})= -\sum_{T\in \mathcal{T}_{h}} \int_{T} p_{h}\nabla\cdot{\mathbf{z}_{h}}\dx+\sum_{e\in \mathcal{E}_{h}}\int_{e}\smean{p_{h}} \sjump{\mathbf{z}_{h}}\ds,
\end{equation}
for all $\mathbf{z}_{h}\in X_{h}$ and $p_{h}\in M_{h}.$ After an integration by parts on the right-hand side of (\ref{38}), we have
\begin{equation}\label{39}
b_{h}(\mathbf{z}_{h},p_{h})= \sum_{T\in \mathcal{T}_{h}} \int_{T} \mathbf{z}_{h}\cdot\nabla{p_{h}}\dx-\sum_{e\in \mathcal{E}_{h}^{i}}\int_{e}\smean{p_{h}} \sjump{\mathbf{z}_{h}}\ds,
\end{equation}
for all $\mathbf{z}_{h}\in X_{h}$ and $p_{h}\in M_{h}.$
We choose $\sigma>0$ large enough such that $a_{h}(\cdot,\cdot)$ is $X_{h}$-elliptic with respect to the norm $\lVert \cdot\rVert_{h}$ on $X_{h}$ which is given by
$$\lVert \mathbf{v}_{h}\rVert_{h}^2=\sum_{T\in \mathcal{T}_{h}} \int_{T} \rvert\nabla{\mathbf{v}_{h}}\lvert^2\dx+\sum_{e\in \mathcal{E}_{h}} \frac{1}{h_{e}}\int_{e} \underline{\sjump{\mathbf{v}_{h}}}^2\ds$$
and the fact that $b_{h}(\cdot,\cdot)$ satisfies the inf-sup condition can be found in \cite{Ern:2012:DGBook}. Assumptions (\ref{Assmption1}) and (\ref{Assmption2})
are the Poincar\'e type inequalities derived in \cite{Brenner2003}.
\begin{theorem}\label{thm 4.1}
	Let $s\in (0,1]$ be the elliptic regularity index. Then it holds
	\begin{align}\label{4141}
	\lVert   \mathbf{y}-\mathbf{y}_h \rVert_{0,\Omega}+&\lVert \bm{\phi}-\bm{\phi}_{h}\rVert_{h}+ \lVert \mathbf{u}-\mathbf{u}_{h}\rVert_{h} \leq C \big( h^s \lVert \mathbf{u}\rVert_{1+s,\Omega}+h^s \lVert \bm{\phi}\rVert_{1+s,\Omega}+h^s \lVert p\rVert_{s,\Omega}+h^s \lVert r\rVert_{s,\Omega}
	\nonumber\\&+ h\inf_{\mathbf{x}_h\in Q_{h}}\lVert \mathbf{y}-\mathbf{x}_h\rVert_{0,\Omega}
	+h\lVert \mathbf{f} \rVert_{0,\Omega}+h\lVert \mathbf{u}-\mathbf{u}_d \rVert_{0,\Omega}+ h\lVert \bm{\phi}\rVert_{1+s,\Omega}+ h \lVert \mathbf{y}\rVert_{1,\Omega} \big).
	\end{align}
\end{theorem}
\begin{proof}  
	The best approximation results \cite[Theorem 3.1]{BCGG:2014:DGStokes}, give
	\begin{equation}\label{40}
	\lVert\bm{\phi}-\bar{P_{h}}\bm{\phi}\rVert_{h} \leq C\big( \inf_{\mathbf{w}_{h}\in Z_{h}}\lVert\bm{\phi}-\mathbf{w}_{h}\rVert_{h}+\inf_{s_{h}\in M_{h}}\lVert r-s_{h}\rVert_{0,\Omega}+ h\lVert \mathbf{u}-\mathbf{u}_d\rVert_{0,\Omega}\big)
	\end{equation}
	and
	\begin{equation}\label{41}
	\lVert \mathbf{u}-\bar{P_{h}}\mathbf{u}\rVert_{h} \leq C\big( \inf_{\mathbf{v}_{h}\in Z_{h}}\lVert \mathbf{u}-\mathbf{v}_{h}\rVert_{h}+\inf_{s_{h}\in M_{h}}\lVert p-s_{h}\rVert_{0,\Omega}+h \lVert \mathbf{f} \rVert_{0,\Omega}+ h\inf_{\mathbf{x}_h\in Q_{h}}\lVert \mathbf{y}-\mathbf{x}_h\rVert_{0,\Omega}\big).
	\end{equation}
	From  Theorem \ref{thm2.2} we get
	$$\lVert   \mathbf{y}-\mathbf{y}_h \rVert_{0,\Omega} \leq C\big(\norm{\bm{\phi}-\Pi_{h}(\bm{\phi})}_{0,\Omega} +\lVert \mathbf{y}-\Pi_{h}\mathbf{y}\rVert_{0,\Omega}+\lVert \bm{\phi}-\bar{P_{h}}\bm{\phi}\rVert_{h}+\lVert \mathbf{u}-P_{h}\mathbf{u} \rVert_{h}\big).$$
A use of the above estimates in  (\ref{thm232}) from Theorem \ref{thm 2.3} results in
	\begin{align}
	\lVert   \mathbf{y}-\mathbf{y}_h \rVert_{0,\Omega}+&\lVert \bm{\phi}-\bm{\phi}_{h}\rVert_{h}+ \lVert \mathbf{u}-\mathbf{u}_{h}\rVert_{h} \leq C \big(\inf_{\mathbf{v}_{h}\in Z_{h}}\lVert \mathbf{u}-\mathbf{v}_{h}\rVert_{h}+\inf_{\mathbf{w}_{h}\in Z_{h}}\lVert\bm{\phi}-\mathbf{w}_{h}\rVert_{h}\nonumber\\&+\inf_{s_{h}\in M_{h}}\lVert p-s_{h}\rVert_{0,\Omega} +\inf_{s_{h}\in M_{h}}\lVert r-s_{h}\rVert_{0,\Omega}+h \lVert \mathbf{f} \rVert_{0,\Omega}+  h\inf_{\mathbf{x}_h\in Q_{h}}\lVert \mathbf{y}-\mathbf{x}_h\rVert_{0,\Omega}\nonumber\\
	& + h\lVert \mathbf{u}-\mathbf{u}_d\rVert_{0,\Omega} +\norm{\bm{\phi}-\Pi_{h}(\bm{\phi})}_{0,\Omega} +\lVert \mathbf{y}-\Pi_{h}\mathbf{y}\rVert_{0,\Omega} \big).\label{a2}
	\end{align}	
 Using (\ref{65}), (\ref{66}) and (\ref{69}) in (\ref{a2}), we have the result (\ref{4141}).
\end{proof}

\begin{theorem}\label{thm 4.2}
	Let $s\in (0,1]$ be the elliptic regularity index. Then there holds
	\begin{align*}
	\lVert p-p_{h}\rVert_{0,\Omega} \leq C& \big( h^s \lVert \mathbf{u}\rVert_{1+s,\Omega}+h^s \lVert \bm{\phi}\rVert_{1+s,\Omega}+h^s \lVert p\rVert_{s,\Omega}+h^s \lVert r\rVert_{s,\Omega}+h^2 \lVert \mathbf{y}\rVert_{1,\Omega}\\
	& +h\lVert \mathbf{f} \rVert_{0,\Omega}+h\lVert \mathbf{u}-\mathbf{u}_d \rVert_{0,\Omega}+ h\lVert \bm{\phi}\rVert_{1+s,\Omega}+ h \lVert \mathbf{y}\rVert_{1,\Omega} \big).
	\end{align*}
	
\end{theorem}
\begin{proof} 
%
The best approximation result from \cite[Theorem 4.1]{BCGG:2014:DGStokes} gives
\begin{align*}
\lVert p-R_{h}p\rVert_{0,\Omega} \leq C \big(\inf_{\mathbf{v}_{h}\in Z_{h}}\lVert \mathbf{u}-\mathbf{v}_{h}\rVert_{h}+\inf_{s_{h}\in M_{h}}\lVert p-s_{h}\rVert_{0,\Omega}+&\lVert p-P^{0}p\rVert_{0,\Omega}\\&+  h\inf_{\mathbf{x}_h\in Q_{h}}\lVert \mathbf{y}-\mathbf{x}_h\rVert_{0,\Omega} +h\lVert \mathbf{f} \rVert_{0,\Omega}\big).
\end{align*}
	Using this  in the pressure estimate
	from Theorem \ref{thm 2.4} and (\ref{40})-(\ref{41}) we have
	\begin{align*}
	\lVert p-p_{h}\rVert_{0,\Omega}\leq C& \big(\inf_{\mathbf{v}_{h}\in Z_{h}}\lVert \mathbf{u}-\mathbf{v}_{h}\rVert_{h}+\inf_{\mathbf{w}_{h}\in Z_{h}}\lVert\bm{\phi}-\mathbf{w}_{h}\rVert_{h}+\inf_{s_{h}\in M_{h}}\lVert p-s_{h}\rVert_{0,\Omega}\\
	& +\inf_{s_{h}\in M_{h}}\lVert r-s_{h}\rVert_{0,\Omega}+h \lVert \mathbf{f} \rVert_{0,\Omega}+ h\inf_{\mathbf{x}_h\in Q_{h}}\lVert \mathbf{y}-\mathbf{x}_h\rVert_{0,\Omega} + h\lVert \mathbf{u}-\mathbf{u}_d\rVert_{0,\Omega}\\
	&+\lVert p-P^{0}p\rVert_{0,\Omega}+\norm{\bm{\phi}-\Pi_{h}(\bm{\phi})}_{0,\Omega} +\lVert \mathbf{y}-\Pi_{h}\mathbf{y}\rVert_{0,\Omega} \big).
	\end{align*}
Use of the approximations (\ref{65}), (\ref{66}) and (\ref{69}) in  the above equation concludes the proof.
\end{proof}


\begin{theorem}\label{thm 4.3}{\bf (A~posteriori error estimator) }
	There holds,
	\begin{align*}
	\lVert \mathbf{y}-\mathbf{y}_h\rVert_{0,\Omega}+\lVert \mathbf{u}-\mathbf{u}_{h}\rVert_{h} +\lVert p-p_{h}\rVert_{0,\Omega}+\lVert \bm{\phi}-\bm{\phi}_{h}\rVert_{h}+&\lVert r-r_{h}\rVert_{0,\Omega}\leq C \big( \eta_{(\mathbf{u},p)}+\eta_{(\bm{\phi},r)}\\
	&+\lVert\bm{\phi}_{h}-\Pi_{h}\bm{\phi}_{h}\rVert_{0,\Omega}\big),
	\end{align*}
	where the estimators are defined by,
	\begin{align*}
	\eta_{(\mathbf{u},p)}^2=&\sum_{T\in \mathcal{T}_{h}} \big(h_{T}^2 \lVert \mathbf{f}+\mathbf{y}_h\rVert_{0,T}^2+\lVert\nabla\cdot{\mathbf{u}_{h}}\rVert_{0,T}^2\big)+\sum_{e\in \mathcal{E}_{h}^{int}} \lVert h^{\frac{1}{2}}\sjump{p_{h}I-\nabla{\mathbf{u}_{h}}}\rVert_{0,e}^2\\
	&+ \eta^2 \sum_{e\in \mathcal{E}_{h}} \lVert h^{-\frac{1}{2}}\underline{\sjump{\mathbf{u}_{h}}}\rVert_{0,e}^2,
	\end{align*}
	and
	\begin{align*}
	\eta_{(\bm{\phi},r)}^2=&\sum_{T\in \mathcal{T}_{h}}\big( h_{T}^2 \lVert \mathbf{u}_{h}-\mathbf{u}_d\rVert_{0,T}^2+\lVert\nabla\cdot{\bm{\phi}_{h}}\rVert_{0,T}^2\big)+\sum_{e\in \mathcal{E}_{h}^{i}} \lVert h^{\frac{1}{2}}\sjump{r_{h}I+\nabla{\bm{\phi}_{h}}}\rVert_{0,e}^2\\
	&+ \sigma^2 \sum_{e\in \mathcal{E}_{h}} \lVert h^{-\frac{1}{2}}\underline{\sjump{\bm{\phi}_{h}}}\rVert_{0,e}^2.
	\end{align*}
\end{theorem}

\begin{proof} Theorem \ref{thm2.7}, \ref{thm2.8}, and \ref{thms1} imply
	\begin{align}\label{too1}
	\lVert \mathbf{y}-\mathbf{y}_h\rVert_{0,\Omega}+&\lVert \mathbf{u}-\mathbf{u}_{h}\rVert_{h} +\lVert p-p_{h}\rVert_{0,\Omega}+\lVert \bm{\phi}-\bm{\phi}_{h}\rVert_{h}+\lVert r-r_{h}\rVert_{0,\Omega}\leq C \big(\lVert R\mathbf{u}-\mathbf{u}_{h}\rVert_{h}\nonumber\\&
	+\lVert R_{0}p-p_{h}\rVert_{0,\Omega}
	+\norm{\bm{\phi}_{h}-\bar{R}\bm{\phi}}_{h}+\lVert \bar{R}_{0}r-r_{h}\rVert_{0,\Omega}+\lVert\bm{\phi}_{h}-\Pi_{h}\bm{\phi}_{h}\rVert_{0,\Omega}\big).
	\end{align}
	Again, the error analysis in \cite[\text{Section 5}]{Carsten & Gudi & Jensen 2009}, \cite[Theorem 3.1]{HSW:2005:Err.Stokes} allows to conclude the following error estimates:
	\begin{equation}\label{42}
	\lVert R\mathbf{u}-\mathbf{u}_{h} \rVert_{h}+\lVert R_{0}p-p_{h}\rVert_{0,\Omega}\leq C \eta_{(\mathbf{u},p)},
	\end{equation}
	\begin{equation}\label{43}
	\lVert\bar{R} \bm{\phi}-\bm{\phi}_{h} \rVert_{h}+\lVert \bar{R_{0}}r-r_{h}\rVert_{0,\Omega}\leq C \eta_{(\bm{\phi},r)}.
	\end{equation}
Now using (\ref{42}) and (\ref{43}) in (\ref{too1}), we complete the proof of the theorem.
\end{proof}


The standard bubble function techniques can deduce the following theorem on local efficiency:
\begin{theorem}\label{thm 4.4} {(\bf Efficiency)}
	Let $\mathcal{T}_{e}$ be the set of two triangles sharing the edge $e\in \mathcal{E}_{h}^{i}.$ Then there holds
	\begin{align*}
	h_{T}\lVert \mathbf{f}+\mathbf{y}_h\rVert_{0,T} &\leq C\big(\lVert \nabla(\mathbf{u}-\mathbf{u}_{h})\rVert_{0,T}+\lVert p-p_{h}\rVert_{0,T}+ osc(\mathbf{f,T})+\lVert \mathbf{y}-\mathbf{y}_h\rVert_{0,T}\big),\\
	h_{T} \lVert \mathbf{u}_{h}-\mathbf{u}_d\rVert_{0,T}&\leq C\big(\lVert \mathbf{u}-\mathbf{u}_{h}\rVert_{0,T}+\lVert \nabla( \bm{\phi}-\bm{\phi}_{h})\rVert_{0,T}+\lVert r-r_{h}\rVert_{0,T}+ osc(\mathbf{u}_d,T)\big),\\
	\lVert h^{\frac{1}{2}}\sjump{p_{h}I-\nabla{\mathbf{u}_{h}}}\rVert_{0,e}&\leq C \sum_{T\in \mathcal{T}_{e}}\big(\lVert \nabla(\mathbf{u}-\mathbf{u}_{h})\rVert_{0,T}+\lVert p-p_{h}\rVert_{0,T}+ osc(\mathbf{f},T)+\lVert \mathbf{y}-\mathbf{y}_h\rVert_{0,T}\big),\\
	\lVert h^{\frac{1}{2}}\sjump{r_{h}I+\nabla{\bm{\phi}_{h}}}\rVert_{0,e}&\leq C \sum_{T\in\mathcal{T}_{e}}\big(\lVert \mathbf{u}-\mathbf{u}_{h}\rVert_{0,T}+\lVert \nabla( \bm{\phi}-\bm{\phi}_{h})\rVert_{0,T}+\lVert r-r_{h}\rVert_{0,T}+ osc(\mathbf{u}_d,T)\big),\\
	\lVert\nabla\cdot{\mathbf{u}_{h}}\rVert_{0,T}&\leq C\lVert \nabla(\mathbf{u}-\mathbf{u}_{h})\rVert_{0,T},\\
	\lVert\nabla\cdot{\bm{\phi}_{h}}\rVert_{0,T}&\leq C \lVert \nabla( \bm{\phi}-\bm{\phi}_{h})\rVert_{0,T},\\
	\lVert\bm{\phi}_{h}-\Pi_{h}\bm{\phi}_{h}\rVert_{0,T}&\leq C \big(\lVert\bm{\phi}_{h}-\bm{\phi}\rVert_{0,T}+\lVert\bm{\phi}-\Pi_{h}\bm{\phi}_{h}\rVert_{0,T}\big).
	\end{align*}
\end{theorem}
\vspace{0.2cm}
\subsubsection{\bf Discrete boundary control problem:}
The model problem in this section is the model problem 2 introduced in the section \ref{sec3}. Set
$$X=[H^1(\Omega)\cap L_0^2(\Omega)]^2,~ W=[L^2(\Omega)]^2, ~M=L^2(\Omega)~ {\rm and} ~Q=[L^2(\Gamma)]^2.$$ The set $Q_{ad}=Q_{b}$, where $Q_{b}$ is defined in section \ref{sec3}. The discrete spaces
$$X_{h}:= [L_{0}^2(\Omega)]^2\cap [\mathbb{P}_{1}(\cT)]^2,~{\rm and}~M_{h}:= \mathbb{P}_{0}(\cT).$$
Define discrete control space  $Q_{h}=\{\mathbf{x}_h\in [L^2(\Gamma)]^2 : \mathbf{x}_h|_{e}\in [\mathbb{P}_{0}(e)]^2\; \text{for all}\; e\in \mathcal{E}^b_{h}\},$
and the admissible control set $Q_{ad}^h=\{\mathbf{x}_h\in Q_{h}:\int_{\Gamma}{\mathbf{x_h}}\ds+\int_{\Omega}{\mathbf{f}}\dx=\mathbf{0},\;\; \mathbf{y}_a\leq \mathbf{x}_h\leq \mathbf{y}_b\}.$

It is clear that $Q_{ad}^h \subseteq Q_{ad}$ and $\Pi_{h}\mathbf{u}\in Q_{ad}^h$ for $\mathbf{u}\in Q_{ad}.$
The operator $E: X\rightarrow Q$ is the trace map  and $E_h:X_{h}\rightarrow Q_{h}$ is defined
by the piecewise (edge-wise) trace, i.e., $E_h
\mathbf{w}|_e=\mathbf{w}_T|_e$ where $\mathbf{w}_T=\mathbf{w}|_T$ and $T$ be the triangle having the edge $e$ on boundary.
The DG bilinear form for the diffusion term is given by
\begin{align}\label{44}
a_{h}(\mathbf{v}_{h},\mathbf{z}_{h})=& \sum_{T\in \mathcal{T}_{h}} \int_{T} \nabla{\mathbf{v}_{h}}:\nabla{\mathbf{z}_{h}}\dx-\sum_{e\in \mathcal{E}^{i}_{h}}\int_{e}\smean{\nabla{\mathbf{v}_{h}}}: \underline{\sjump{\mathbf{z}_{h}}} \ds-\sum_{e\in \mathcal{E}^{i}_{h}}\int_{e}\smean{\nabla{\mathbf{z}_{h}}} : \underline{ \sjump{\mathbf{v}_{h}}}\ds\nonumber\\&+\sum_{e\in \mathcal{E}^{i}_{h}} \frac{\sigma}{h_{e}}\int_{e} \underline{\sjump{\mathbf{v}_{h}}}: \underline{\sjump{\mathbf{z}_{h}}}\ds,
\end{align}
for all  $\mathbf{v}_{h}$,~$\mathbf{z}_{h}\in X_{h}$ and $\sigma> 0$ be a real number. The DG bilinear form for the pressure term is given by
\begin{equation}\label{45}
b_{h}(\mathbf{z}_{h},p_{h})= -\sum_{T\in \mathcal{T}_{h}} \int_{T} p_{h}\nabla\cdot{\mathbf{z}_{h}} \dx+\sum_{e\in \mathcal{E}^{i}_{h}}\int_{e}\smean{p_{h}}  \sjump{\mathbf{z}_{h}} \ds,
\end{equation}
for all $\mathbf{z}_{h}\in X_{h}$ and $p_{h}\in M_{h}.$ After integration by parts on the right-hand side of (\ref{45}) we have
\begin{equation}\label{46}
b_{h}(\mathbf{z}_{h},p_{h})= \sum_{T\in \mathcal{T}_{h}} \int_{T} \mathbf{z}_{h}\cdot\nabla{p_{h}}\dx-\sum_{e\in \mathcal{E}_{h}}\int_{e}\sjump{p_{h}} \smean{\mathbf{z}_{h}}\ds,
\end{equation}
for all $\mathbf{z}_{h}\in X_{h}$ and $p_{h}\in M_{h}.$
We choose $\sigma>0$ large enough such that $a_{h}(\cdot,\cdot)$ is $X_{h}$-elliptic with respect to the norm $\lVert \cdot \rVert_{h}$ on $X_{h}$ which is given by
$\lVert \mathbf{v}_{h}\rVert_{h}^2=\sum_{T\in \mathcal{T}_{h}} \int_{T} \rvert\nabla{\mathbf{v}_{h}}\lvert^2 \dx+\sum_{e\in \mathcal{E}^i_{h}} \frac{1}{h_{e}}\int_{e} \underline{\sjump{\mathbf{v}_{h}}}^2\ds$ and $b_{h}$ satisfies the inf-sup condition. The inequality (\ref{Assmption1}) follows
from the results Poincar\'e-Friedrichs  type inequalities in \cite{Brenner2003}. The estimate in (\ref{Assmption2})
follows from the well-known trace inequality and \cite{Brenner2003}.


\begin{theorem}\label{thm 4.5}
	Let $s\in (0,1]$ be the elliptic regularity index. Then there holds
	\begin{align}\label{4143}
	\lVert   \mathbf{y}-\mathbf{y}_h \rVert_{0,\Gamma}&+\lVert \bm{\phi}-\bm{\phi}_{h}\rVert_{h}+\lVert \mathbf{u}-\mathbf{u}_{h}\rVert_{h} \leq C \big( h^s \lVert \mathbf{u}\rVert_{1+s,\Omega}+h^s \lVert \bm{\phi}\rVert_{1+s,\Omega}+h^s \lVert p\rVert_{s,\Omega}+h^s \lVert r\rVert_{s,\Omega}\nonumber\\&+\sum_{i=1}^{m}h^{1+s} \lVert \mathbf{y}\rVert_{\frac{1}{2}+s,\Gamma_{i}}
	+h\lVert \mathbf{f} \rVert_{0,\Omega}+h\lVert \mathbf{u}-\mathbf{u}_d \rVert_{0,\Omega}+ h\lVert \bm{\phi}\rVert_{1+s,\Omega}+\sum_{i=1}^{m} h^{\frac{1}{2}+s} \lVert \mathbf{y}\rVert_{\frac{1}{2}+s,\Gamma_{i}} \big).
	\end{align}	
\end{theorem}

\begin{proof}
	From Theorem \ref{thm 2.3} we have
	$$ \lVert \bm{\phi}-\bm{\phi}_{h}\rVert_{h}\leq C \big(\lVert \bm{\phi}-\bar{P_{h}}\bm{\phi}\rVert_{h}+\norm{\bm{\phi}-\Pi_{h}(\bm{\phi})}_{0,\Gamma} +\lVert \mathbf{y}-\Pi_{h}\mathbf{y}\rVert_{0,\Gamma}+\lVert \mathbf{u}-P_{h}\mathbf{u} \rVert_{h}\big),$$
	and
	$$\lVert \mathbf{u}-\mathbf{u}_{h}\rVert_{h} \leq C \big(\lVert \mathbf{u}-P_{h}\mathbf{u} \rVert_{h}+\norm{\bm{\phi}-\Pi_{h}(\bm{\phi})}_{0,\Gamma}+\lVert \mathbf{y}-\Pi_{h}\mathbf{y}\rVert_{0,\Gamma}+\lVert \bm{\phi}-\bar{P_{h}}\bm{\phi}\rVert_{h}\big).$$
	From \cite[Theorem 3.1]{BCGG:2014:DGStokes}, we have
	\begin{equation}\label{48}
	\lVert\bm{\phi}-\bar{P_{h}}\bm{\phi}\rVert_{h} \leq C\big( \inf_{\mathbf{w}_{h}\in Z_{h}}\lVert\bm{\phi}-\mathbf{w}_{h}\rVert_{h}+\inf_{s_{h}\in M_{h}}\lVert r-s_{h}\rVert_{0,\Omega}+ h\lVert \mathbf{u}-\mathbf{u}_d\rVert_{0,\Omega}\big)
	\end{equation}
	and
	\begin{equation}\label{49}
	\lVert \mathbf{u}-\bar{P_{h}}\mathbf{u}\rVert_{h} \leq C\big( \inf_{\mathbf{v}_{h}\in Z_{h}}\lVert \mathbf{u}-\mathbf{v}_{h}\rVert_{h}+\inf_{s_{h}\in M_{h}}\lVert p-s_{h}\rVert_{0,\Omega}+h \lVert \mathbf{f} \rVert_{0,\Omega}+h^{\frac{1}{2}} \inf_{\mathbf{x}_h\in Q_{h}}\lVert \mathbf{y}-\mathbf{x}_h\rVert_{0,\Gamma}\big).
	\end{equation}
	From  Theorem \ref{thm2.2}, we have
	$$\lVert   \mathbf{y}-\mathbf{y}_h \rVert_{0,\Gamma} \leq C\big(\norm{\bm{\phi}-\Pi_{h}\bm{\phi}}_{0,\Gamma} +\lVert \mathbf{y}-\Pi_{h}\mathbf{y}\rVert_{0,\Gamma}+\lVert \bm{\phi}-\bar{P_{h}}\bm{\phi}\rVert_{h}+\lVert \mathbf{u}-P_{h}\mathbf{u} \rVert_{h}\big).$$
	The above estimates yield
	\begin{align}\label{4142}
	\lVert   \mathbf{y}-\mathbf{y}_h \rVert_{0,\Gamma}+\lVert \bm{\phi}-\bm{\phi}_{h}\rVert_{h}+& \lVert \mathbf{u}-\mathbf{u}_{h}\rVert_{h} \leq C \big(\inf_{\mathbf{v}_{h}\in Z_{h}}\lVert \mathbf{u}-\mathbf{v}_{h}\rVert_{h}+\inf_{\mathbf{w}_{h}\in Z_{h}}\lVert\bm{\phi}-\mathbf{w}_{h}\rVert_{h}\nonumber\\&+\inf_{s_{h}\in M_{h}}\lVert p-s_{h}\rVert_{0,\Omega} +\inf_{s_{h}\in M_{h}}\lVert r-s_{h}\rVert_{0,\Omega}+ h^{\frac{1}{2}} \inf_{\mathbf{x}_h\in Q_{h}}\lVert \mathbf{y}-\mathbf{x}_h\rVert_{0,\Gamma}\nonumber\\
	&+h\lVert \mathbf{f} \rVert_{0,\Omega} +h\lVert \mathbf{u}-\mathbf{u}_d\rVert_{0,\Omega} +\norm{\bm{\phi}-\Pi_{h}\bm{\phi}}_{0,\Gamma} +\lVert \mathbf{y}-\Pi_{h}\mathbf{y}\rVert_{0,\Gamma} \big).
	\end{align}
	Now we can apply elliptic regularity to derive concrete error estimates. Note that by well-posedness of the problem, $\mathbf{u},\bm{\phi}\in [H^1(\Omega) \cap L^2_{0}(\Omega)]^2$ and $p\in L^2(\Omega)$. The elliptic regularity of polygonal domains implies that $\bm{\phi} \in [H^{1+s}(\Omega)]^2\cap X$ and $r \in H^s(\Omega)$  for some $ s\in (0,1]$, which depends on the interior angles of the domain $\Omega.$ We know that,
	$\mathbf{y}={\Pi_{[\mathbf{y}_a,\mathbf{y}_b]}}\big(-\frac{\bm{\phi}}{\lambda}\big)\;\; \text{on}\; \partial\Omega.$
	Hence the control $\mathbf{y}|_{\Gamma_{i}} \in [H^{\frac{1}{2}+s}(\Gamma_{i})]^2$ for $1\leq i\leq m$ and $\mathbf{u}\in [H^{1+s}(\Omega)]^2\cap X$ and $p \in H^s(\Omega),$ where $m$ is the number of boundary edges.
	Also we have the following estimates,
	\begin{equation}\label{67}
	\norm{\bm{\phi}-\Pi_{h}\bm{\phi}}_{0,\Gamma}\leq C h^{\frac{1}{2}+s} \lVert \bm{\phi}\rVert_{1+s,\Omega},~~~~~~~~~\text{and}~~~~~~\norm{\mathbf{y}-\Pi_{h}\mathbf{y}}_{0,\Gamma}\leq C \sum_{i=1}^{m} h^{\frac{1}{2}+s} \lVert \mathbf{y}\rVert_{\frac{1}{2}+s,\Gamma_{i}}.
	\end{equation}
	\noindent
Finally, substitution  of the estimates from (\ref{65}), (\ref{66}) and (\ref{67}) in (\ref{4142}) leads to (\ref{4143}), and this concludes the proof.	
\end{proof}

\begin{theorem}\label{thm 4.6}
	Let $s\in (0,1]$ be the elliptic regularity index. Then there holds
	\begin{align*}
	\lVert p-p_{h}\rVert_{0,\Omega} \leq C & \big( h^s \lVert \mathbf{u}\rVert_{1+s,\Omega}+h^s \lVert \bm{\phi}\rVert_{1+s,\Omega}+h^s \lVert p\rVert_{s,\Omega}+h^s \lVert r\rVert_{s,\Omega}+\sum_{i=1}^{m}h^{1+s} \lVert \mathbf{y}\rVert_{\frac{1}{2}+s, \Gamma_{i}}\\
	& +h\lVert \mathbf{f} \rVert_{0,\Omega}+h\lVert \mathbf{u}-\mathbf{u}_d \rVert_{0,\Omega}+ h\lVert \bm{\phi}\rVert_{1+s,\Omega} \big).
	\end{align*}
	
\end{theorem}

\begin{proof}
	From Theorem \ref{thm 2.4}, we have pressure estimate
	\begin{align*}
	\lVert p-p_{h}\rVert_{0,\Omega}\leq &~ C\big(\lVert p-R_{h}p\rVert_{0,\Omega}+ \lVert \mathbf{u}-P_{h}\mathbf{u} \rVert_{h}+\lVert \bm{\phi}-\bar{P_{h}}\bm{\phi}\rVert_{h}\\&+\norm{\bm{\phi}-\Pi_{h}(\bm{\phi})}_{0,\Gamma} +\lVert \mathbf{y}-\Pi_{h}\mathbf{y}\rVert_{0,\Gamma}\big).
	\end{align*}
	As a consequence of the error analysis in \cite[Theorem 4.1]{BCGG:2014:DGStokes} for Neumann boundary value  problem, we have
	\begin{align*}
	\lVert p-R_{h}p\rVert_{0,\Omega} \leq &~C \Big(\inf_{\mathbf{v}_{h}\in Z_{h}}\lVert \mathbf{u}-\mathbf{v}_{h}\rVert_{h}+\lVert p-P^{0}p\rVert_{0,\Omega}+\inf_{q_{h}\in M_{h}}\lVert p-q_{h}\rVert_{0,\Omega}\\&+ h^{\frac{1}{2}} \inf_{\mathbf{x}_h\in Q_{h}}\lVert \mathbf{y}-\mathbf{x}_h\rVert_{0,\Gamma}+h\lVert \mathbf{f} \rVert_{0,\Omega}\Big).
	\end{align*}
	Using the (\ref{48}), (\ref{49}) and the above estimate we get
	\begin{align*}
	\lVert p-p_{h}\rVert_{0,\Omega}\leq C & ~\Big(\inf_{\mathbf{v}_{h}\in Z_{h}}\lVert \mathbf{u}-\mathbf{v}_{h}\rVert_{h}+\inf_{\mathbf{w}_{h}\in Z_{h}}\lVert\bm{\phi}-\mathbf{w}_{h}\rVert_{h}+\inf_{s_{h}\in M_{h}}\lVert p-s_{h}\rVert_{0,\Omega}\\
	&+\inf_{s_{h}\in M_{h}}\lVert r-s_{h}\rVert_{0,\Omega}+h\lVert \mathbf{f} \rVert_{0,\Omega}+h^{\frac{1}{2}} \inf_{\mathbf{x}_h\in Q_{h}}\lVert \mathbf{y}-\mathbf{x}_h\rVert_{0,\Gamma}
	+h\lVert \mathbf{u}-\mathbf{u}_d\rVert_{0,\Omega}\\
	&+\lVert p-P^{0}p\rVert_{0,\Omega}+\norm{\bm{\phi}-\Pi_{h}\bm{\phi}}_{0,\Omega} +\lVert \mathbf{y}-\Pi_{h}\mathbf{y}\rVert_{0,\Gamma} \Big).
	\end{align*}
	By using (\ref{65}),  (\ref{66}) and (\ref{67}) we have the pressure estimate.
\end{proof}

\begin{theorem}\label{thm 4.7} {(\bf A posteriori error estimator)}
	There holds,
	\begin{align}\label{too3}
	\lVert \mathbf{y}-\mathbf{y}_h\rVert_{0,\Gamma}+\lVert \mathbf{u}-\mathbf{u}_{h}\rVert_{h} +\lVert p-p_{h}\rVert_{0,\Omega}+\lVert \bm{\phi}-\bm{\phi}_{h}\rVert_{h}+&\lVert r-r_{h}\rVert_{0,\Omega}\leq C \big( \eta_{(\mathbf{u},p)}+\eta_{(\bm{\phi},r)}\nonumber\\
	&+\lVert\bm{\phi}_{h}-\Pi_{h}\bm{\phi}_{h}\rVert_{0,\Gamma}\big).
	\end{align}
	where the estimators are defined by,
	\begin{align*}
	\eta_{(\mathbf{u},p)}^2=&\sum_{T\in \mathcal{T}_{h}}\big( h_{T}^2 \lVert \mathbf{f}\rVert_{0,T}^2+\lVert\nabla\cdot{\mathbf{u}_{h}}\rVert_{0,T}^2 \big)+\sum_{e\in \mathcal{E}_{h}^{i}} \lVert h^{\frac{1}{2}}\sjump{p_{h}I-\nabla{\mathbf{u}_{h}}}\rVert_{0,e}^2\\
	&+\sum_{e\in \mathcal{E}_{h}^{b}} \lVert h^{\frac{1}{2}}(p_{h}I-\nabla{\mathbf{u}_{h}}+\mathbf{y}_hI)\rVert_{0,e}^2+\sigma^2 \sum_{e\in \mathcal{E}_{h}^{i}} \lVert h^{-\frac{1}{2}}\underline{\sjump{\mathbf{u}_{h}}}\rVert_{0,e}^2,
	\end{align*}
	and
	\begin{align*}
	\eta_{(\bm{\phi},r)}^2=&\sum_{T\in \mathcal{T}_{h}} \big(h_{T}^2 \lVert \mathbf{u}_{h}-\mathbf{u}_d\rVert_{0,T}^2+\lVert\nabla\cdot{\bm{\phi}_{h}}\rVert_{0,T}^2)+\sum_{e\in \mathcal{E}_{h}^{i}} \lVert h^{\frac{1}{2}}\sjump{r_{h}I+\nabla{\bm{\phi}_{h}}}\rVert_{0,e}^2\\
	&+\sum_{e\in \mathcal{E}_{h}^{b}} \lVert h^{\frac{1}{2}}(r_{h}I+\nabla{\bm{\phi}_{h}})\rVert_{0,e}^2+\sigma^2 \sum_{e\in \mathcal{E}_{h}^{i}} \lVert h^{-\frac{1}{2}}\underline{\sjump{\bm{\phi}_{h}}}\rVert_{0,e}^2.
	\end{align*}
\end{theorem}
\begin{proof}From Theorem \ref{thm2.7}, \ref{thm2.8} and \ref{thms1}, we get
	\begin{align}\label{too2}
	\lVert \mathbf{y}-\mathbf{y}_h\rVert_{0,\Gamma}+&\lVert \mathbf{u}-\mathbf{u}_{h}\rVert_{h} +\lVert p-p_{h}\rVert_{0,\Omega}+\lVert \bm{\phi}-\bm{\phi}_{h}\rVert_{h}+\lVert r-r_{h}\rVert_{0,\Omega}\leq C \big(\lVert R\mathbf{u}-\mathbf{u}_{h}\rVert_{h}\nonumber\\&
	+\lVert R_{0}p-p_{h}\rVert_{0,\Omega}
	+\norm{\bm{\phi}_{h}-\bar{R}\bm{\phi}}_{h}+\lVert R_{0}r-r_{h}\rVert_{0,\Omega}+\lVert\bm{\phi}_{h}-\Pi_{h}\bm{\phi}_{h}\rVert_{0,\Gamma}\big).
	\end{align}
	As a consequence of the error analysis in \cite[\text{Section 5}]{Carsten & Gudi & Jensen 2009},~\cite[Theorem 3.1]{HSW:2005:Err.Stokes} for Neumann boundary value problem conclude the following error estimates:
	 \begin{equation}\label{50}
	\lVert R\mathbf{u}-\mathbf{u}_{h} \rVert_{h}+\lVert R_{0}p-p_{h}\rVert_{0,\Omega}\leq C \eta_{(\mathbf{u},p)},
	\end{equation}
	\begin{equation}\label{51}
	\lVert\bar{R} \bm{\phi}-\bm{\phi}_{h} \rVert_{h}+\lVert \bar{R_{0}}r-r_{h}\rVert_{0,\Omega}\leq C \eta_{(\bm{\phi},r)}.
	\end{equation}
	Now, the substitution of (\ref{50}) and (\ref{51})  in (\ref{too2}) concludes the   proof of (\ref{too3}).
\end{proof}
\begin{theorem}\label{thm 4.8} {\bf(Efficiency)}
	Let $\mathcal{T}_{e}$ be the set of two triangles sharing the edge $e\in \mathcal{E}_{h}^{i}.$ Then there hold
	\begin{align*}
	h_{T} \lVert \mathbf{f}\rVert_{0,T}&\leq C \big(\lVert \nabla(\mathbf{u}-\mathbf{u}_{h})\rVert_{0,T}+\lVert p-p_{h}\rVert_{0,T}+ osc(\mathbf{f},T)\big),\\
	h_{T} \lVert \mathbf{u}_{h}-\mathbf{u}_d\rVert_{0,T}&\leq C\big(\lVert \mathbf{u}-\mathbf{u}_{h}\rVert_{0,T}+\lVert \nabla( \bm{\phi}-\bm{\phi}_{h})\rVert_{0,T}+\lVert r-r_{h}\rVert_{0,T}+ osc(\mathbf{u}_d,T)\big),\\
	\lVert h^{\frac{1}{2}}\sjump{p_{h}I-\nabla{\mathbf{u}_{h}}}\rVert_{0,e}&\leq C \sum_{T\in \mathcal{T}_{e}}\big(\lVert \nabla(\mathbf{u}-\mathbf{u}_{h})\rVert_{0,T}+\lVert p-p_{h}\rVert_{0,T}+ osc(\mathbf{f},T)+\lVert \mathbf{y}-\mathbf{y}_h\rVert_{0,T}\big),\\
	\lVert\nabla\cdot{\mathbf{u}_{h}}\rVert_{0,T}&\leq C\lVert \nabla(\mathbf{u}-\mathbf{u}_{h})\rVert_{0,T},\\
	\lVert\nabla\cdot{\bm{\phi}_{h}}\rVert_{0,T}&\leq C \lVert \nabla( \bm{\phi}-\bm{\phi}_{h})\rVert_{0,T},\\
	\lVert\bm{\phi}_{h}-\Pi_{h}\bm{\phi}_{h}\rVert_{0,T}&\leq C \big(\lVert\bm{\phi}_{h}-\bm{\phi}\rVert_{0,T}+\lVert\bm{\phi}-\Pi_{h}\bm{\phi}_{h}\rVert_{0,T}\big).
	\end{align*}
	Further for any boundary edge $e\in \mathcal{E}_{h}^{b}$, there hold
	\begin{align*}
	\lVert h^{\frac{1}{2}}(p_{h}I-\nabla{\mathbf{u}_{h}}+\mathbf{y}_hI)\rVert_{0,e}&\leq C\big(\lVert \nabla(\mathbf{u}-\mathbf{u}_{h})\rVert_{0,T}+\lVert p-p_{h}\rVert_{0,T}+ osc(\mathbf{f},T)+\lVert \mathbf{y}-\mathbf{y}_h\rVert_{0,T}\big),\\
	\lVert h^{\frac{1}{2}}(r_{h}I+\nabla{\bm{\phi}_{h}})\rVert_{0,e}&\leq C \big(\lVert \mathbf{u}-\mathbf{u}_{h}\rVert_{0,T}+\lVert \nabla( \bm{\phi}-\bm{\phi}_{h})\rVert_{0,T}+\lVert r-r_{h}\rVert_{0,T}+ osc(\mathbf{u}_d,T)\big).\\
	\end{align*}
\end{theorem}
\begin{remark}
	The analysis can be extended to the three dimensions also for the simplicity we strict ourselves to two dimensions.
\end{remark}

\section{Numerical Experiments}\label{sec5}
\begin{figure}[h!]
	\centering
		\includegraphics[width=18cm,height=2.4in]{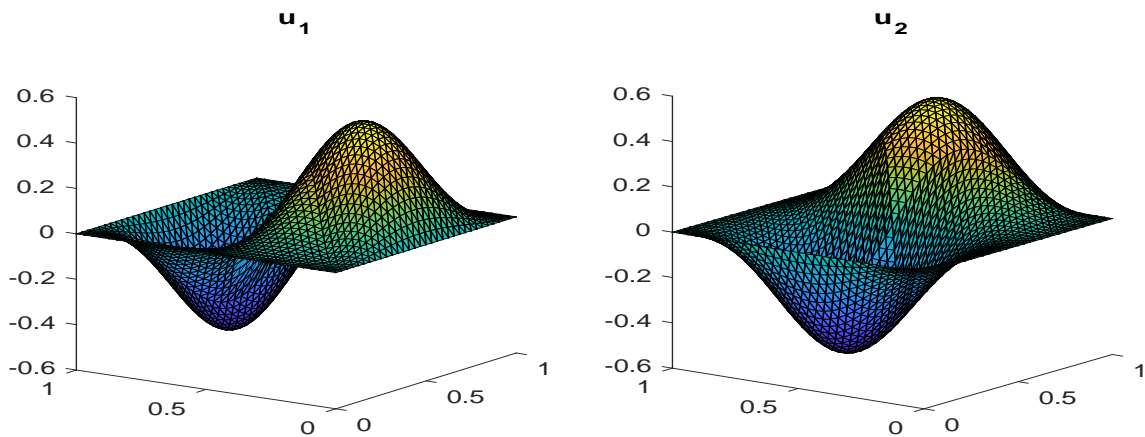}
				\caption{Discrete velocity ${\bf u}_h=(u_1, u_2)$ for Example \ref{ex1} }
		\label{fig6.1}
		\end{figure}

\begin{figure}[h!]
		\centering
	\includegraphics[width=18cm,height=2.9in]{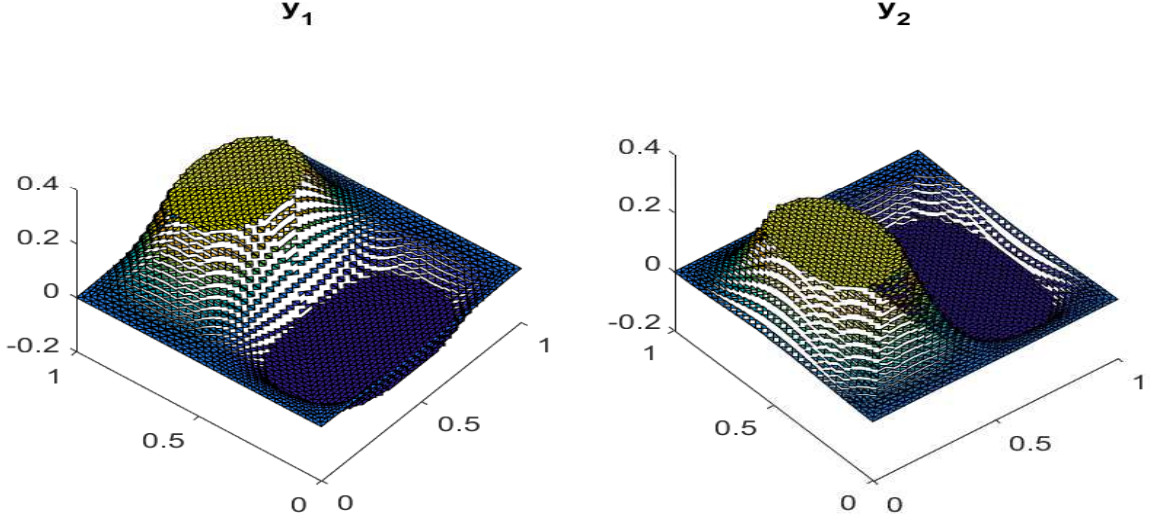}
		\caption{Discrete control ${\bf y}_h=(y_1,y_2)$ of Example \ref{ex1}}
			\label{fig6.2}
	\end{figure}

This section  presents some numerical experiments to illustrate the theoretical results derived in the article. The abstract framework of {\it a~priori} and {\it a~posteriori} error analysis is applicable for the set of discrete spaces $X_h \times M_h\times Q_{ad}^h$ for the approximation of velocity,  pressure, and control in conforming, nonconforming  FEM and discontinuous Galerkin methods as discussed in Section \ref{sec4}. Here in the following the numerical experiments, we have considered   ${\rm CR}_{0}^1(\mathcal T_h)\times \big(\mathbb P_0(\mathcal T_h)\cap L_0^2(\Omega)\big) \times \big(\mathbb P_0(\mathcal T_h)\cap Q_d\big) $ spaces for the approximations.
\begin{example}\label{ex1}Consider the optimal control problem (\ref{mp1})-(\ref{mp3})  with the domain $\Omega=(0, 1)^2$ and the exact solution
\begin{align}
{\bf u}= \bm{\phi}= \left(
\begin{array}{c}
\sin^2(\pi x) \sin(\pi y) \cos(\pi y) \\
-\sin^2(\pi y) \sin(\pi x) \cos(\pi x)
\end{array}
\right),~~~ p=r = \sin(2\pi x)\sin(2 \pi y),
\end{align}
and ${\bf y}=\Pi_{[a,b]} \big( -\frac{1}{\lambda} {\bf u}\big)$, where $\mathbf{y}_a=-0.1,~\mathbf{y}_b=0.25$.  The data of the problem are  chosen such that
\begin{align}
\mathbf{f} = -\Delta \mathbf{u}+\nabla{p}-\mathbf{y},~{\rm and}~\mathbf{u}_d= \mathbf{u}+\Delta \bm{\phi}+\nabla{r}.
\end{align}
\end{example}
In this numerical simulation, the  ${\rm CR}_{0}^1$/$\mathbb P_0$ pair is used for the approximations of state and adjoint state velocity and pressure variables, and piecewise constant space  for the control variable. For the computation of the discrete solution, the primal-dual algorithm \cite[pp. 100]{Troltzzsch:2005:book} is used. The discrete approximations for the state velocity variables using nonconforming finite elements
	are shown in Figure \ref{fig6.1} and the discrete approximations for the control variable ${\bf y}$ using piecewise constant elements
	are shown in Figure \ref{fig6.2}.

Table \ref{tabla12} displays the errors and convergence rates of FE approximations. The linear  convergence is observed for error in approximation of state and adjoint state velocity  in energy norm, and also for state pressure, adjoint pressure and control variables in $L_2$-norm. Moreover, we have also observed the quadratic convergence in $L^2$-norm for state and adjoint state velocity variables.
\begin{table}[h!]
    \begin{center}
    \footnotesize
\begin{tabular}{ c||c|c|c|c|c|c|c|c|c|c }
    $h$ &$\norm{{\bf u}-{\bf u}_h}_h$ & CR&  $\norm{p-p_{h}}$& CR&  $\norm{\bm{\phi}-\bm{\phi}_{h}}_h$  & CR& $\norm{r-r_{h}}$ & CR& $\norm{{\bf y}-{\bf y}_h}$ &CR \\
    \hline\hline
           0.2500 &   0.8877&      0 &    0.9362  &     0  &  0.8888  &     0  &  0.9360   &   0   & 0.0985 &        0 \\
           0.1250 &   0.5350&    0.73&    0.3511  &  1.41 &   0.5351  &  0.73 &   0.3510  &  1.41   & 0.0556  &  0.82\\
           0.0625 &   0.2680 &   0.99 &   0.1682  &  1.06 &   0.2680 &   0.99 &   0.1682  &  1.06  &  0.0296  &  0.91\\
           0.0312 &   0.1346 &   0.99 &   0.0819  &  1.03 &   0.1346 &   0.99 &   0.0819  &  1.03  &  0.0152  &  0.96\\
           0.0156 &   0.0674 &   0.99 &   0.0406  &  1.01 &   0.0674 &   0.99 &   0.0406  &  1.01  &  0.0076 &   0.98\\
            0.0078 &   0.0337  & 0.99&    0.0202  &  1.00 &   0.0337&    0.99  &  0.0202   & 1.00 &   0.0038   & 0.99\\
    \hline
\end{tabular}

\vspace{0.2cm}
\caption{Errors and  convergence rates (CR) for the Example \ref{ex1}.}
\label{tabla12}
\end{center}
\end{table}

\begin{example}\label{ex2}  Consider the optimal control problem (\ref{mp1})-(\ref{mp3})  with the  L-shaped domain $\Omega=(-1, 1)^2\setminus ((0,1)\times (-1,0))$ and the exact solution
    \begin{align*}
    {\bf u}= r^\alpha \left(
    \begin{array}{c}
     (1+\alpha)\sin(\theta)\omega(\theta)+\cos(\theta)\omega'(\theta)\\
-(1+\alpha)\cos(\theta)\omega(\theta)+\sin(\theta)\omega'(\theta)
    \end{array}
    \right),\\
    p=- r^{\alpha-1}((1+\alpha)^2 \omega'(\theta)+\omega{'''}(\theta))/(1-\alpha),
    \end{align*}
    where
    \begin{align*}
        \omega(\theta)=&1/(1+\alpha)\sin(\alpha+1)\theta)\cos(\alpha w)-\cos((\alpha+1)\theta)\\ &+1/(1+\alpha)\sin(\alpha-1)\theta)\cos(\alpha\omega)-\cos((\alpha-1)\theta)
    \end{align*}
and $\alpha=856399/1572864$ and $w=3\pi/2$. The adjoint variables ${\rm \phi},~r$ are considered as same as in Example \ref{ex1}  and  ${\bf y}=\Pi_{[a,b]} \big( -\frac{1}{\lambda} {\bf \phi}\big)$, where $\mathbf{y}_a=-0.1,~\mathbf{y}_b=0.25$. The data of the problem is  chosen such that
    \begin{align}
    \mathbf{f} = -\Delta \mathbf{u}+\nabla{p}-\mathbf{y},~{\rm and}~\mathbf{u}_d= \mathbf{u}+\Delta \bm{\phi}+\nabla{r}.
    \end{align}
\end{example}
 This problem is defined on the L-shaped domain, and the solution $({\bf u},p)$ has a singularity at the origin. It is known that for this problem the uniform refinements will not provide an optimal convergence rate. We have similar observation from Figure \ref{fig_sol2}, for uniform refinements convergence rate with respect to the number of degrees of freedom (Ndof) is 0.25 (that is, respect to the mesh-size $h\equiv {\rm Ndof}^{-1/2} $). Hence, we have to use the adaptive algorithm to get the optimal convergence. The adaptive algorithm contains a loop:~ {\it
 Solve $\rightarrow$ Estimate $\rightarrow$ Mark $\rightarrow$ Refine}.
\begin{figure}[h!]
	\centering
	\centering
	\includegraphics[width=18cm,height=2.9in]{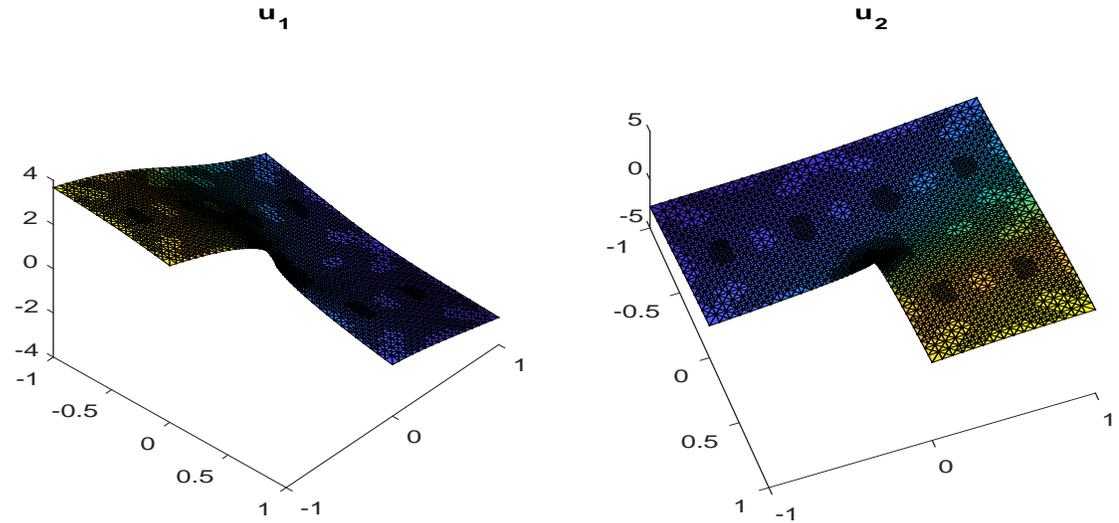}
	
	\caption{Discrete solution ${\bf u}_h=(u_1,u_2)$ of Example \ref{ex2}}
\label{6.3}
\end{figure}
\begin{figure}[h!]
	\includegraphics[height=2.5in]{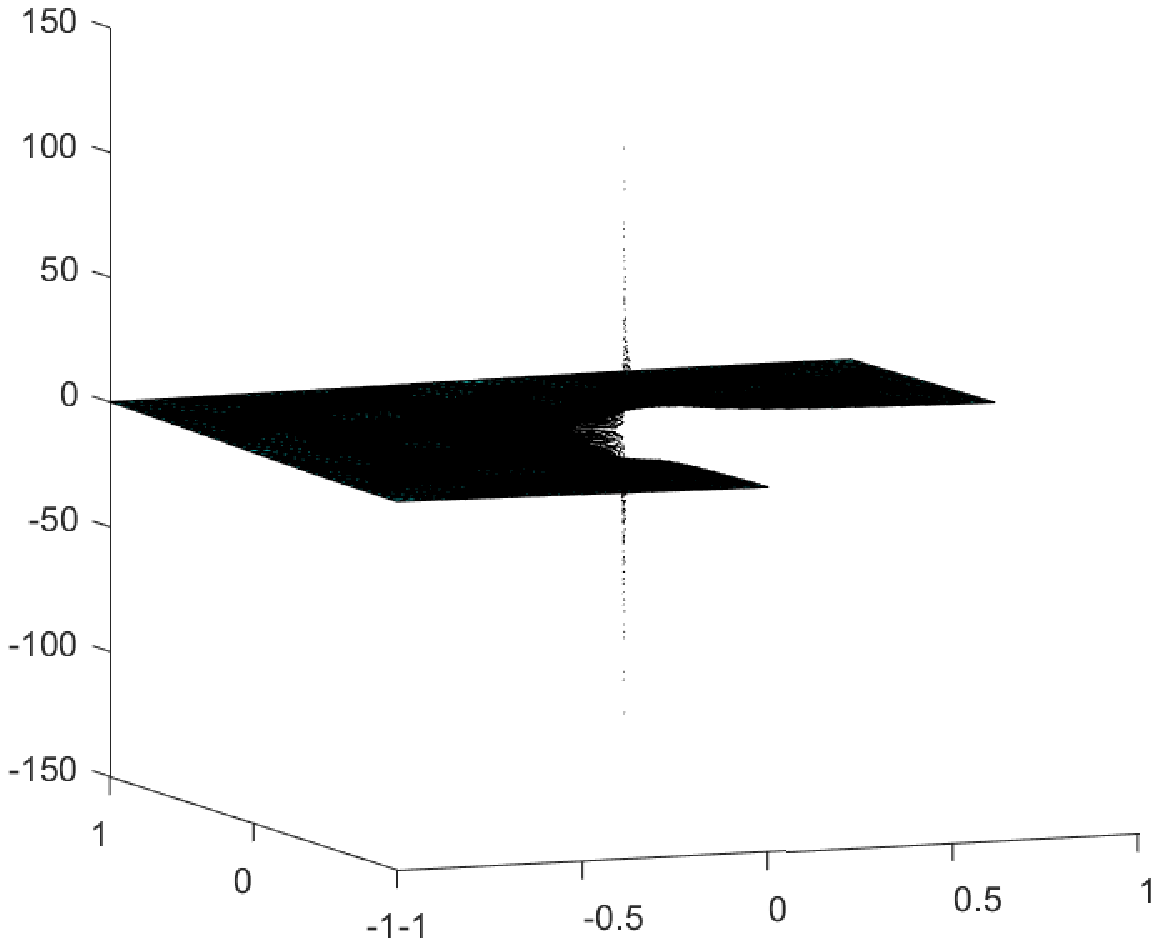}
		\includegraphics[height=2.5in]{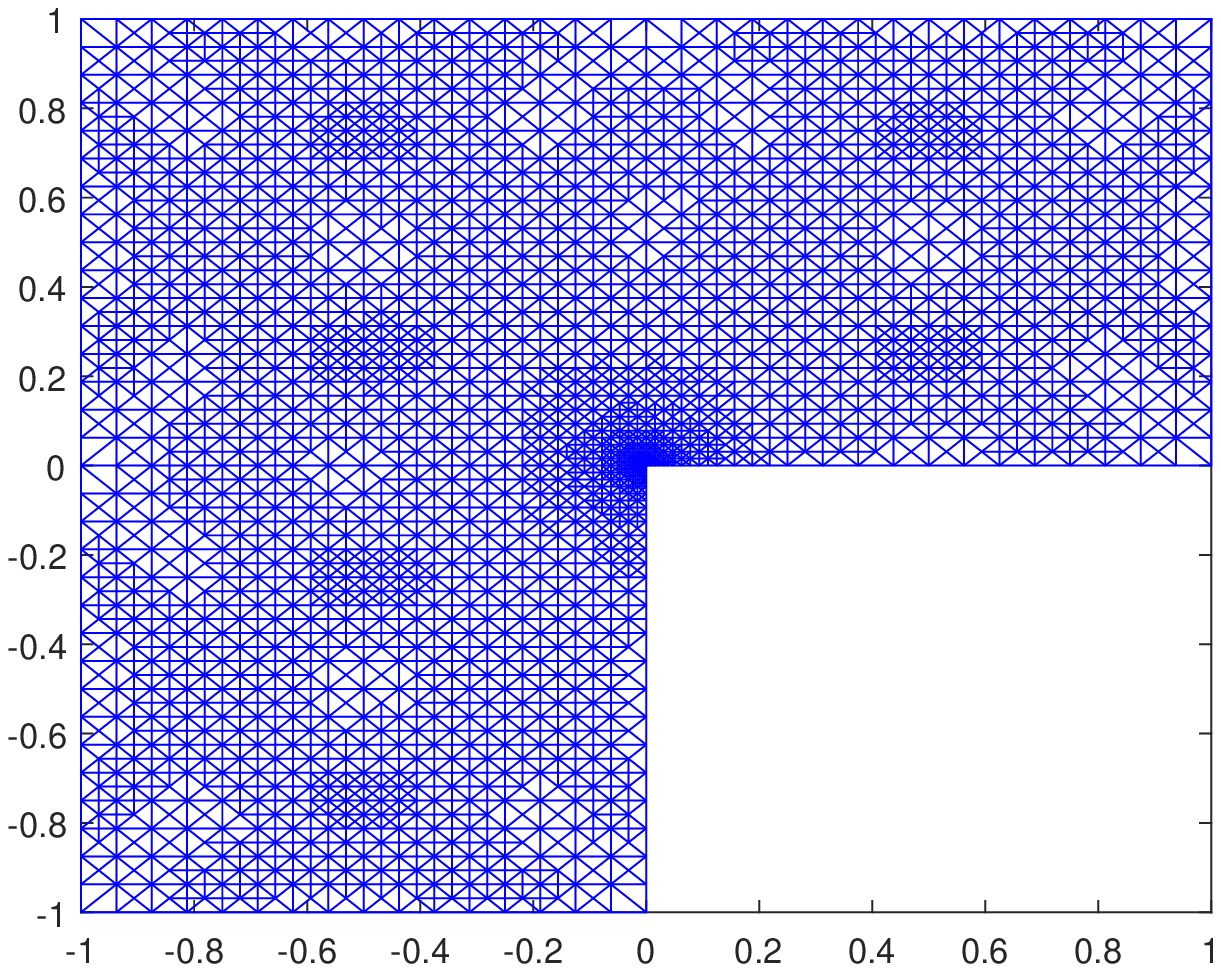}
	\caption{Discrete pressure $p_h$ and adaptive mesh-refinement of Example \ref{ex2}}
\label{fig6.4}
\end{figure}

\noindent
First, we compute the discrete solutions using the primal-dual algorithm.  Then, in the second step using the discrete solution we compute  the error
estimator $\eta$ (as defined in Theorem (\ref{thm 4.11})) over each element.  We use  the D\"{o}rlfer marking technique
{\cite{dofler}} with bulk parameter $\theta=0.3$ for the mark step and
the newest vertex bisection algorithm  for  mesh-refinements. Figure {\ref{6.3}} displays the discrete approximation to velocity ${\bf u}=(u_1,u_2)$ and the left-hand side image from Figure \ref{fig6.4} show the discrete approximation to the pressure $p$. The right-hand image of Figure  \ref{fig6.4} shows adaptive mesh generated from  several iterations of adaptive refinements. In figure \ref{fig6.4}, due to a singularity at the origin of the state velocity and pressure variables, more mesh-refinements are observed at the origin, while other  refinements are results of the adjoint variable estimator.

\begin{figure}[h!]
	   \centering
	 \includegraphics[]{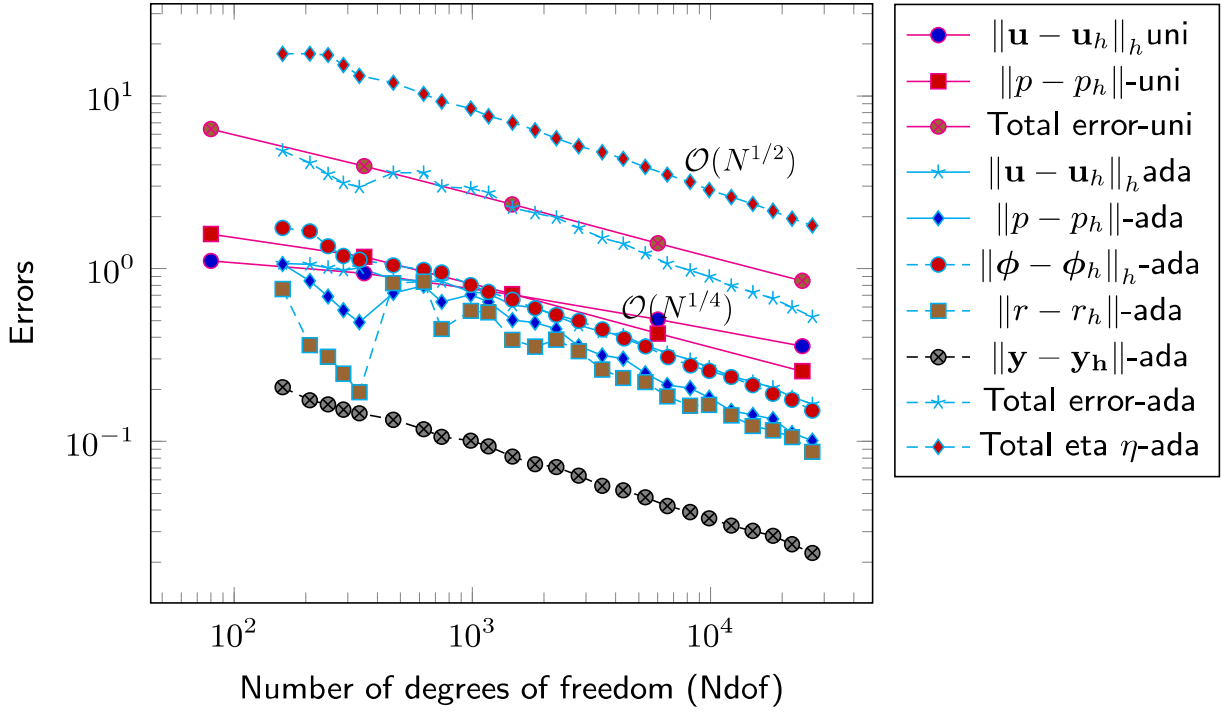}
 \caption{Convergence history with uniform and adaptive refinements for  Example \ref{ex2}}
   \label{fig_sol2}
\end{figure}
\noindent
Figure \ref{fig_sol2} depicts convergence rates for errors and estimators  with uniform and adaptive refinements. The optimal convergence  is achieved using the adaptive algorithm for the error in energy norm in the state and adjoint state velocity approximation, in $L^2$-norm  of control, pressure and adjoint pressure variables. Hence, the optimal convergence  for the  {a~posteriori} estimator and the total error which is the combination of all errors term as in the left-hand side of (\ref{terror}).
Here, the optimal convergence means 
order 0.5 with respect to Ndof.


\section*{Conclusions}\label{sec6}
In this paper, we have developed an abstract framework for discontinuous finite element methods  error analysis of both the distributed control and Neumann boundary control problems governed by the stationary Stokes equation, with control constraints. This framework will also work for linear elliptic and mixed optimal control problems with control constraints. The abstract analysis provides the best approximation results, which will be useful in the convergence analysis of adaptive methods and delivers a reliable and efficient a posteriori error estimators. Numerical experiments illustrate the theoretical findings. The results in the article will not directly cover the analysis of nonlinear mixed elliptic optimal control problems;  however, they will be useful to analyze the nonlinear problems.

\section*{Acknowledgment} The first author gratefully acknowledges financial support from the National Board for
Higher Mathematics (NBHM), Government of India. 

\end{document}